\RequirePackage{fix-cm}
\documentclass[smallcondensed]{svjour3-ours}     
\smartqed  
\usepackage{graphicx}
 \usepackage{mathptmx}      

\usepackage{amsmath, amsfonts, amssymb, amscd, enumerate, enumitem, color, cancel, soul}
\let\lpolish\l
\newtheorem{theo}{Theorem}[section]

\renewcommand\thesection{\arabic{section}.}
\renewcommand\thesubsection{\arabic{section}.\arabic{subsection}}

\newtheorem{lem}[theo]{Lemma}
\newtheorem{sublem}[theo]{Sublemma}
\newtheorem{prop}[theo]{Proposition}

\newtheorem{cor}[theo]{Corollary}
\newtheorem{rem}[theo]{Remark}
\newtheorem{definit}[theo]{Definition}
\newenvironment{pf}{\noindent{\it Proof. }}{$\square$\par\medskip}
\newenvironment{pfns}{\noindent{\it Proof. }}{\par\medskip}

\renewcommand{\=}{:=}
\newcommand{\beq}{\begin{equation}}
\newcommand{\eeq}{\end{equation}}

\renewcommand{\a}{\alpha}
\renewcommand{\b}{\beta}

\renewcommand{\d}{\delta}

\newcommand{\f}{\varphi}
\newcommand{\g}{\gamma}
\newcommand{\h}{\eta}

\renewcommand{\k}{\kappa}
\renewcommand{\l}{\lambda}
\renewcommand{\o}{\omega}
\newcommand{\q}{\vartheta}
\renewcommand{\r}{\rho}
\newcommand{\s}{\sigma}
\renewcommand{\t}{\tau}

\newcommand{\z}{\zeta}
\newcommand{\D}{\Delta}

\renewcommand{\O}{\Omega}

\font\sixrm=cmr6

\renewcommand{\sc}{\sixrm}

\newcommand{\bB}{\mathbb{B}}
\newcommand{\bC}{\mathbb{C}}

\newcommand{\bE}{\mathbb{E}}
\newcommand{\bF}{\mathbb{F}}

\newcommand{\bR}{\mathbb{R}}

\newcommand{\bJ}{\mathbb{J}}

\newcommand{\bN}{\mathbb{N}}

\newcommand{\bX}{\mathbb{X}}
\newcommand{\bY}{\mathbb{Y}}
\newcommand{\bW}{\mathbb{W}}

\newcommand{\gv}{\hbox{\bf v}}

\renewcommand\sl{\mathfrak{sl}}

\newcommand{\cA}{\mathcal{A}}

\newcommand{\cC}{\mathcal{C}}
\newcommand{\cD}{\mathcal{D}}

\newcommand{\cF}{\mathcal{F}}

\newcommand{\cH}{\mathcal{H}}
\newcommand{\cI}{\mathcal{I}}

\newcommand{\cL}{\mathcal{L}}

\newcommand{\cS}{\mathcal{S}}

\newcommand{\cU}{\mathcal{U}}
\newcommand{\cV}{\mathcal{V}}
\newcommand{\cW}{\mathcal{W}}

\newcommand{\cZ}{\mathcal{Z}}

\newcommand{\Jst}{J_{\operatorname{st}}}

\newcommand{\p}{\partial}

\DeclareMathOperator\Aut{Aut\;}
\DeclareMathOperator\Id{Id}
\newcommand\Hom{\operatorname{Hom}}

\newcommand{\wt}{\widetilde}
\newcommand{\wh}{\widehat}

\newcommand{\OB}{\overline \bB{}^n}
\newcommand{\TOB}{\wt{\overline \bB}{}^n}

\renewcommand{\Re}{\operatorname{Re}}
\renewcommand{\Im}{\operatorname{Im}}
\def\sideremark#1{\ifvmode\leavevmode\fi\vadjust{
\vbox to0pt{\hbox to 0pt{\hskip\hsize\hskip1em
\vbox{\hsize3cm\tiny\raggedright\pretolerance10000
\noindent #1\hfill}\hss}\vbox to8pt{\vfil}\vss}}}

\newcounter{ssig}
\newcounter{ttig}
\newcommand{\ve}{\varepsilon}

\begin{document}

\title{Propagation of regularity\\ for Monge-Amp\`ere  exhaustions and\\ Kobayashi metrics
}

\titlerunning{Propagation of regularity}        

\author{Giorgio Patrizio         \and
        Andrea Spiro 
}


\institute{G. Patrizio \at
             Dip. Matematica e Informatica ``U. Dini'',
Universit\`a di Firenze,
\& Istituto Nazionale di
Alta Matematica
``Francesco Severi''.
              \email{patrizio@math.unifi.it}           
           \and
           A. Spiro \at
             Scuola di Scienze e Tecnologie,
Universit\`a di Camerino,
Via Madonna delle Carceri,
I-62032 Camerino (Ma\-cerata).
\email{andrea.spiro@unicam.it}
}


\date{
}

\maketitle

\begin{abstract}
We prove that if a smoothly bounded strongly pseudoconvex domain $D \subset \bC^n$, $n \geq 2$,  admits at least one  Monge-Amp\`ere exhaustion  smooth up to the boundary  (i.e.  a plurisubharmonic  exhaustion $\t: \overline D \to [0,1]$,  which is $\cC^\infty$ at all points  except possibly  at the unique  minimum point $x$ and   with  $u \= \log \t$  
satisfying  the homogeneous complex Monge-Amp\`ere equation), then there exists a  bounded open neighborhood  $\cU\subset D$ of the minimum point $x$,  such that for each $y \in \cU$ there exists a Monge-Amp\`ere exhaustion with minimum at  $y$.     
This yields that for each such domain  $D$, the restriction  to the subdomain $\cU\subset D$ of the   Kobayashi  pseudo-metric $\k_D$  is  a  smooth Finsler metric for $\cU$ and 
each  pluricomplex Green function of $D$ with  pole at a point $y \in \cU$ is of class $\cC^\infty$.  The boundary of the maximal open subset having all such properties is also explicitly characterized.\par
The result is a direct consequence  of a  general theorem on abstract complex manifolds with boundary,  with Monge-Amp\`ere exhaustions 
of regularity $\cC^{k}$ for some $k \geq 5$. In fact, analogues of the above  properties  hold for each bounded strongly pseudoconvex complete circular  domain with boundary of such weaker regularity.  
\end{abstract}

\section{Introduction}
Motivated by  
Lempert's results  on smoothly bounded convex domains and the  several following developments,  the second author introduced in \cite{Sp} the  notion of {\it Lempert manifold}, that is    a complex manifold $(M,J)$   such that: 
\begin{itemize}
\item[i)] the   Kobayashi pseudo-metric $\k$  of $M$ is 
a smooth strongly pseudoconvex  Finsler metric, i.e. 
\begin{itemize}
\item[a)] It is a strictly positive and 
$\cC^\infty$  on $TM^o=TM\setminus\{Zero\;section\}$; 
\item[b)] 
 $\k(\lambda v) = |\lambda|\k(v)$ for any $\lambda \in \bC^*$ and $v\in TM^o$;
\item[c)]  The set of tangent vectors
$\cI_x = \{ v\in T_xM: \k(v) < 1\}$
is a strongly pseudoconvex  domain of $\bC^n \simeq T_xM$ for each $x \in M$;
\end{itemize}
\item[ii)]  all exponentials $\exp_x: T_x M \to M$, $x \in M$, are  diffeomorphisms and 
the  distance function $d$,  determined by $\k$,  is complete.
\item[iii)] for each    $w\in T^\bC_x M^o$, the $2$-dimensional real submanifold  $\exp(\bC w)$ of $M$ 
is a complex curve, which is  totally geodesic  and with induced metric of constant holomorphic curvature  $-4$. 
\end{itemize}
This definition was designed  to capture  most of the properties of  the smoothly  bounded strictly     convex domains.  Indeed, by \cite{Le,Le1,Pa,AP,AP1},   {\it each such domain  is a Lempert manifold}.  It is also quite  remarkable  that,   in contrast with the fact that on each smoothly bounded  strictly  convex domain (and each domain biholomorphic to one of them) the Kobayashi pseudo-metric $\k$ has  so many  and so nice properties,  up to now  no other  explicit class  of  domains,   for which  $\k$ is smooth on the whole $TM^o$, is known.  \par
\medskip
A closely related notion introduced by the first author in \cite{Pt} is the class of {\it manifolds of circular type} which  naturally includes  all  smoothly bounded  strictly  convex domains and all
smoothly bounded  strongly  pseudoconvex  circular domains.
Manifolds of circular type are characterized by  the existence of  a {\it Monge-Amp\`ere exhaustion}, i.e.  a plurisubharmonic  exhaustion $\t: M \to [0,1)$,  which is smooth at all points  except possibly  at the minimal set $  \{\t = 0\}$ and   with  $u \= \log \t$  
satisfying  the homogeneous complex Monge-Amp\`ere equation.  The minimal set always consists of only one point $x_o$, called {\it center}, and the 
 function  $\k=\k^{(\t)}: T_{x_o} M \setminus \{0\} \to [0, + \infty)$ for $(x_o, v)\in T_{x_o} M$, defined  by 
\beq \label{kob-tau} 
\kappa(x_o, v) \=\underset{{t_o\to 0}} \lim \left.\frac{d }{dt} \sqrt{\tau(\gamma(t))}\right|_{t = t_o},
\eeq
 where $\g$ is any smooth curve with  $\g(0) = x_o,\ \dot \g(0) = v$,
coincides with the  Kobayashi 
  metric of $M$  at the point $x_o$.\par
 \smallskip
 Apparently strictly  convex  domains and strongly pseudoconvex circular domains are   manifolds of circular type for unrelated reasons. Indeed, for a  smoothly bounded strictly  convex domain $\O \subset \bC^n$, for {\it any}  $x_o \in \O$, the exhaustion  
$\t$ defined by $\t(y)= (\tanh \delta(y,x_o))^{2}$, 
  with $\d(y, x_o)$  the  Kobayashi distance from $y$ to $x_o$,  is  a  Monge-Amp\`ere exhaustion  centered at $x_o$. This is a consequence of the very special properties of the Kobayashi metric and distance and their peculiar relation with the pluricomplex Green potential, discovered by Lempert for this class of domains. 
On the other hand, a smoothly bounded, strongly pseudoconvex complete circular domain  $D \subset \bC^n$ is  of circular type because  its  squared Minkowski function $\mu^2_D$ is  a Monge-Amp\`ere exhaustion centered at $0$. This is a simple consequence of the   symmetry properties of such domains, namely of the $S^{1}$ action that defines  them and  has the origin as the unique fixed point. While it is known that  the pluricomplex Green potential exists at every point and it is   in general only $C^{1,1}$, there is no simple reason why such a domain should admit a Monge-Amp\`ere exhaustion as defined above,  centered at some other point but the origin. Indeed, this would imply strong regularity and the least possible degeneracy of the Levi form of the 
pluricomplex Green potential. \par
\smallskip
Important steps towards a deeper understanding of the geometry of the set of centers of domains of circular type are provided by  
the results of  Y.-M. Pang (\cite{Pa}). Indeed Pang detected a  special condition on  holomorphic stationary disks,  called {\it coercitivity property}, which is  invariant under biholomorphisms (in contrast with the  strict convexity condition for  domains), it characterizes isolated local extremal disks and, most importantly, {\it it  is stable under small deformations}.   As an important  application  of such results and especially of the quoted stability  of coercitivity property, Pang  managed to prove that  {\it any smoothly  bounded strongly pseudoconvex complete circular domain $D$ 
admits  an open subset $\cU \subset D$ with the property that each $y \in \cU$ is the center of an appropriate smooth Monge-Amp\`ere exhaustion}. \par
\medskip
 In this  paper,  we show that the existence of such   ``clouds'' of centers  is  actually a  general property, provided that the considered domain  admits at least one Monge-Amp\`ere exhaustion $\t$, which  is smoothly extendible up to the boundary.  More precisely, we prove the following\par
\begin{theo} \label{main0} Let $ D \subset \bC^n$ be a  smoothly bounded strongly pseudoconvex domain,  admitting   a  Monge-Amp\`ere exhaustion $\t: D \to [0,1)$,  centered at $x_o \in D$, which  is smoothly  extendible up  to the boundary. Then
there exists  an open neighborhood $\cU\subset D$ of $x_o$, such that for each  $x \in \cU$ there exists a  Monge-Amp\`ere exhaustion $\t^{(x)}: M \to [0,1)$,  smoothly extendible up to the boundary and centered at $x$. 
The family of exhaustions $\t^{(x)}$, considered as collections of maps parameterized by the uniquely associated centers   $x\in \cU$, depends smoothly  on such parameter. \par
Further, in the cases in which  the  maximal connected neighborhood  $\cU^{(\text{\rm max})}$ of $x_o$  with the above property  is strictly  smaller than $D$, the  intersection   $\p\cU^{(\text{\rm max})} \cap f(\D)$   with each   stationary disk $f(\D) \subset D$  through $x_o$  coincides with the  set of  points,   where an appropriate matrix valued function  is singular (see Proposition \ref{prop5.4} below for details)
\end{theo}
In this result, the ``maximality'' of $\cU^{(\text{\rm max})}$  means that for 
each point  $y \in \p \cU^{(\text{\rm max})} \subset \overline D$ there is no  Monge-Amp\`ere exhaustion centered at $y$  that is smooth up to the boundary. Actually, as mentioned in the statement,  our proof   gives also a constructive method to determine the maximal open set $\cU^{(\text{\rm max})}$ and, in particular, its boundary.  It thus provides a biholomorphically invariant    characterization of the domains for which $D =\cU^{(\text{\rm max})}$, i.e. for the domains for which each point is the center of a smooth Monge-Amp\`ere exhaustion, hence of domains, on which the  Kobayashi metric is a smooth complex  Finsler metric. \par
In Theorem \ref{main0} the smoothness properties are deliberately not specified. Indeed whenever the boundary of the domain $D$ and the Monge-Amp\`ere exhaustion $\tau$ are  either $\cC^{k}$, for some $k\geq 5$,  or  $\cC^{\infty}$ or  $\cC^{\omega}$, then our  result yields that the exhaustions $\t^{(x)}$, $x \in \cU^{(\text{\rm max})} $, are  at least  $\cC^{k-3,\alpha}$, $\a \in (0,1)$,  in the first case, or $\cC^{\infty}$ or $\cC^{\omega}$ respectively  in the other cases, up to the boundary. \par
\smallskip
The above theorem    has the   noticeable consequence  that  {\it for  any smooth domain of circular type $D \subset \bC^n$ admitting at least one    Monge-Amp\`ere exhaustion that is smooth up to the boundary},   {\it the Kobayashi pseudo-metric $\k$  is    a   strongly pseudoconvex Finsler  metric on a (possibly non-connected) open subdomain of $D$}.\par 
 \par  \medskip
  Of course, all this can be rephrased in terms of the  pluricomplex Green functions and, in particular, 
  gives new information on their  regularity.  In fact,
if $\t^{(x)}: D \to [0, 1)$ is a  Monge-Amp\`ere exhaustion, with center  $x$ and smooth up to the boundary,   the   associated function
 \beq \label{primabis}  u^{(x)} = \log \t^{(x)}: \overline D \longrightarrow [-\infty, 0]\eeq
 is the pluricomplex Green function for $\overline D$ with pole at  $x$, it is   of class $C^{\infty}$ and satisfies the  homogeneous complex Monge-Amp\`ere equation on $D\setminus \{x\}$ with the least possible degeneracy, i.e. the annihilator of the form $\partial \bar \partial u^{(x)}$ has rank $1$ at every point of  $D\setminus \{x\}$.
 Conversely, if 
the pluricomplex Green function $u^{(x)}$  with pole at $x$ 
 is of class $C^{\infty}$ and satisfies the homogeneous complex  Monge-Amp\`ere equation on $D\setminus \{x\}$ with the least possible degeneracy, then $\t^{(x)} = \exp({u^{(x)}})$ is a Monge-Amp\`ere exhaustion centered at $x$.   \par 
  By general results on hyperconvex manifolds
 (\cite{De,Gu,Bo}), it is known that  for each point $x$ of a smoothly bounded strongly pseudoconvex domain $D \subset \bC^n$ there exists   a unique pluricomplex Green function $u^{(x)}$ with pole at $x$. Its  regularity is at least  $\cC^{1,1}$ up to the boundary but, in general,  of not higher regularity.  Our  result implies  the following curious phenomenon of propagation of regularity: {\it if a smoothly bounded strongly pseudoconvex domain $D \subset \bC^n$  admits at least one pluricomplex Green function $u^{(x)}$  of class $\cC^\infty$, then necessarily  {\rm all}  other pluricomplex Green functions $u^{(x')}$ with poles at the points  $x'$ in an appropriate  subdomain $\cU^{\text{\rm max}}\subset D$,  are of class $\cC^\infty$.}\par 
 \smallskip
While, by   Lempert's results  (\cite{Le}), it is known that on a  smoothly bounded {\it strictly convex} domain   $D \subset \bC^n$ the pluricomplex Green function $u^{(x)}$ with pole at $x$ is of class  $\cC^\infty$  for each $x\in D$, there are {\it no explicit examples} of smoothly bounded strongly pseudoconvex domain $D \subset \bC^n$  with   a  point $x_o$ for which  pluricomplex Green function $u^{(x_o)}$ is  $\cC^\infty$ and    such that the  surrounding  ``cloud'' $\cU^{\text{(max)}}\subset D$, made of the points $x$  for which 
 also  their pluricomplex Green functions $u^{(x)}$ are   $\cC^\infty$,  is strictly smaller than in $D$. On the other hand,  in Proposition 5.3, we provide specific  (closed) necessary conditions  for the  ``cloud''  $ \cU^{\text{(max)}}$ to be  properly contained in $D$. At the end of the paper we discuss how such conditions can be used to  construct explicit examples with $ \cU^{\text{(max)}} \subsetneq D$  and to determine  biholomorphically invariant sufficient conditions  for   propagations of regularity to the whole domain to occur. 
\par  
Theorem \ref{main0}  is   obtained as a consequence   of  a more general  result (Theorem \ref{main1})   concerning abstract   complex manifolds with boundary, equipped with   Monge-Amp\`ere exhaustions of  regularity $\cC^{k}$, $k \geq 5$.  Such result  shows that    the above  phenomena of propagation of regularity    hold also for the $\cC^{k}$-analogues of   domains of circular type,  in particular for  strongly pseudoconvex circular domains with boundaries of such regularity.    We remark that the  assumption  $k \geq 5$   is rather technical and determined by  the  special approach we exploit,  namely the use of  the so-called {\it manifolds in normal form}.  It is reasonable to expect that such lower bound for $k$ is not sharp and that  the  propagation phenomena for Monge-Amp\`ere exhaustions and pluricomplex Green functions should  occur  in full generality also if  $k \geq 3$. This is actually true if the considered complex  manifold is already a manifold in normal form (Theorem \ref{main}).
We also notice that Theorem \ref{main1} holds also in the $C^{\infty}$-smooth category and that its $C^{\omega}$ version follows from ancillary results needed as part of its proof.
\par
\smallskip
The structure of the paper is the following. In \S 2, we introduce the new notion of complex manifolds with boundary of circular type of class $\cC^{k,\a}$, a minor 
modification of the original definition of manifolds of circular type,  and we state the main results of our paper. In \S 3, we study  one-parameter families of Monge-Amp\`ere exhaustions of a  manifold of circular type in normal form, with centers at the point of a given  segment of the manifold.  We  then prove the existence of a  bijection between such families and  the so-called {\it abstract fundamental pairs}. These are   one-parameter families of  pairs, formed by non-standard complex structures  and vector fields satisfying appropriate  conditions.  Section \S 4 is dedicated to a detailed study of the vector fields occurring in the abstract fundamental pairs and in \S 5 it is given 
a crucial result on the existence of abstract fundamental pairs and the proof of the main theorem. \par
\section{Manifolds  of circular type}
\setcounter{equation}{0}
\subsection{Complex manifolds  with boundary of class $\cC^{k,\a}$}
\label{preliminaries}
Given  $k \in \bN \cup\{0, \infty\}$ and  $\a \in [0,1]$, a (real or complex) tensor field $T$ on  a (real) smooth manifold with boundary $\overline M = M \cup \p M$
 is said  to be  {\it of class $\cC^{k,\a}$ on $\overline M$} if all components  of $T$ in any system of  coordinates  of the structure of the  manifold with boundary $\overline M$ are of class $\cC^k$ and with $\a$-H\"older continuous $k$-th order derivatives up to the boundary. 
 \par 
 In this paper we need to exploit in many ways the relation provided by   the Newlander-Nirenberg Theorem between the notions of 
  {\it complex manifold $M$ of  dimension $n$}, i.e. a topological $2n$-dimensional manifold equipped with a complete  atlas $\cA$ of $\bC^n$-valued  homeomorphisms between open sets of $M$ and of $\bC^n$ with  holomorphic overlaps, and of {\it  integrable  complex structure   of class $\cC^{k,\a}$},  $k\geq 1$, $0 < \a < 1$,   on   a {\it real} $2n$-manifold  $M$, which is 
  a tensor field  $J$ of type  $(1,1)$ 
 of class $\cC^{k,\a}$ with $J_x^2 = - \Id_{T_x M}$ at each point  $x$ and vanishing Nijenhuis tensor $N_J$  i.e. such that for all vector fields $X,Y$ on $M$
 $$N_J(X, Y) = [X, Y] - [JX, JY] + J[X, JY]  +  J[JX, Y]=0.$$   
A $\cC^{k,\a}$-complex structure $J$ on $M$ induces the  direct sum decompositions 
$T^\bC_x M = T^{1,0}_x M\oplus T^{0,1}_x M$ of the  complexified tangent spaces into 
 $\pm i$-eigenspaces of the  linear maps $J_x:T_x M \to T_x M$.  
If $T^{1,0} M\subset T^{\bC} M$ is the complex bundle with fiber $T^{1,0}_x M$, $x \in M$,   a
    {\it set of $J$-holomorphic coordinates} is any  $\cC^1$-map
$F = (F^1, \ldots, F^n): \cU \subset M\longrightarrow \bC^n$, which is  homeomorphic onto its image and  satisfies, for any choice of   local generators $\{X_i\}_{1 \leq i \leq n}$ of $T^{1,0}_x M$:
 \beq \overline{X_i}(F^j) = 0\qquad \text{for any}\ 1 \leq i, j \leq n.\eeq 
 \par
By the  Newlander-Nirenberg Theorem (\cite{HT,Ma,NN,NW,We}), if $J$ is    $\cC^{k,\a}$   with   $k\geq 1$ and $0 < \a < 1$,  then  there exists  an atlas  $\cA_J$ of  $J$-holomorphic coordinates $F = (F^i)$ over $M$,    which makes  $(M, \cA_J)$  a complex manifold. Moreover, all $\bC$-valued functions $F^i$,  that are components of some set of 
 $J$-holomorphic  coordinates, are of class at least $\cC^{k+1, \a}$ relatively to any set of real coordinates of the  original real manifold structure   of  $M$. \par
Conversely,  if $(M, \cA)$ is a complex manifold of dimension $n$ (hence, a  manifold of real dimension $2n$), the real 
 tensor field $J$ having the form $ J = i \frac{\partial}{\partial  z^j} \otimes d  z^j -  i \frac{\partial}{\partial \overline z^j} \otimes d \overline z^j $  in any chart $(z^i)$ of the complex atlas $\cA$, is easily seen to be  a complex structure, for which the charts of $\cA$ are $J$-holomorphic coordinates. Such  a tensor field $J$  is of class $\cC^\o$ relatively to the charts in $\cA$ and, hence,   of class $\cC^{k,\a}$ relatively to any other atlas  $\cA'$ of real coordinates that   overlaps  with those of   $\cA$  in a $\cC^{k+1,  \a}$ fashion. 
\par
\medskip
For  $k  \geq 1$ and  $\a \in (0,1)$, if  $M$ is  a real $2n$-manifold  and  $J$ a  complex structure     of class $\cC^{k,\a}$ on $M$, we  call the pair $(M, J)$ a {\it  complex manifold of class $\cC^{k+1, \a}$}  and denote by $\cA_J$ the 
$\cC^{k+1,\a}$ atlas of $J$-holomorphic coordinates.
\par
\medskip
Let us now consider  a convenient analogue    of  real manifolds with boundary.  \par
\begin{definit} \label{definition21} {\rm Let  $\overline M = M \cup \p M$ be  a real $2n$-manifold with boundary.  A  {\it  complex structure on $\overline M$ of class $\cC^{k,\a}$, with $k\geq 1$, $\a \in (0,1)$},  is a triple $(J, \cD, J^\cD)$  where $J$ is a complex structure  on $M$ of class $\cC^{k,\a}$, $\cD$ is a $\cC^\infty$ codimension one  distribution  on $\p M$,   smoothly extendible to  a $ 2(n-1)$-dimensional  distribution   on a  tubular neighborhood $\cU \subset \overline M$ of $\p M$ and  $J^\cD$ is a tensor field  of type $(1,1)$ on the distribution $\cD$ of $\p M$,  with components of class $\cC^{k,\a}$  in any smooth frame field for the spaces $\cD_x$, $x \in \p M$, of the distribution $\cD$
subject to  the following  conditions: 
\begin{itemize}
\item[i)] the  smooth extension of $\cD$ on a  tubular neighborhood $\cU \subset \overline M$ of $\p M$, can be taken to be $J$-invariant; 
\item[ii)]  the restrictions   $J|_{\cD_x} $, $x \in \cU \setminus \p M$,  together with the tensors $J^\cD_y$, $y \in \p M$,   form a $(1,1)$-tensor field on the smooth extension of  $\cD$ in (i), with  components of bounded $\cC^{k,\a}$-norm   in any smooth frame field for the spaces of $\cD$. 
\end{itemize}
A manifold with boundary $\overline M$,  endowed  with a triple  $(J, \cD, J^\cD)$,   is  briefly   denoted  by  the pair $(\overline M, J)$.   We  may also write    $J^\cD = J|_{\cD_{\p M}}$ and we call $(\cD, J^\cD)$ the {\it CR structure of the boundary of $(\overline M, J)$}. In fact, by  boundary regularity assumptions,   the pair $(\cD, J^\cD)$ is an integrable CR structure on $\p M$.  We  say  that  the boundary  $\p M$  is {\it strongly  pseudoconvex} if the boundary CR structure has strictly positive Levi forms at all points. }
\end{definit}
Note that if $(N, \bJ)$ is   a complex manifold of class  $\cC^{k,\a}$, with $k\geq 1$ and $\a \in (0, 1)$, any relatively compact strongly pseudoconvex domain $D \subset N$ with smooth  boundary   is  a manifold with boundary,  equipped with a  complex structure $(J, \cD, J^\cD)$ of class $\cC^{k,\a}$. Indeed,  Definition \ref{definition21} is designed  to  capture  precisely the  properties of such domains.    
\par 
\subsection{Monge-Amp\`ere exhaustions} 
\label{section31}
\begin{definit} \label{MAdefinition} {\rm Let $(\overline M = M \cup \p M, J)$ be a complex manifold of class $\cC^\infty$ with a strongly pseudoconvex boundary $\p M$. 
Given    $k\geq 2$ and $\a \in (0,1)$, we call   {\it Monge-Amp\`ere exhaustion for $\overline M$ of class $\cC^{k ,\a}$ }  a continuous exhaustion $\t: \overline M \to [0, 1]$,  with   $\p M = \{x \in M\ :\  \t(x) = 1\}$,     such that:
\begin{itemize}[itemsep=8pt plus 5pt minus 2pt, leftmargin=18pt]
\item[i)]   The level set $\{\t = 0\} \subset M$  consists of  a single point $x_o$, called {\it center},  and    the pull back  $p^*(\t)$ of $\t$
 on the blow up $p: \wt {\overline M} \to \overline M$  at $x_o $,   is  of class $\cC^{k,\a}$;
\item[ii)] On  the complement $M \setminus \{x_o\} = \{0 < \t < 1\}$  of the center, the exhaustion $\t$ is a solution to the  differential problem
\beq \label{MAdifferentialproblem} \left\{\begin{array}{l} 2 i \partial \overline\partial \t = d d^c \t > 0\ ,\\[6pt]
2 i \partial \overline \partial \log \t = d d^c \log \t \geq 0\ , \\[6pt]
(d d^c \log \t)^n = 0\ \ (\text{Monge-Amp\`ere Equation)}\ \ ;\end{array}\right.
\eeq
\item[iii)]  In some (hence, in   {\it any}) system of  complex coordinates $z = (z^i)$ centered at $x_o$, the exhaustion $\t$ has a logarithmic singularity at $x_o$, i.e. 
$$\log \t(z) = \log\|z\| + O(1)\ .$$
\end{itemize}
}
\end{definit}
\smallskip
\par
Basic examples of Monge-Amp\`ere exhaustions of class $\cC^{k,\a}$ are given by the Minkowski functionals $\mu_D$ of the strongly pseudoconvex complete circular 
domains $D \subset \bC^n$ with boundary of class $\cC^{k, \a}$,  $k\geq 2$, $\a \in (0, 1)$.   In fact, by definition, $\mu_D$ is the function
$$\mu_D : \bC^n \longrightarrow  [0, + \infty)\ ,\qquad \mu_D(z) = \left\{\begin{array}{ll} 0 & \text{if}\ z = 0\ ,\\[7pt]  1/t_z & \text{if}\ z \neq 0\ ,
\end{array}\right.$$
where  $t_z = \sup\{t \in \bR :  t z \in D\}$, so that   $D = \{\mu_D < 1\}$ because $D$ is balanced.   Being $D$ strongly pseudoconvex,   
 $\t \= \mu_D^2|_D$
 is a Monge-Amp\`ere exhaustion  of  class $\cC^{k ,\a}$,   centered at $x_o = 0$. \par
 Other  crucial examples are given by (reparametrizations of) the Kobayashi distance functions of  the strictly  
convex domains   in $\bC^n$. 
Indeed, by   the results of  Lempert in  \cite{Le}, if $D\subset \bC^n$ is a strictly convex domain with a $\cC^{k + 2, \a}$-boundary,   $k \geq 2$,  each  point  $x_o  \in  D$ is the center of a  Monge-Amp\`ere exhaustion $\t^{(x_o)}: \overline D \longrightarrow [0, 1]$ of class $\cC^{k,\a}$,   where $\t^{(x_o)}(x)$ is   the  squared  hyperbolic tangent of the Kobayashi distance between $x$  and  $x_o$. The same is true for strictly linearly convex domains in $\bC^n$ in the  $C^{\infty}$ and the $C^{\omega}$ case  (\cite{Le1}).
\par
\medskip
Let   $\t: \overline M \to [0,1]$ be a  Monge-Amp\`ere $\cC^{k,\a}$-exhaustion    with center  $x_o$  and 
 \beq\label{defcappa}  \kappa:  T_{x_o} M \simeq \bC^n \longrightarrow \bR_{\geq 0}\ ,\qquad \kappa(v) \=\underset{{t_o\to 0}} \lim \left.\frac{d }{dt} \sqrt{\tau(\gamma_t)}\right|_{t = t_o}\ ,\eeq
 where  $\gamma_t$ is any smooth curve with $\gamma_0 = x_o$ and $\dot \gamma_0 = v$.  
Such  $\k$ is well defined and coincides with the  Kobayashi infinitesimal metric of $M$ at $x_o$ (\cite{Pt3}), it is  of class $\cC^{k,  \a}$ on $ T_{x_o} M \setminus \{0\}$
 and satisfies   $\kappa (\lambda v) = |\lambda| \kappa(v)$ for any $\lambda \in \bC\ $ so that the {\it (closed) indicatrix  at $x_o$ of  $((M,J), \t)$} defined by  
\beq \label{defindicatrix} \overline \cI \= \{v \in T_{x_o} M: \kappa(v) \leq 1 \} \subset T_{x_o}M  \simeq \bC^n\ ,\eeq
 is a complete circular  domain.
\par
\medskip
Let $p: \wt{\overline \cI} \to \overline \cI$ and $p': \wt{\overline M} \to \overline M$  be the blow-ups of 
$\overline \cI$ and $\overline M$   at $0$ and $x_o$, respectively.
A word by word repetition of the arguments   in \cite{Pt1},  where it was assumed  $\t$ to be of class $\cC^\infty$,
shows that {\it if the Monge-Amp\`ere exhaustion $\t$ is  of class $\cC^{k , \a}$, for some $k \geq 3$ and $\a \in (0,1)$,  there exists a unique $\cC^{k-2, \a}$-diffeomorphism $\Psi: \widetilde{\overline \cI}  \to \widetilde{\overline M}$}  
 satisfying   the following three conditions:
\begin{itemize}
\item[i)]  $\Psi|_{p^{-1}(0)} = \Id_{p'^{-1}(0)}$;  
\item[ii)] For each $t \in (0,1]$ the map 
$$\Psi^{(t)}: \partial \cI \rightarrow  \widetilde{\overline M}\ ,\qquad  \Psi^{(t)}([v], z) \= \Psi([v], t z)\ ,$$ 
 is a  diffeomorphism  between 
 $ \partial \cI$ and 
 $ \{\tau =  t^2\}$ mapping 
  the real distribution of   the CR structure of $\partial  \cI$ into the real distribution 
of   the CR structure of  $ \{\tau =  t^2\}$; 
\item[iii)] For each $([v], z) \in \partial  \cI$,   the map 
$$\widetilde f^{([v], z)}  :\overline \Delta \longrightarrow \widetilde{\overline M} \ ,\qquad \widetilde f^{([v], z)}( \zeta) \= \Psi([v],  \zeta z)\ ,$$  
is proper holomorphic and  injective and 
$\widetilde f^{([v], z)}(\Delta \setminus \{0\})$ is an integral leaf of the distribution $\cZ$, given  by  the spaces   $\cZ_x = \ker d d^c \t_x$. 
\end{itemize}
Such a map $\Psi$  is called {\it circular representation of  $(\overline M, J)$} determined by  $\t$. 
\par
\medskip
The proof shows also that  this circular representation $\Psi$  is such that: 
\begin{itemize}
 \item[a)] The projection  onto $\overline{\cI}$ of the pulled-back exhaustion  $\t \circ \Psi: \wt{\overline{\cI}} \to [0,1]$    coincides with  the  Kobayashi infinitesimal metric  $\k$   and  is 
  of class $\cC^{k -2, \a}$ on $\overline \cI \setminus \{0\}$; 
 \item[b)] If  $k \geq 4$, the pulled-back  complex structure  $\wt J' = \Psi^*(J)$  on $\wt{\overline \cI}$ is integrable and of 
class  $\cC^{k-3, \a}$, i.e.  $(\wt{\overline \cI}, \wt J')$   is a $\cC^{k -2, \a}$ complex manifold. 
\end{itemize}
By classical results  on blow-ups and Remmert reductions, the complex manifold $(\wt{\overline \cI}, \wt J')$  has a blow-down, which can be naturally identified with the manifold with boundary
$\overline \cI$, equipped with  an appropriate non-standard  atlas $\cA$ of complex charts. Denoting by $J'$  the tensor field,    for which   $\cA$  is the atlas of $J'$-holomorphic coordinates, we  conclude that  {\it the $\cC^{k -2,\a}$-diffeomorphism $\Psi: \wt{\overline \cI} \to \wt{\overline M}$ determines  a $(J', J)$-biholomorphism   $\Psi:  (\overline \cI, J') \longrightarrow (\overline M, J)$  and  induces  the $\cC^{k-2,  \a}$ Monge-Amp\`ere exhaustion $\k \= \t \circ \Psi$   on   $ (\overline \cI, J')$.} \par
 \smallskip
We stress that, despite of the fact that  the tensor field $J'$ on $\overline \cI$  has  smooth components in each chart of the atlas  $\cA$, such tensor field has in general  non-smooth components  in the standard coordinates of $\cI \subset \bC^n$. Indeed  they  are  in most cases not even differentiable at $0$. Nonetheless $J'$  has   the same  regularity of  $\wt J'$  at the points of  $\overline{\cI} \setminus \{0\} = \wt{\overline \cI} \setminus p^{-1}(0)$ -- we refer to \cite{PS2} for a more detailed discussion of all this (here and in all what follows, for any blow up we tacitly  identify the complementary region of the exceptional divisor with its  image in the blow down). \par
\subsection{Manifolds of circular type   with boundary}
\label{section23}
The discussion of the previous section leads to  the following notion.
\begin{definit} \label{closeddef}  {\rm Given    $k\geq 2$, $\a \in (0,1)$, we call  {\it  manifold of circular  type with boundary of  class $\cC^{k,\a}$ }  (or, simply,   {\it  $\cC^{k,\a}$-manifold of circular type})
a  pair  $((\overline M, J), \t)$ formed by 
\begin{itemize}
\item 
a complex manifold $(\overline M, J)$ of class $\cC^{k, \a}$ with strongly pseudoconvex boundary,   diffeomorphic to the closed unit ball $\OB = \bB^n \cup \p \bB^n$ of $\bC^n$; 
\item a  Monge-Amp\`ere  exhaustion  $\t: \overline M \longrightarrow [0,1]$  of class $\cC^{k,\a}$ in each set of  $J$-holomorphic coordinates.
\end{itemize}
}
\end{definit}
\par 
The  main example  of a  $\cC^\infty$-manifold of circular type  to keep in mind is   the pair $((\OB, \Jst), \t_o))$, where $\OB \subset \bC^n$ is the  closed unit ball centered at $0$,   equipped with    the standard complex structure $\Jst$  of $\bC^n$ and
  the standard Monge-Amp\`ere exhaustion 
$\t_o(x) = \|x\|^2$ where $\|\cdot\| $  is the Euclidean norm of $\bC^n $.
This basic example comes with two special distributions in the tangent space $ T (\bB^n\setminus \{0\})$, which  have direct  analogues in any other manifold of circular type. The first is the $\Jst$-invariant  distribution  $\cZ \subset T (\bB^n\setminus \{0\})$ of the spaces 
 \beq \label{distrbis} \cZ_x\= \{v \in T_x \bB^n:(\p \bar \p \log \t_o)(v, \cdot) = 0\} = 
 Span_{\bC}\left(\left. \hskip-.1cm z^i \frac{\p}{\p z^i}\right|_x \right)\hskip-.1cm.
 \eeq
   The second is the   distribution $\cH$ of the complementary subspaces
\beq \label{distrter} \cH_x = (\cZ_x)^\perp = \{v \in T_x\bB^n\,:\ < v,\cZ_x> = 0\}\ , \eeq 
where $< \cdot, \cdot> $ denotes the standard  Euclidean inner product of $\bC^n$. \par
 These  distributions  are not defined at $0$
 but they both admit  smooth extensions  at all points of  the exceptional divisor of the blow up $p: \wt \bB^n \longrightarrow \bB^n$ at $0$. 
 Note also that $\cZ$ is integrable, with integral leaves given by the straight  disks   through $0$. On the contrary,  $\cH$  is not integrable.  In fact, the  restriction of $\cH$ to each sphere $S_r \= \{\|x\| = r\}$, $0 < r \leq 1$,  coincides with the  contact distribution underlying the CR structure of such sphere. 
Both distributions  $\cZ$ and   $\cH$ are   $\Jst$-invariant, so that   the  complex distributions  $\cZ^{10}, \cZ^{01} \subset \cZ^\bC$  and $\cH^{10}, \cH^{01} \subset \cH^\bC$,  
 of the $(+i)$- and $(-i)$-eigenspaces of $\Jst$, are well defined at each point of $\wt \bB^n$.  The complex distributions  $\cH^{10}, \cH^{01}$ admit smooth extensions up to  the boundary, where they determine the standard CR structure of $\p \bB^n$. \par
\medskip
We now introduce a crucial  class of deformations of   $((\OB, \Jst), \t_o)$.\par
\begin{definit} \label{definition5.1}  {\rm
 Let $k \geq 2$ and $\a \in (0,1)$. We call {\it L-complex structure of class $\cC^{k-1, \a}$}  a complex structure $\wt J$ of class $\cC^{k-1,\a}$ on the blow up $\TOB$  of $(\OB, \Jst)$ at the origin,   satisfying the following conditions:
\begin{itemize}
\item[i)] it leaves invariant all spaces of the distributions $\cZ$ and  $\cH$; 
\item[ii)] $\wt J|_{\cZ} = \Jst|_{\cZ}$; 
\item[iii)] there exists a  homotopy $\wt J_t$, $t \in [0,1]$,  of class $\cC^{k-1, \a}$ in the parameter $t$, between  $\wt J_{t = 0} = \Jst$ and $\wt J_{t = 1}= \wt J$, such that each $\wt J_t$ is a  complex structures of class $\cC^{k-1, \a}$   satisfying (i) and (ii).  
\end{itemize}
If $\wt J$ is  an L-complex structure of class $\cC^{k-1,\a}$ on $\TOB$,  let $J$ be the (non-standard) complex structure on $\OB$ that makes $(\OB, J)$ the blow-down of the complex manifold  with boundary $(\TOB, \wt J)$ (it exists by the above mentioned  facts on Remmert reductions).  By the results in \cite{PS2} and  the following remark,  each  pair of the form $((\OB, J), \t_o)$,  with $J$ coming from an $L$-complex structure  $\wt J$ as above,    is a $ \cC^{k,\a}$-manifold of circular type,   called  {\it in normal form}. }
 \end{definit}
\begin{rem} \label{remark25}{\rm The fact that  $\t_o$ is of class  $ \cC^{k,\a}$ in the  charts of $(\OB\setminus\{0\}, J)$  can be easily checked as follows. By construction, the components in standard coordinates of the complex structure $J$  are  of class $\cC^{k-1, \a}$  at all points of $\OB \setminus \{0\} = \TOB \setminus p^{-1}(0)$. 
 This implies that   each chart of $J$-holomorphic coordinates of $(\OB \setminus \{0\}, J)$   overlaps in a $\cC^{k, \a}$ way  with the standard coordinates (see \S \ref{preliminaries}).  Being $\t_o= \|\cdot \|^2$ of class $\cC^\infty$ in the standard coordinates outside the origin, it must be  of class $\cC^{k,\a}$  when it is expressed in  $J$-holomorphic coordinates. }
 \end{rem} \par
\medskip
The interest for the manifolds in  normal forms comes from  the fact that, for $k \geq 4$, any $\cC^{k,\a}$-manifold of circular type   is   biholomorphic with a  $\cC^{k-2,\a}$-manifold of circular type in normal form. 
This property was  proven  in  \cite{PS}  just  for the case of $\cC^\infty$-manifolds of circular type, but a careful   check of all arguments in that paper shows that they actually go through for each $\cC^{k,\a}$-manifold of circular type, provided that  $k \geq 4$. More precisely: \par
\smallskip
 (1) By the remarks in \S \ref{section31},  if $((\overline M, J), \t)$ is  a  $\cC^{k,\a}$-manifold of circular type  with $k\geq 4$, $\a \in (0, 1)$,  the  circular representation  determines  a $(J,J')$-biholomorphism 
from $(\overline M, J)$ into the $\cC^{k- 2, \a}$-manifold of circular type $((\overline \cI, J'), \k)$,   given by the closed  indicatrix   $\overline \cI$ of the Kobayashi infinitesimal  metric $\k$ at the center  $x_o$,  an appropriate complex structure $J'$ of class $\cC^{k - 3, \a}$ on $\overline{\cI}$ and the exhaustion   $\k$, which is   $\cC^{k, \a}$ in the complex charts of $(\overline \cI, J')$. 
\par
\smallskip
 (2) Since  
 $d d^c \k$ 
is of class $\cC^{ k- 2,  \a}$ and $k \geq 4$, 
 all arguments of the proofs of  Moser's Theorem in \cite{Mo} and of  Lemma  3.5 and Thm. 3.4 in \cite{PS} remain valid and yield to the existence 
of:
\begin{itemize}
\item[i)]  A smooth family of $\cC^{k- 2, \a}$ diffeomorphisms $\Phi_t: \wt{\bC^n} \to \wt{\bC^n}$  from  the blow up $\wt{\bC^n}$ of $\bC^n$ at $0$ into itself, satisfying all   conditions  of Lemma 3.5 in \cite{PS}; 
\item[ii)] A smooth isotopy of L-complex structures $J_t = \Phi_{t*}(\Jst)$ of class $\cC^{ k-3, \a}$ between 
$J_0 = \Jst$ and a (non-standard) complex structure $J_1 = J''$  that 
makes  $((\TOB, J''), \t_o)$ biholomorphic to the blow up $(\wt{\overline \cI}, J')$ of $(\overline{\cI}, J')$ at the origin. Such a complex structure $J''$ projects  onto a (non-standard) complex structure $
J''$ on $\OB$ that makes $((\OB, J''), \t_o)$ a $\cC^{k- 2,\a}$-manifold of circular type in normal form,  biholomorphic to $((\overline{\cI}, J'), \k^2)$.
\end{itemize}
\par
\smallskip
Combining   (1) and (2) (as it is done in \cite{PS}, Thm. 3.4),   one gets 
 \begin{theo} \label{existenceanduniqueness} Let $k \geq 4$, $\a \in (0,1)$ and $((\overline M, J), \t)$ be  a  $\cC^{k,\a}$-manifold of circular type with  center  $x_o$.  Then, there is a $(J, J')$-biholomorphism 
 $\Phi: (\overline M, J) \to (\OB, J')$  from $((\overline M, J), \t)$ to a  $\cC^{k-2,\a}$-manifold in normal form $((\OB, J'), \t_o)$ with the following properties: 
 \begin{itemize} 
 \item[a)] $\Phi(x_o) = 0$ and $\t = \t_o \circ \Phi = \Phi^*(\t_o)$; 
 \item[b)] $\Phi$ maps the integral leaves of the  distribution on $M \setminus \{x_o\}$
 \beq \label{distquater} \cZ^{J,\t} \= \hskip-10pt \bigcup_{x \in M \setminus \{x_o\}} \cZ^{J,\t}_x\ , \quad\cZ^{J,\t}_x \= \{v \in T_x M\,:\,d d^c_J \log \t(v, \cdot) = 0\}\eeq
 into the integral leaves of the   distribution $\cZ$ of $\bB^n$, that is into the straight disks  through the origin of $\bB^n$. 
 \end{itemize}
 \end{theo}
 A biholomorphism $\Phi: (\overline M, J) \longrightarrow (\OB, J')$,   mapping a $\cC^{k,\a}$-manifold of circular type into one  in normal form,    satisfying    (a) and (b)  of the above theorem, is called {\it normalizing map}.
  \par
 \subsection{Statement of  the main result}
All results  of this paper are consequence of the following theorem on manifolds in normal form, the  proof of which is divided    in the remaining three sections.
  \begin{theo}\label{main} Let $((\OB, J), \t_o)$ be a  manifold of circular type in normal form  of class $\cC^{k}$ with $k \geq 3$. Then, for each $v \in \bC^n$ with $\|v\| = 1$ there exists a  value $\l_v \in (0, 1]$ such that: 
\begin{itemize}
\item[a)] For each $\l \in (0,\l_v)$ and the corresponding parameterized segment $x_t \= t  \l v$ between $0$ and $x_1 = \l v$,  there is a one-parameter family of Monge-Amp\`ere exhaustions $\t^{(x_t)}: \OB \longrightarrow [0,1]$, with $t \in [0,1]$, centered  at  $x_t$  and of   class $\cC^{k-1,\a}$ for each  $\a \in (0,1)$. The restrictions of  the   $\t^{(x_t)}$, $t \in I \subset [0,1]$,  on any  open subset $\cV \subset D$ non containing the centers $x_t$, $t \in I$,   gives a real function on $\cV \times I$ which is  of class $\cC^{k-1,\a}$ in all of its arguments, i.e. also with respect to $t$. 
\item[b)] If $\l_v  < 1$, there is no  Monge-Amp\`ere exhaustion of class $\cC^{k-1,\a}$ for any $\a \in (0,1)$, which is  centered at  the point $\wh x = \l_v v$.
\end{itemize}
The set   $\cS\= \{z \in \OB\ :\ z = \l_v v\ , \ \|v\| = 1\}$   is the boundary  of an open neighborhood $\cU \subset \overline \bB^n$ of $0$ and  is   the projection onto $\TOB \setminus p^{-1}(0)$ of the intersection of two   submanifolds, uniquely determined by the complex structure $J$,   of the bundle  $\cH^{\bC *} \otimes \cH^{\bC}|_{\TOB \setminus p^{-1}(0)} \to \TOB \setminus  p^{-1}(0)$.
 \end{theo}
 Due  to Theorem \ref{main}, this  immediately implies our main result, which is 
  \begin{theo}\label{main1} Let  $k \geq 5$  and $((\overline M, J), \t)$ be a $\cC^{k}$-manifold of circular type. Then there exists an open neighborhood $M'$ of the center $x_o$ of $\t$ with the property that for each $x \in M'$  there exists a Monge-Amp\`ere exhaustion $\t^{(x)}: \overline M \longrightarrow [0,1]$, of class $\cC^{k-3,\a}$ for any $\a \in (0,1)$, with center $x$. The dependence of the exhaustions $\t^{(x)}$  on the center $x$ is also of class $\cC^{k-3,\a}$. 
 \end{theo}
By uniqueness of the Monge-Amp\`ere exhaustions with  given centers,  Theorem \ref{main1} shows that if a manifold with boundary of circular type $\overline{M}$  admits a Monge-Amp\`ere exhaustion of class  $\cC^\infty$, hence of class $\cC^{k}$ for each $k$, it is equipped with 
$\cC^\infty$ Monge-Amp\`ere exhaustions,  centered at all  points of an appropriate open subset  $M'$.  
 In fact, as we shall point out later on,  the determination of the maximal open set  $M'$ with the properties described in Theorem \ref{main1} is independent of $k$ 
as long as $k\geq 5$. Furthermore some of the steps in the proof of Theorem \ref{main} show that its $\cC^\omega$ version holds.
Thus, we also have:
 \begin{theo}\label{main2} Let   $((\overline M, J), \t)$ be a $\cC^{\infty}$-manifold (resp. $\cC^{\omega}$-manifold) of circular type. Then there exists an open neighborhood $M'$ of the center $x_o$ of $\t$ with the property that for each $x \in M'$  there exists a Monge-Amp\`ere exhaustion $\t^{(x)}: \overline M \longrightarrow [0,1]$, of class $\cC^{\infty}$ (resp. $\cC^{\omega}$), with center $x$. The dependence of the exhaustions $\t^{(x)}$  on the center $x$ is also of class $\cC^{\infty}$ (resp. $\cC^{\omega}$). 
 \end{theo}
Direct consequences of this   are  Theorem  \ref{main0} and the properties of pluricomplex Green functions, which have been  
discussed in the Introduction.  
 \par
\section{One-parameter  families of  Monge-Amp\`ere exhaustions}
\label{section3}
\setcounter{equation}{0}
\subsection{Quasi-diffeomorphisms and  quasi-regular vector fields}
\label{sect41}
Let $k\geq 2$,  $\a \in (0,1)$ and assume that  $(\overline M, J)$ and $(\overline{M'}, J')$  are two complex $n$-manifolds with boundary of class $\cC^{k,\a}$ and let
$$\pi: \wt{\overline M} \longrightarrow \overline M\ ,\qquad \pi': \wt{\overline M'} \longrightarrow \overline{M'}$$ 
be  the   blow-ups of  $\overline M$ and $\overline{M'}$ at some fixed interior points $x, x'$ with {\it exceptional divisors}  $\pi^{-1}(x) \subset \wt{\overline M}$ and $\pi'{}^{-1}(x') \subset \wt{\overline {M'}}$   both biholomorphic to $\bC P^{n-1}$. 
\begin{definit} {\rm
A  $\cC^{k,\a}$-diffeomorphism $\wt F: \wt{\overline M}\to \wt{\overline{M'}}$  is called {\it tame} if it 
maps  diffeomorphically the exceptional divisor of $\wt{\overline M}$  onto the exceptional divisor of  $\wt{\overline{M'}}$, so that 
$\wt F$ induces   a homeomorphism $F: \overline M \to \overline{M'}$,   mapping $x$ into $x'$, which  is of class $\cC^{k,\a}$ on $\overline M \setminus \{x\}$. \par
We call any such homeomorphism a  {\it quasi-diffeomorphism between $\overline M$, $\overline{M'}$ of class $\cC^{k,\a}$, pointed at  $x$}. The associated tame diffeomorphism  $\wt F$,  inducing  $F$,  is called {\it    tame lift of $F$}.
}
 \end{definit}
  \par
 Among the various  reasonable  infinitesimal counterparts of the notion of quasi-diffeomorphism, the
  one that better  fits with our purposes  (and which we formally introduce at the end of this section)  is rooted  in the following observations.   Let $x_t:[0,1] \longrightarrow M = \overline M \setminus \p M$ be a $\cC^2$ curve entirely included in the interior $M$.  For each $t \in [0,1]$,  consider   the  blow-up $\pi_t: \wt{\overline M}_{t} \longrightarrow \overline M$ of $\overline M$ at   $x_t$.
A one-parameter family $F_t: \overline M \longrightarrow \overline M$ of quasi-diffeomorphisms of class $\cC^{k,\a}$ pointed at the $x_t$  induced by  a  family of   tame diffeomorphisms  
$\wt F_t:   \wt{\overline M}_{t} \longrightarrow \wt{\overline M}_{1}$ 
of class $\cC^{1,\b}$,    $\b > 0$, in  the coordinates of  $\bigcup_{t \in [0,1]} \wt{\overline M}_{t}  \simeq  \wt{\overline M}_{1}\times [0,1]$, 
is  called   {\it $\cC^{1,\b}$-family of  $\cC^{k,\a}$ quasi-diffeomorphisms, pointed  at the  points  $x_t$}. 
\par
 For any  such  family,  at a fixed  $t \in [0,1]$, we may consider 
the vector field $X_{t}$, defined at the points  $x \in M \setminus \{x_1\}$  by 
\beq \label{quasi-field} X_{t} \big|_x\= \frac{dF_{t+s}(F^{-1}_{t}(x))}{ds}\bigg|_{ s= 0}\ .\eeq
It is a vector field  of class $\cC^{k,\a}$ on $\overline M  \setminus \{x_1 = F_t(x_t)\}$, whose  restriction to $\p M$ is always tangent to the boundary and of class $\cC^{k,  \a}$. Moreover, if we extend $X_t$  at the point $x_1$  by setting  (for the existence of this derivative, see the argument below)
\beq \label{limitX} X_{t}|_{x_1} \= \frac{d F_{t+s}(x_t)}{ds}\bigg|_{s = 0}\ ,\eeq
we get a vector field over the whole $\overline M$ with the following property: {\it for each sequence $y_k \to x_1$, with $y_k \in \overline M \setminus \{F_t(x_{s}), s \in   [0,1]\}$,  }  
$$\lim_{k \to \infty} X_t|_{y_k} =  X_{t}|_{x_1}\ .$$
Indeed, this can   
  be checked as follows.  At each point $y_k$,   the vector $X_t|_{y_k}$  is the tangent vector at $s = 0$ of   
  the curve  $\h_s^{(t| y_k)}  \= F_{t + s}(F^{-1}_{t}(y_k))$ with $s$  in some  fixed small interval $[- \ve, \ve]$.  Note that, by  the assumption on the sequence $y_k$, for each $s \in [- \ve, \ve]$, the point $\h_s^{(t| y_k)}$ satisfies 
  $$ F_{t + s}^{-1}( \h_s^{(t| y_k)})  = F^{-1}_{t}(y_k) \neq x_{t+s}\qquad \Longrightarrow \qquad  \h_s^{(t| y_k)} \neq  F_{t + s}(x_{t+s}) = x_1\ .$$
  This implies that  each curve $\h_s^{(t| y_k)}$  admits a unique  $\cC^{1,\b}$  lifted curve  $\wt \h_s^{(t|y_k)}$ on  the blow up  $\wt{\overline M}_{1} $. Moreover, 
by the regularity assumptions on the family of diffeomorphisms $F_t$,  when   $y_k \to x_1$,  the  lifted curves  $\wt \h_s^{(t|y_k)}$ tend uniformly in   $\cC^1$-norm to a   curve 
$\wt \h$ in $\wt{\overline M}_{1}$.  Such a limit curve     projects onto  the $\cC^{1, \b}$ curve  
$$ \h_s  = F_{t+s}(F^{-1}_{t}(x_1))= F_{t+s}(x_t)\ .$$
 This implies that the vectors $X_{t}|_{y_k}$ (that are  the  tangent vectors at $s = 0$ of the curves  $\h^{(t|y_k)}_s$)  tend to the  tangent vector of $\h_s$ at $s = 0$, that is to    \eqref{limitX}.  \par
\smallskip
All this motivates the following
\begin{definit} \label{quasi-vector} {\rm Let $(\overline M , J)$  be a  complex manifold of dimension $n$ with boundary of class $\cC^{k,\a}$   and $x \in M= \overline M \setminus \p M$. A  vector field $X$ on $\overline M $ is called {\it quasi-regular of class $\cC^{k,\a}$,  pointed at $x$}, if:
\begin{itemize}
\item[i)]  It is of class  $\cC^{k, \a}$ on $\overline M \setminus \{x\}$ and tangent to $\p M$ at  boundary points; 
\item[ii)] There is an open and dense subset of $\cU \subset \overline M \setminus \{x\}$,  such that $\lim_{k \to \infty} X|_{y_k} = X|_x$ for each sequence $y_k \in \cU$ converging  to $ x$. 
\end{itemize}
}
 \end{definit}
\subsection{Curves of centers  and  families of quasi-diffeomorphisms} \label{section32}
Consider now a manifold of circular type in normal form  $((\OB, J), \t_o)$  of class $\cC^{k,\a}$, with  $k\geq 1$, $\a \in (0, 1)$. Let also  $v \in \bB^n \setminus \{0\}$ and $x_t$ the straight curve
\beq \label{segment}x_t = t v: [0, 1] \longrightarrow \OB \ .\eeq
 If $(\OB, J)$ admits  Monge-Amp\`ere exhaustions $\t_t$   of class $\cC^{k, \a}$, $k \geq 4$, $\a \in (0,1)$, centered at  the points $x_t$, then it also has   a one-parameter family of normalizing maps  
$\Phi_t: ((\OB, J), \t_t) \longrightarrow ((\OB, J_t), \t_o)$
for   appropriate  non-standard complex structures $J_t$. On the base of this, we  consider the following notion.
\begin{definit}  {\rm Let $k\geq 2$ and $\a \in (0,1)$. We call   {\it curve of Monge-Amp\`ere quasi-diffeomorphisms of class $\cC^{k,\a}$,  guided by  the curve \eqref{segment}}, 
a  $\cC^{k,\a}$-family of  quasi-diffeomorphisms $\Phi_t$ of $(\OB, J)$,  each of them   pointed  at    $x_t = t v$, $t \in [0,1]$, 
with the property that $\Phi_t(x_t) = 0$ for each $t$  and  satisfying the following condition:  each pushed-forward complex structure  $\wt J_t = \Phi_{t*}(\wt J)$ 
of the complex structure $\wt J$ of the blow-up $\pi_t: \TOB_{t} \to \OB$ at the point $x_{t}$,  is an L-complex structure, so that  its projected complex structure $J_t$   
determines  a normal form $((\TOB, J_t), \t_o)$. }
\end{definit} 
We have the following\par
  \begin{prop} \label{theorem43}  If $\Phi_t$ is a curve of Monge-Amp\`ere quasi-diffeomorphisms of class $\cC^{k,\a}$,   guided by  \eqref{segment}, then for each $t \in [0, 1]$  the function
 $$\t_t: \OB \longrightarrow [0,1]\ ,\qquad \t_t(x) \= \left\{\begin{array}{ll} 0& \text{if}\ x = x_t\\[8pt]
 ( \t_o \circ \Phi_t)(x)& \text{otherwise}
  \end{array}\right.
  $$
 is a Monge-Amp\`ere exhaustion of $(\OB, J)$ of class $\cC^{k,\a}$,  centered at $x_t$. 
 \end{prop} 
 \begin{pf} By  the  properties of $\t_o$,  each exhaustion  $\t_t$  satisfies (i) and (iii) of Definition \ref{MAdefinition}.  Hence,   
we only need to show that each $\t_t$ satisfies  also   (ii)  on $(\OB, J)$. This is in turn
 equivalent to prove that $\t_o = \t_t \circ \Phi^{-1}_t$ satisfies  \eqref{MAdifferentialproblem}  on the complex manifold $(\OB, J_t)$.  For this, we need to  the following  lemma.\par
 \begin{lem}\label{crucial5}  If $J$ is the complex structure of a normal form $((\OB, J), \t_o)$
  then 
 $d d^c_{J} \t_o = d d^c_{\Jst} \t_o$ at all points different from $0$.
 \end{lem}
 \begin{pfns}  Let $\wt J$ be  the L-complex structure on $\TOB$ which induces the complex structure $J$ on $\OB$. We recall that $\wt J$ preserves both distributions $\cZ$,  $\cH$ and that $J|_{\cZ} = \wt J|_{\cZ} = \Jst|_{\cZ}$. Hence,  if  we denote by $(\cdot)^\cZ$ the natural projection of  each tangent space $T_x \bB^n = \cZ_x \oplus \cH_x$ onto the subspace  $\cZ_x$, we have that
 \begin{multline}
 \nonumber d d^c_{J} \t_o(X, Y) = -  X(J Y (\t_o)) + Y(JX (\t_o)) + J([X, Y])(\t_o) \\
 \phantom{vvvv}=  -  X(J (Y)^\cZ (\t_o)) + Y(J (X)^\cZ (\t_o)) + J([X, Y]^\cZ)(\t_o) \\
  \phantom{vvvvvvv}= -  X(\Jst (Y)^\cZ (\t_o)) + Y(\Jst (X)^\cZ (\t_o)) + \Jst([X, Y]^\cZ)(\t_o)  \\
   =   d d^c_{\Jst} \t_o(X, Y). \phantom{vvvvvvvvvvvvvvvvvvvvvvvvvvvvvvvvvvvvvvv}  \qed
 \end{multline}
\end{pfns}
\par
\medskip
 By this lemma, the $2$-form $d d^c_{J_t} \t_o$ is nowhere degenerate and hence it is positive definite  
 on $(\OB \setminus \{0\}, J_t)$, since $J_t$ is isotopic to $\Jst$. 
 Similarly, one has that $d d^c_{J_t} \log \t_o \geq 0$ and that $d d^c_{J_t} \log \t_o$ satisfies the Monge-Amp\`ere equation on $\OB \setminus \{0\}$. 
 Thus all conditions  of \eqref{MAdifferentialproblem} are satisfied. 
\end{pf}
 \begin{prop} \label{easylemma}  If $\Phi_t$ is a curve of Monge-Amp\`ere quasi-diffeomorphisms of class $\cC^{k, \a }$, $k\geq 2$, $\a \in (0,1)$, the corresponding  family of complex structures  $J_t \=\Phi_{t*}(J)$    is such that, for each $t \in [0,1]$,  on  $\OB \setminus\{0\}$ one has 
\beq \label{crucial2}  \frac{d J_t}{dt}  = - \cL_{X_t} J_t \qquad\text{where}\ X_t\ \text{is defined by}\quad X_t |_x\=   \left.\frac{d  \Phi_{s}}{ds} \right|_{(\Phi^{-1}_t(x),t)}\ .\eeq
\end{prop}
\begin{pf}  
First of all, consider  the one-parameter family of quasi-regular vector fields   $Y_t$,  defined  at each point  $x \in \OB \setminus \{x_t\}$ by 
$$Y_t |_x{\=} \left.\frac{d  \Phi^{-1}_{s}}{ds} \right|_{(\Phi_t(x),t)}\ .$$
Note that,  since $\Phi_s \circ \Phi^{-1}_s = \Id_{\OB}$ for each $s$,   by taking the derivatives of both sides with respect to $s$  at a fixed  $x$,   we have
$$0 = \frac{d \Phi_s}{d s}\bigg|_{(\Phi^{-1}_{t}(x), t)} + \Phi_{t *}\left(\frac{d \Phi^{-1}_s}{d s}\bigg|_{(x, t)}\right) = X_{t}|_{x} +  \Phi_{t*}\left(Y_{t}|_{\Phi^{-1}_t(x)}\right),$$
from which we infer that, for each given $t$ and $x$,  $Y_{t}|_{x}$ and $X_{t}|_{x}$ are  related by 
\beq\label{crucial1} X_{t}|_{x} = - \Phi_{t*}(Y_{t})|_{x}\ .\eeq
We now recall that  each complex structure   $J_t$,  $t\in [0,1]$,  is well defined and of class $\cC^{k',\a}$, with $k' = k -1\geq 1$, $\a \in (0,1)$.  Actually, by  the considered regularity assumptions at $t = 0$ and $t = 1$, we may  assume that $J_t$ is well defined and $\cC^{k', \a}$ for each $t$ in a slightly larger open interval of $[0,1]$, say $(-\ve, 1 + \ve)$. 
Let us pick a fixed value $t \in [0, 1]$ and, for each $s$ close to $t$,  consider the  the maps $F_{s} \= \Phi_s \circ \Phi^{-1}_t$, $F^{-1}_{s} \= \Phi_t \circ \Phi^{-1}_s$. 
By the definition of  $J_s$, we have
 \beq \label{FJ} J_{s} = \Phi_{s*}(J) = F_{s*} (\Phi_{t*}(J) )=  F_{s*} (J_t)\ .\eeq
We also remark that, from   \eqref{crucial1}, 
 \beq
 \begin{split}\label{3.5} & \frac{d F_s}{ds}\bigg|_{(F^{-1}_t(x),  t)} \!\!\!= \frac{d \Phi_s}{ds}\bigg|_{(\Phi^{-1}_t(x),  t)} = X_{t}|_{x}\ ,\\
& \frac{d F^{-1}_s}{ds}\bigg|_{(x, t)} \!\!\! =  \Phi_{t*} \left(\frac{d \Phi^{-1}_s}{ds}\bigg|_{(x, t)}\right) = \Phi_{t*} \left( Y_t|_{\Phi^{-1}(x)}\right) = \Phi_{t*}( Y_t)|_{x}  = - X_t|_x\ .
 \end{split}\eeq
 Consider now  a  (local) system of real coordinates    $\xi = (x^i)$  on an open set,  on which the tensor field $J_s$, the vector field $Y_s$ and the map  $\Phi_s$ are of class $\cC^{k', \a}$ for all $s$ in a small open interval $(t - \d, t + \d) $,  entirely included in $(-\ve, 1 + \ve)$.  In these  coordinates, $J_t$,      $X_t$  and the maps $F_{s}$, $F^{-1}_{s} $,  have the form
$$J_t = J^{\ i}_{t j} \frac{\partial}{\partial x^i}\otimes dx^j\ ,\qquad X_t = X^i_t \frac{\p}{\p x^i}\ ,$$
  $$F_{s} \= \big(F^1_s(x^1, \ldots, x^{2n}), \ldots, F^{2n}_s(x^1, \ldots, x^{2n})\big)\ ,$$
$$F^{-1}_{s} = \big( (F^{-1}_s)^{1}(x^1, \ldots, x^{2n}), \ldots, (F^{-1}_s)^{2n}(x^1, \ldots, x^{2n})\big)\ .$$
We now observe that,  by \eqref{FJ}  and  the classical coordinate expressions of  push-forwards, the components $J^{\ i}_{sj}$ of $J_s|_x$ have the form  
\beq
\label{santa}
J^{\ i}_{sj}|_x = \frac{\partial F^i_s}{\partial x^\ell}\bigg|_{F^{-1}_{s}(x)}\!\! \frac{\partial (F^{-1}_s)^{m}}{\partial x^j}\bigg|_{x} \!\!J^{\ \ell}_{t m}|_{F^{-1}_{s}(x)}\ .
\eeq
So, differentiating \eqref{santa} with respect to  $s$ at $s = t$, from  \eqref{3.5} we get 
\beq
\label{anti-vigilia}
\frac{d}{ds}(J^{\ i}_{tj}|_x )\bigg|_{t}=   \frac{\partial X^i_t}{\partial x^\ell}\bigg|_{x} J^{\ \ell}_{tj}|_x- \frac{\partial X^m_t}{\partial x^j}\bigg|_{x} J^{\ i}_{tm}|_x - \frac{\partial J^{\ i}_{tj}}{\partial x^r} \bigg|_{x} X^r_t|_x  \ .
\eeq
As the right hand side  is the coordinate expression of  $-\cL_{X_t}J_t|_x$, the claim follows.
\end{pf}
The one-parameter family of pairs $(J_t, X_t)$, given  by  the  complex structures $J_t = \Phi_{t*}(J)$ and   the vector fields   \eqref{crucial2},  is called {\it fundamental pair of the Monge-Amp\`ere quasi-diffeomorphisms    $\Phi_t$}.  
Note that, if we extend each vector field $X_t$  at $0$ by setting 
\beq X_t|_0 \=  \left.\frac{d  \Phi_{t +s}(x_{t})}{ds} \right|_{s = 0} \ ,\qquad x_t = t v\ ,\eeq
then {\it each $X_t$ is  a quasi-$\cC^{k,\a}$-regular,  pointed at $0$} (Definition \ref{quasi-vector}). \par
 \subsection{Abstract fundamental pairs and associated curves of exhaustions}
 Motivated by the correspondence of the previous section between  curves of Monge-Amp\`ere quasi-diffeomorphisms and families of pairs $(J_t, X_t)$
 of complex structures and quasi-regular vector fields, we now introduce the following
\begin{definit} {\rm Let $((\OB, J), \t_o)$ be a manifold in normal form of class $\cC^{k,\a}$, $k\geq 2$, $\a \in (0,1)$, 
and a one-parameter family of vectors   $v_t  \in \bC^n \setminus \{0\}$, of class $\cC^{1,  \b}$, $\b \in (0,1)$ in the parameter $t \in [0,1]$. We call
 {\it  abstract fundamental pair guided by    $v_t$} a pair of isotopies $(J_t, X_t)$, of class at least $\cC^{1,\b}$, $\b >0$,  in the parameter $t \in [0,1]$,  where 
 \begin{itemize}
 \item[1)]  $J_t$ is an isotopy of  complex structure of class $\cC^{k, \a}$ of a manifold in normal form $((\OB, J_t), \t_o)$,  
 \item[2)] $X_t$ is an isotopy of  quasi-regular vector fields   $X_t$ on $\OB$  of  class $\cC^{k,\a}$, pointed at  $0$,  with $X_t|_0 = v_t$, 
 \end{itemize}
  satisfying  the  differential condition  
\beq 
\label{42}
 \frac{d J_t}{dt} = - \cL_{X_t} J_t\qquad \text{with initial condition}\ J|_0 = J\ .
\eeq }
 \end{definit}
\  \par
 \medskip
From the discussion of \S \ref{section32}, it is clear that the fundamental pair $(J_{t}, X_{t})$ of a curve $\Phi_t$  of Monge-Amp\`ere quasi-diffeomorphisms
satisfies all conditions of the above definition. Our main goal now is to show that  a converse  is also true, namely  that  each abstract   fundamental pair satisfying appropriate conditions  is the fundamental pair of   a 
curve of Monge-Amp\`ere quasi-diffeomorphisms  and generates a one-parameter family  of Monge-Amp\`ere exhaustions, with
 centers along a fixed straight line. \par
\medskip
The  proof of this crucial result requires a few preliminaries. 
For each   $x \in \bB^n$  we fix once and for all  an automorphism $\bF_x \in \Aut(\OB)$ of the standard closed  ball $(\OB, \Jst)$, mapping $x$ into $0$.  We  select such a distinguished automorphism $\bF_x$ in such a way  that $\bF_0 = \Id_{\OB}$ and that $\bF_x$   depends smoothly on  the point $x$. 
For each vector $0 \neq v \in \bB^n$,  we consider the points $x_t$ of the parameterized straight line $x_t = t v$, 
and for each $t \in [0,1]$, we set
\beq \label{lallalla} \bF_t \=  \bF_{ tv}\qquad \text{and}\qquad \gv_t \= \frac{d \bF_{t +s}(x_{t})}{ds}\bigg|_{s = 0}   \in T_0 \bB^n\  .\eeq
Finally, for each $t \in [0,1]$,   we consider the blow up   $\pi_t:  \TOB_t  \to \OB$ of $\OB$ at the point $x_t$ 
and the  lifted map $\wt \bF_{t}:  \TOB_t \to \TOB$  
of the automorphism  $\bF_{t}$ between $\TOB_t$ and the blow up $\TOB$ of the unit ball at the origin. By  assumptions, each map $\wt \bF_t$ transforms the exceptional divisor of $ \TOB_t $ (which projects onto $x_{t}$) onto the exceptional divisor of  $\TOB$ (which projects onto $0$).  All blow-ups $ \TOB_t $  of the complex manifold $(\OB, J)$ are diffeomorphic one to the other 
and, in particular, to $ \TOB_{t = 0} = \TOB$. On the other hand, the complex structure of $ \TOB$ is of class $\cC^{k,\a}$ with respect to 
the standard coordinates of $\TOB$ and this implies that the complex coordinates of  $ \TOB$ are of class $\cC^{k+1, \a}$ with respect to the standard one.
Thus each $ \TOB_t $ is 
 $\cC^{k+1,\a}$-diffeomorphic to the standard  
blow up $\TOB$ and, consequently, that each map  $\bF_{t}$ is a  quasi-regular diffeomorphisms of class $\cC^{k+1,\a}$ for the complex manifold $(\OB, J)$, pointed at $x_t$.\par
\medskip
 This settled,
we may now state the following  
important lemma. \par
\begin{lem} \label{lemma3.8} Let $0 \neq v \in \bB^n$ and $\gv_t \in \bC^n$, $t \in [0,1]$,  the  one-parameter family of vectors,  defined in \eqref{lallalla}.  
If $(J_t, X_t)$ is an abstract fundamental 
pair of class $\cC^{k,\a}$, $k\geq 2$, $\a \in (0,1)$,  of $((\OB, J), \t_o)$, with  
$X_t|_0 = \gv_t$ for all $t$, 
 then there exists a curve of Monge-Amp\`ere quasi-diffeomorphisms $\Phi_t$ of class $\cC^{k,\a}$, guided by  the curve  $x_t = t v$, $t \in [0,1]$, of which $(J_{t}, X_{t})$ is the associated fundamental pair. 
\end{lem}
\begin{pf}  We   observe that
any  curve of Monge-Amp\`ere quasi-diffeomorphisms  of class $\cC^{k, \a}$, guided by  the straight curve  $x_t = t v$, $t \in [0,1]$, can be expressed as a composition of the form $\Phi_t =\Phi^o_t \circ \bF_{t}$,
where $\Phi^o_t$ is a quasi-regular  $\cC^{k, \a}$ diffeomorphism, pointed at $0$ and mapping $0$ into itself. 
Hence,  the curve of quasi-diffeomorphisms, of which we need to prove the existence, corresponds to a one-parameter family 
of tame diffeomorphisms of the blow-up at the origin $\wt{\Phi^o_t}: \TOB\longrightarrow \TOB$,
whose associated quasi-regular maps $\Phi^o_t$  are solutions  to the problem
\begin{align}
\label{312} &   \left.\frac{d  \Phi^o_{s} \circ \bF_{s}(x)}{ds} \right|_{s = t} = X_t |_{\Phi^o_t(\bF_t(x))}\ , \ x \in \OB\setminus\{0\}\ ,\quad \text{with}\ \Phi^o_0 = \Id
\ ,\\[10pt]
\label{313} & (\Phi^o_t \circ \bF_t)_*(J) = J_t\ ,\\[10pt]
\label{314} &  \left.\frac{d  \Phi^o_s \circ \bF_t(x_{t})}{ds} \right|_{s = 0} = \gv_t\ .
\end{align}
Being $X_t$  quasi-regular with  $X_t|_0 = \gv_t$, by continuity,   each one-parameter family of quasi-regular maps $\Phi^o_t$ that solve   \eqref{312}  necessarily satisfies also  \eqref{314}, proving that the latter is redundant. Moreover,  if $\Phi^o_t$ is a  solution to \eqref{312}, then,   by Proposition \ref{easylemma},  the fundamental pair  of the curve of quasi-diffeomorphisms $\Phi^o_t \circ \bF_t$ is a solution to the same differential problem satisfied by $J_t$ with the same initial condition. By uniqueness of such solution,   we get that  also \eqref{313} is necessarily satisfied. In conclusion  we  only need to show the existence of a one-parameter family of $\cC^{k,\a}$ quasi-diffeomorphisms   satisfying the differential problem  \eqref{312}. \par
\smallskip
In turn, this is equivalent to the problem on $(\OB \setminus  \{0\}) \times [0,1]$,
\beq\label{312bis}  \left.\frac{d  \Phi^o_{s} (y)}{ds} \right|_{s = t} +  \frac{d \bF_s(\bF^{-1}_t(y))}{ds}\bigg|_{s = t} \cdot \Phi^o_t \bigg|_y= X_t |_{\Phi^o_t(y)}\ ,\qquad  \Phi^o_0 = \Id\eeq
(here and below, for simplifying some formulas, we sometimes use the notation  $V{\cdot}f$ for   the directional derivative $V(f)$ of a function $f$ along a vector field $V$). 
Let the maps $\Phi^o_t$ be the values at $t$ of  a  map 
$\Phi^o: [0,1] \times \OB \setminus  \{0\} \longrightarrow \OB $
and $V$,  $W$  denote  the vector fields 
$$V|_{(t, y)}  \= \frac{\p}{\p t} +   \frac{d \bF_s(\bF^{-1}_t(y))}{ds}\bigg|_{s = t}\ ,\qquad W_{(t, y)} \= X_t |_{\Phi^o_t(y)}$$ 
taking values in 
 $T([0,1] \times \OB \setminus  \{0\}) = \bR + T( \OB \setminus  \{0\})$ and  $T( \OB \setminus  \{0\})$, 
respectively. In this way    \eqref{312bis}  can be equivalently formulated as 
\beq V {\cdot} \Phi^o |_{(t, y)} - W|_{(t,y)} = 0\ , \qquad \quad \Phi^o(0, y) = y\ .\eeq
Now, a map $\Phi^o: [0,1] \times \OB \setminus  \{0\} \longrightarrow \OB $ is  a solution to this problem if and only if its  graph 
$$S = \{\ (t, y, x) \ :\ x = \Phi^o(t, y)\ \} \subset  \left([0,1] \times \OB \setminus  \{0\}  \right) \times \OB  $$
is a $(2n+1)$-dimensional real manifold, tangent at all points to the vector field
$$T_{(t, y, x)}\= \frac{\p}{\p t} + W^i|_{(t,y)}\frac{\p}{\p y^i} + \left(V {\cdot} \Phi^o{}^j |_{(t, y)}\right)\frac{\p}{\p x^j}\ .$$
 The required graph $S$ can be directly obtained by taking   the union  
 $$S = \bigcup_{s \in [0,1]} \Psi^T_s(S_0)\ ,\qquad \Psi^T_s = \text{flow of} \ T\ $$
 of the  images $\Psi^T_s(S_0)$,  under the diffeomorphisms of the flow of $T$,  
 of  the transversal $2n$-submanifold
$$S_0 = \{ t, y, x) \in ([0,1] \times \OB \setminus  \{0\}) \times \OB\ :\ t = 0\ ,\ x = y\ \}\ .$$
 By classical properties of  flows, each  diffeomorphism $\Psi^T_s$ is of class $\cC^{k, \a}$, so that $S$ is the graph of the unique solution of \eqref{312bis} of class $\cC^{k, \a}$ on  $\OB \setminus  \{0\}$.\par
 \smallskip  By the regularity assumptions on the maps $\bF_t$ and  the limiting behavior of $X_t$ at $x \to  0$, we see that each  map $\Phi^o_t = \Phi^o(t,\cdot)$  has a  unique  extension to a  map  $\wt{\Phi^o_t}: \TOB \to \TOB$, mapping the exceptional divisor into itself.  On the other hand, since 
  $\Phi_t = \Phi^o_t \circ \bF_t$ is a $(J, J_t)$-holomorphic map of $\bB^n \setminus \{0\}$,   the unique continuous extension at $x_t$ is necessarily $(J, J_t)$-holomorphic also at $x_t$,  i.e. it  is a  holomorphic maps in the usual  sense  in $J$- and $J_t$-holomorphic coordinates.   
  The regularity assumptions on the lifts of  $J$,  $J_t$  at the levels of the blow ups, imply
  that the continuous extension of $\wt \Phi^o_t$ is of class $\cC^{k,\a}$ on 
   $\TOB$ for each $t$. \par
\medskip
In order to conclude the proof, it remains to show that each map $\wt \Phi_t = \wt \Phi^o_t \circ \wt \bF_t$ is a $\cC^{k,\a}$-diffeomorphism or, equivalently, 
that each Jacobian $J (\wt \Phi_t)|_y$ is invertible  at all  $y \in \TOB$.  In fact,  this would directly imply that the induced maps $\Phi_t: \OB \to \OB$ constitute a one-parameter family of Monge-Amp\`ere quasi-diffeomorphism of class $\cC^{k,\a}$, guided by the curve $x_t = tv$.\par
\smallskip
Note that, having proven that each  $\Phi_t(x)$ is well defined and smooth at each point $(t, y) \in [0, 1] \times \TOB$, all Jacobians $J(\Phi_t)|_y$ are bounded for each such $(t, y)$. Hence, in order to prove what we need, we simply have to show that all these Jacobians have non-vanishing determinants. Since this is surely true for the points $(t, y) \in \{0\} \times \TOB$ and  thus, by compactness of $\TOB$,  for all points of $[0, \ve) \times \TOB$,  with $\ve$ sufficiently small, the 
claim   is proven if  
the set 
$$A \= \{ t \in [0, 1]\ : \ \text{for each}\ s \in [0, t),\ J (\wt \Phi^o_s)|_y\ \text{is invertible at each}\ y \in \TOB\ \}$$
coincides with $[0,1]$. By the above remark  and the semi-continuity of the rank,  $A$ is  non-empty and open and it remains to 
prove  that $A$ is closed, i.e. that  for each   sequence $t_k \in A$  converging  to $t_\infty \= \sup A$, one has that 
  $\lim_{k \to \infty} \left|\det J(\wt \Phi_{t_k})|_y\right| > 0$ for all $y \in \TOB$ or, equivalently,  
\beq \label{toprove} \lim_{k \to \infty} \left|\det J(\wt \Phi^{-1}_{t_k})|_{y}\right| < \infty \qquad \text{for each} \ y \in \TOB\ .\eeq
 Since each quasi-diffeomorphism $\wt \Phi_{t_k} = \wt \Phi^o_{t_k} \circ \wt \bF_{t_k}$, $t \in A$,  satisfies \eqref{313},    for  each  straight disk  $\overline{\D^{(c)}} = \{\ \z c\ , \z \in \overline \D\} \subset \OB$, $c = (c^a) \neq 0$, the restriction $\wt \Phi^{-1}_{t_k}|_{\overline {\D^{(c)}}}$ maps properly and   biholomorphically  $\overline{\D^{(c)}}$  into $\OB$. Therefore,   on  each such disc, the  sequence of holomorphic maps $\wt \Phi^{-1}_{t_k}|_{\overline {\D^{(c)}}}$ converges uniformly
to a proper holomorphic map $f_\infty:  \overline{\D^{(c)}} \to \TOB$ passing  through the  point $\wt \Phi_{t_\infty}^{-1}(0) = t_\infty v  \in \bB^n$. This implies that for each $y \in \overline{\D^{(c)}}$, the restrictions of the push-forward  $\wt \Phi^{-1}_{t_k*}|_y$ to the tangent space $T_y \D^{(c)}$ have uniformly bounded components. Since the tangent spaces  $T_y \D^{(c)}$  coincides with the spaces $\cZ_y$ of the distribution $\cZ$, we conclude that 
the restrictions   $\wt \Phi^{-1}_{t_k*}|_{\cZ}$  are uniformly bounded and that  the proof of \eqref{toprove} reduces to check that 
  the restricted linear map  $\wt \Phi^{-1}_{t_k*}|_{\cH_y}$, $y \in \TOB$,  have uniformly bounded components.  \par
\smallskip
To prove this, we  first  focus on the limit behavior of the linear maps $\wt \Phi^{-1}_{t_k*}|_{\cH_y}$ at the points $y \in \p \bB^n$.   The  exhaustion $\t_\infty = \t_o \circ \Phi_{t_\infty}$ is  a regular defining functions for $(\OB, J)$, i.e.  with $d \t_\infty \neq 0$ at all  points of $\p \bB^n$.  This can be checked by simply noticing  that  the derivatives of $\t_\infty$ at  boundary  points  along   vectors  that are tangent to the limit disks  $\overline{\D^{(c)\infty}} \= \lim_{k \to \infty} \wt \Phi^{-1}_{t_k} (\overline {\D^{(c)}})$  (which are transversal to $\p \bB^n$ by Hopf Lemma) are equal to the derivatives
of $\t_o$ at boundary points along vectors that are tangent to the radial disks $\overline{\D(c)}$  and are therefore non vanishing. Being $\t_\infty$ a regular defining function,  the  Levi forms $d d^c_J \t_\infty(\cdot, J \cdot) |_{\cH_y \times \cH_y}$ and $d d^c_J \t_\infty(\cdot, J \cdot)|_{\cH_y \times \cH_y}$  at  each $y \in \p \bB^n$  are    multiples one  of the other  by  the values at $y$ of a nowhere vanishing continuous function. This  implies that  
 the  limit (and bounded) $2$-form 
 $$d d^c_J \t_\infty(\cdot, J \cdot)|_{\p \bB^n}  = \lim_{k \to \infty} \left.\Phi^*_{t_k} \left(d d^c_{J_{t_k}} \t_o(\cdot, J_{t_k} \cdot)\right) \right|_{\p \bB^n} =  \left.\Phi^*_{t_\infty} \left(d d^c_{J_{t_\infty}} \t_o(\cdot, J_{t_\infty} \cdot)\right) \right|_{\p \bB^n} $$  is strictly positive  at all points of $\p \bB^n$ and that $\ker \Phi_{t_\infty*}|_{\cH_y} = \{0\}$  for each $y \in \p \bB^n$. Thus that  the linear maps $\wt \Phi^{-1}_{t_k*}|_{\cH_y}$, $y \in \p \bB^n$, have  uniformly  bounded limits and, consequently,   uniformly  bounded components.\par
 We  now   show that indeed the coordinate components of all linear maps $\wt \Phi^{-1}_{t_k*}|_{\cH_y}$, $y \in \TOB$, 
  are   uniformly bounded. 
Consider at first  a fixed closed radial disk $\overline{\D^{(c_o)}}$, the associated lifted disk $\wt{\overline{\D^{(c_o)}}}$ in $\TOB$, the  corresponding  limit disk $\overline{\D^{(c_o)\infty}} \= \lim_{k \to \infty} \wt \Phi^{-1}_{t_k} (\overline {\D^{(c_o)}})$ in $\OB$  and, finally,  the corresponding lifted closed disk  $\wt{\overline{\D^{(c_o)\infty}}}$  in the blow up $\TOB$. Since the   lifted disks $\wt{\overline{\D^{(c)}}}$ form a regular foliation of  $\TOB$,   we may consider  one-parameter families of maps 
$f_i^{(s)}: \D \to \TOB$, $i = 1, \ldots, n-1$,  with the following properties:
\begin{itemize}
\item[i)] Each $f^{(s)}_i$ is a $J$-holomorphic parametrization of a radial disc $\wt{\overline{\D^{(c(s))}}}$ (recall  that any such map is also  $J_t$ holomorphic for each $t \in [0,1]$, by the properties  of the $L$-complex structures $J_t$); 
\item[ii)] The map $f^{(0)}$ is the standard holomorphic parametrization of $\wt{\overline{\D^{(c_o)}}}$; 
\item[iii)]  The vector fields $Y_i$,   defined at the points of the disk $\D^{(c_o)} = f^{(0)}(\D)$,  by the formula
$$Y_i(\z)\= \frac{d f^{(s)}_i(\z)}{ds}\bigg|_{s = 0}\ ,$$
 have $J_{t_\infty}$-holomorphic components $Y^{(10) \infty}_i$ that are linearly independent and generate the $J_{t_\infty}$-holomorphic  space $\cH^{10}_{J_{t_\infty}}|_{y(\z)} \subset \cH^\bC|_{y(\z)}$ at all  points of the disk $f^{(0)}_i (\z)\in \D^{(c_o)}$.  
 \end{itemize}
 Note that, since the complex structures $J_{t_k}$ converge to the complex structure $J_{t_\infty}$ for $k \to \infty$,  condition (iii)  implies also  that, for all sufficiently large $k$, 
 the $J_{t_k}$-holomorphic components  $Y^{(10)k}_i \in \cH^{10}_{J_{t_k}}$ of the vector fields $Y_i$ are linearly independent and generate  the spaces $\cH^{10}_{J_{t_k}}|_y$ at each point $y \in \D^{(c)}$. \par
 \smallskip
By construction,   the components in  $J$-holomorphic coordinates  of the restricted linear maps  $\wt \Phi^{-1}_{t_k*}|_{\cH_y}$ at the points of the disk $y \in \Phi^{-1}_{t_k}(\D^{(c_o)})$,   are uniformly bounded  if and only if the sequences of vector fields 
$\wh Y^{(10)k}_{i} \= \wt \Phi^{-1}_{t_k*}(Y^{(10)k}_i)$,  defined at the points of the disks $\wt \Phi^{-1}_{t_k} (\wt{\overline {\D^{(c_o)}}})$, 
have  uniformly bounded components at all points of their domains.  To prove this claim,  we first  observe that the  components of these vector fields in   $J$-holomorphic coordinates  are necessarily holomorphic functions. Indeed, the components of the vector fields $Y^{(10)k}_i$, which take values  at the points of $\D^{(c_o)}$,  are actually  derivatives with respect to the parameter    $s$  of   the expressions in $J_{t_k}$-holomorphic coordinates of  the $J_{t_k}$-holomorphic disks $f^{(s)}_i$ of $(\TOB, J_{t_k})$, which are in fact $J_t$-holomorphic disks for any $t$, as we remarked.   
Since $J = \wt \Phi^{-1}_{t_k*}(J_{t_k})$, this shows that
that  the components in $J$-holomorphic coordinates  of the vector fields $\wh Y^{(10)k}_{i} = \wt \Phi^{-1}_{t_k*}(Y^{(10)k}_i)$  are     derivatives in   $s$  of  coordinated expressions of the   $J$-holomorphic disks $\wt \Phi^{-1}_{t_k} \circ f^{(s)}_i$.   From this  the   holomorphicity  of the components of the vector fields $\wh Y^{(10)k}_{i}$    follows immediately.\par
\smallskip
  Observe that  the limits  for $k \to \infty $ of the restricted vector fields   $\wh Y^{(10)k}_{i}|_{\p \bB^n}$ coincide with  the vector fields $Y^{(10)\infty}_i|_{\p \bB^n} = \lim_{k\to \infty} \Phi^{-1}_{t_k*}|_{\p \bB^n} (Y^{(10) k}_i|_{\p \D})$, which are well-defined and bounded because of  previous discussion on the limit behavior of the    $ \Phi^{-1}_{t_k*}|_y$, $y \in \p \bB^n$. This remark guarantees the existence of  uniform upper bounds for the components of the restricted vector fields  $\wh Y^{(10)k}_{i}|_{\p \bB^n}$.  The  holomorphic dependence on $\z \in \D$ of the components of the vector fields $\wh Y^{(10)k}_{i}$  implies  the existence of uniform  bounds  for such components  at {\it all} points of the disk  $\Phi^{-1}_{t_k}(\D^{(c_o)})$, proving the desired claim. 
\par
 Since the above arguments  imply that the derived upper bounds depend continuously on the parameter  $c_o$ of the radial disk $\D^{(c_o)}$,   
 the compactness of $\TOB$ yields the existence of 
a  uniform bound for   all components in $J$-holomorphic coordinates of all linear maps  $\wt \Phi^{-1}_{t_k*}|_{\cH_y}$,  $y \in \TOB$. This concludes the proof.  
  \end{pf}
\begin{rem} \label{whatweneed} {\rm Proposition \ref{theorem43} and Lemma \ref{lemma3.8}      
imply  that if there is  an abstract $\cC^{k,\a}$ fundamental pair $(J_t, X_t)$   with  $X_t|_0 = \gv_t$, $k \geq 2$, $\a \in (0,1)$,   then there is also a one-parameter family of Monge-Amp\`ere exhaustions $\t_t$ of class $\cC^{k,\a}$, whose centers  are at the  points  of the curve $x_t = tv$, $t \in [0,1]$.  Due to this,  the proof of Theorem \ref{main} reduces   to determine the exact set of  vectors   $0 \neq v \in \bB^n$ for which there exists an associated  fundamental pair  $(J_t, X_t)$   with  $X_t|_0 = t v$.}
\end{rem}
\par
\section{Special  vector fields of  manifolds  in normal form}
\label{section4}
\setcounter{equation}{0}
The main purpose of  this section is to prove
that  any  vector field $X_t$,  
in an abstract fundamental pair $(J_t, X_t)$, is subjected to very strong constraints, 
forcing $X_t$ to be in a very  small class of quasi-regular vector fields,  called   {\it special vector fields}. 
 We shall need convenient  complex coordinates for  $\bC^n$. \par
\subsection{Generalized polar coordinates}
\label{generalizedpolar}
Let 
 $\xi:  \bC P^{n-1} \setminus \{z^n = 0\} \to \bC^{n-1}$  be  the  standard affine coordinates
$$\xi([z^1:\dots:z^{n}]) \= \bigg(w^1 := \frac{z^1}{z^n}, \ldots, w^{n-1}:= \frac{z^{n-1}}{z^n}\bigg)$$
and  define $\f: \bC^{n-1} \times S^1\longrightarrow S^{2n-1} \setminus\{ z^n = 0\}$ by
$$
\f(w^1, \ldots, w^{n-1}, e^{i \theta}) = e^{i \theta} \frac{1}{\sqrt{ \textstyle \sum_{\ell = 1}^{n-1} |w^\ell |^2 +1}}\left( w^1, \ldots, w^{n-1}, 1\right).
$$
This is a (real) diffeomorphism onto $S^{2n-1} \setminus\{ z^n = 0\}$ with the  useful  property that,  for each $(w^i) \in \bC^n$ and $\q \in \bR$, the corresponding point  
$\f(w^i, e^{i\theta})$ is   in    the complex line $[w^1: \dots: w^{n-1}: 1] = \xi^{-1}(w^1, \ldots, w^{n-1})$. 
If $\|z\|$ denotes the euclidean norm of  $z$ and $\h: \bC^n \setminus \{z^n = 0\} \longrightarrow \bC^n$ is   defined by  
\begin{equation}
 \label{coordtransf} \h(z^1, \ldots, z^n) :=  \left( w^1:= \frac{z^1}{z^n}, \ldots, w^{n-1}:= \frac{z^{n-1}}{z^n}, \z:= \|z\| \frac{z^n}{|z^n|}\right),
 \end{equation}
 one has 
 $\h|_{ \bC^{n-1} \times S^1} = \f^{-1}$.  This implies that $\h$ is  a (real) diffeomorphism from  $\bC^n \setminus \{z^n = 0\} $ onto its image, whose inverse map is given by 
$$\h^{-1}(w^i, \z) :=   \frac{\z}{\sqrt{  \textstyle \sum_{\ell = 1}^{n-1} |w^\ell |^2 +1}}\left( w^1, \ldots, w^{n-1}, 1\right)\ .$$
We call $\h$  a {\it generalized system of polar coordinates on $\bC^n \setminus \{z^n = 0\}$}.   Note that a similar  construction can be performed if we replace $\{z^n = 0\}$  by  any other   affine hyperplane $\pi \subset \bC^n$  through the origin.  The corresponding system of   (generalized) polar coordinates on $\bC^n \setminus \pi$ will be denoted by $\h_{(\pi)}$. In case $\pi$ is a coordinate hyperplane $\pi_i \= \{z^i = 0\}$,   we  simply denote it by  $\h_{(i)}$.\par
\smallskip
For  later purposes,  it is convenient to have  the  explicit  expressions of the coordinate  vector fields $\z\frac{\p}{\partial\z} $ and $\frac{\p}{\p w^\a}$, $1 \leq \a \leq n-1$,   
of the  polar coordinates $\h = (w, \z)$,   in terms of the standard coordinate vector fields $\frac{\p}{\p z^j}$, $\frac{\p}{\p \bar z^j}$  of $\bC^n$. Using   just the definitions, one can directly see that 
\beq
\begin{split} \label{incoord}
 \z \frac{\partial}{\partial \z} &= z^i \frac{\partial}{\partial z^i}\ ,\\
 \frac{\p}{\p w^\a} &= z^n \frac{\partial}{\partial z^\a}- \frac{1}{2} \frac{\overline z^\a z^n}{\|z\|^2} \left(z^i \frac{\p}{\p z^i} + \overline z^i \frac{\p}{\p \overline z^i}\right) =\\
& =  z^n \frac{\partial}{\partial z^\a}- \frac{1}{2} \frac{\overline w^\a}{1 + \sum_{\b = 1}^{n-1}|w^\b|^2} \left(\z \frac{\p}{\p \z} + \overline \z\frac{\p}{\p \overline \z}\right)\ .
\end{split}
 \eeq
\subsection{Adapted  polar  frame fields} Consider a fixed system of polar coordinates, say $\h_{(n)} = (\z, w^\a)$,  and a circular domain 
$D \subset \bC^n$ of class $\cC^{k, \a}$ determined by the  Minkowski functional $\mu = \mu_D$.  By  the  homogeneity property  of  the Minkowski functionals, the expression of $\mu$ in polar coordinates  has the form 
$\mu( \z, w^\a,\bar \z, \overline{w^\a}) = |\z| \r(w^\a, \overline{w^\a})$ for some  $\r > 0$ of  class $\cC^{k, \a}$. Associated with the considered  polar coordinates and the circular domain $D$, there is  the set $(Z, e_\a)$ of complex vector fields  defined by 
\beq\label{vectorfields}
\begin{split} 
&Z := \z \frac{\p}{\p \z}\ ,\\
&e_\a  = e^{(\r)}_\a := \frac{\p}{\p w^\a} -    \frac{\partial \r^2}{\partial w^\a} \z \frac{\p}{\p \z} + \frac{1}{2} \frac{\overline{w^\a}}{1 + \sum_{\b = 1}^{n-1}|w^\b|^2}  \left(\bar \z \frac{\p}{\p \bar \z} - \z \frac{\p}{\p \z}\right) = \\
&= \frac{\p}{\p w^\a} -    \frac{\partial \log \r^2}{\partial w^\a} \z \frac{\p}{\p \z} + \frac{1}{2} \frac{\p \log\left(1 +  \sum_{\b = 1}^{n-1}|w^\b|^2\right)}{\p w^\a}  \left(\bar \z \frac{\p}{\p \bar \z} - \z \frac{\p}{\p \z}\right) = \\
 & \qquad  =  z^n \frac{\partial}{\partial z^\a} -\bigg( \frac{\partial \r^2}{\partial w^\a} +\frac{\overline{w^\a}}{1 + \sum_{\b = 1}^{n-1}|w^\b|^2} \bigg) z^i \frac{\partial}{\partial z^i}\ .
\end{split}
\eeq
We call  it {\it adapted polar frame field, associated with the circular domain $D$ and the  coordinates $(\z, w^\a)$}. 
The fields of these frames
satisfy the equations
\beq Z(\mu^2) = \mu^2\ ,\qquad e_\a(\mu^2) =  e_\a(|\z|^2 \rho^2) =  0\ ,\qquad \Jst e_\a = i e_\a\ ,\eeq
 where we denote by $\Jst$ the standard complex structure of $\bC^n$.
This means that $\Re(Z)$ and $\Im(Z)$ are  generators for the  distribution $\cZ^{\Jst, \mu^2}$, defined in  \eqref{distquater}, of the manifold of circular type $(\overline D, \Jst, \t = \mu^2)$   and that the vector fields $(\Re(e_\a), \Im(e_\a))$ are generators for the $\Jst$-invariant distribution   $\cH^{\Jst, \mu^2}$, which is   $dd^c_{\Jst} \t$-orthogonal and a complementary distribution of  $\cZ^{\Jst, \mu^2}$. We recall that a normalizing map $\Phi: (\overline D, \Jst) \to (\OB, J)$  maps the distributions $\cZ^{\Jst, \mu^2}$ and $\cH^{\Jst, \mu^2}$ onto the distributions $\cZ$ and $\cH$ of $\OB$ discussed  in \S 2.3. \par
It is useful to  explicitly express  the Lie brackets between such  generators:
\begin{multline} \label{1.4} [Z, e_\a] = 0, \ \   [Z, \overline{e_\a}] = 0,\ \  [e_\a, e_\b] = 0,\ \ 
[e_\a, \overline{e_\b}] =  g_{\a \bar \b} \left(\z \frac{\p}{\p  \z} - \bar \z \frac{\p}{\p \bar \z}\right)\ ,
\end{multline}
where
\begin{equation} \label{att}
\begin{split} g_{\a \bar \b}&:=  \frac{\d_{\a\b} (1 + \sum_{\g = 1}^{n-1}|w^\g|^2)- w^\b \overline{w^\a} }{(1 + \sum_{\g = 1}^{n-1}|w^\g|^2)^2} +  \frac{\partial^2\log  \r^2}{\partial w^\a \partial \overline{w^\b}} = \\
& = \frac{\p^2}{\p w^\a \p w^{\bar \b}} \left(\log \r^2 +   \log\left(1 + \textstyle \sum_{\g = 1}^{n-1}|w^\g|^2\right)\right)\ .
\end{split}
 \end{equation}
The coefficients \eqref{att} are strictly related with the components of the $2$-form $d d^c_{\Jst} \t$, with $\Jst$ standard complex structure of $\bC^n$. In fact,  from \eqref{1.4},  
\beq dd^c_{\Jst} \t(e_\a,  e_{\bar \b}) =  \Jst [e_\a, e_{\bar \b}](\t) =  2 i   |\z|^2 \r^2 \, g_{\a \bar \b} = 2 i \t g_{\a \bar \b} \ .\eeq
\subsection{Special  vector fields of a manifold in normal form}
We now want to characterize   the quasi-regular vector fields $X$,   pointed at $0$,  of a manifold in normal form $((\OB, J), \t_o)$, whose (local) flows
preserve  the restriction $J|_{\cZ}$ of $J$ to  the distribution $\cZ$.
For stating our  result,   we need to fix some notation. \par
\smallskip
Consider the Kobayashi infinitesimal metric $\k$ at the center $x_o = 0$ of  $((\OB, J), \t_o)$,  the (closed) indicatrix $(\overline \cI, \Jst) \subset T_0 \OB \simeq \bC^n$, determined by $\k$, and the corresponding circular representation $\Psi: \wt{\overline \cI} \to \TOB$ (see \S \ref{section31}). Note that the standard coordinates 
of $\OB$ in general {\it do not overlap} in a $\cC^{k,\a}$ fashion with the atlas of complex manifold  structure of $(\OB, J)$. On the contrary, by the properties of the circular representations, any chart of the form  $\h_{(\pi)} \circ \Psi^{-1}$  with $\h_{(\pi)} = (\z, w^\a)$ generalized  polar coordinates on  $\wt{\overline \cI} \to \bC^n$,  does   overlap in a  $\cC^{k-2, \a}$-way with the atlas of the complex manifold of $(\OB, J)$.  We call this kind of  coordinates {\it adapted polar coordinates} of $((\OB, J), \t_o)$. The   corresponding  polar frame fields $(Z, e_\a = e_\a^{(\r)})$ where $\r$ is such that $\k = |\z|^2 \r(w, \bar w)$, are called {\it associated adapted polar frame fields}. From now on,  the only sets of  coordinates we consider on $(\OB, J)$ are either adapted polar coordinates or   coordinates  of the form $\xi \circ \Psi^{-1}$, with $\xi = (z^i)$ standard coordinates of $\overline \cI \subset \bC^n$. We also denote by $\Jst$ the complex structure on $\TOB$ induced via $\Psi: \wt{\overline \cI} \to \TOB$ from the standard complex structure of $\wt{\overline \cI} \subset  \bC^n$.  {\it This complex structure should not be confused with the classical complex structure $J_o$  of $\OB$}, where  $\OB$ is considered as a domain of $\bC^n$. In fact, the components  in a set of  adapted polar coordinates of the tensor field $J_o$ are possibly not even continuous in $x_o = 0$,  
while the components of $\Jst$ in such coordinates are actually constant.  \par
\medskip
 Each {\it real} vector field on $\OB \setminus\{0\}$ of class $\cC^{k, \a}$, $k \geq 1$, $\a \in (0,1)$, admits a unique expansion  of the form  
  \beq\label{expansion} X = X^0 Z + X^\a e_\a + \overline{X^0} \overline{Z} + \overline{X^\a} \overline{e_{\a}} \eeq
  where $Z$, $e_\a$ are  the complex vector fields \eqref{vectorfields} of some adapted polar  frame field. In  \eqref{expansion},   the components $X^0$, $X^\a$  are functions of class $\cC^{k, \a}$ of  the polar coordinates on the open subset $\{\z \neq 0\}$. \par 
  \begin{lem} \label{lemmaA2} Let  $ v = (v^i)  \neq 0$ and $X$  a quasi-regular vector field,   pointed at $0$,     of class $\cC^{k, \a}$ on $(\OB, J)$ and 
  with $X|_0  \=  \Re(v^i \frac{\p}{\p z^i})$.
 If 
 \beq \label{conditiona} ((\cL_X J)|_{\cZ})^\cZ = 0\ ,\eeq
where $({\cdot})^\cZ: T(\bB^n\setminus \{0\} \to \cZ$ is the projection onto $\cZ$,  
 then, for each domain   $\cU_{(i)}$ of affine coordinates $\xi: \cU_{(i)} =  \bC P^{n-1} \setminus \{z^i = 0\} \to \bC^{n-1}$,  there exist two real functions $\k,\r : \cU_{(i)}  \to \bR$
 of class $\cC^{k, \a}$, the first unconstrained, the second with $\r > 0$ and $d d^c \r >0$ at all points, such that: \\[5pt]
  (1) The  components    $X^\a = X^\a(w, \bar w, \z, \bar \z)$ are necessarily of the form
 \beq \label{goccia} X^\a = \frac{Y^\a(w, \bar w) }{ \z} + \wt Y^\a(w, \bar w,  \z,  \bar \z)\eeq 
  with $Y^\a(w, \bar w)$ defined by 
  \beq \label{C6bis} Y^\a(w, \bar w)  = (v^\a - v^n w^\a) \sqrt{  \textstyle \sum_\g |w^\g|^2 +1}\eeq
    and  $\wt Y^\a = \wt Y^\a(w, \bar w,  \z, \bar \z)$ functions  such that, for each  sequence $x_k \in  \OB \setminus\{0\}$ converging to $0$, the limits  $\lim_{k \to \infty}(\wt Y^\a \z)|_{x_k} $ are  $0$; \\[5pt]
 (2)  The component $X^0$  in \eqref{expansion} has  necessarily  the form
 \beq \label{623} X^0 =  \frac{Y^0(w, \bar w)}{\z} +  i \kappa(w,\bar w) - \overline{Y^0(w, \bar w)}\z\ , \eeq
  where $Y^0$  is  the  complex function on $\bC P^{n-1}$ 
  \beq \label{C7bis} Y^0(w, \bar w) \= v^n \sqrt{  \textstyle \sum_\g |w^\g|^2 +1} + Y^\a \frac{\p \left( \log \r^2 + \log (1 +  \sum|w^\g|^2)\right)}{\p w^\a}\bigg|_{(w, \bar w)} \ .\eeq
 \end{lem} 
  \begin{pf}
Condition \eqref{conditiona}  is actually equivalent to the pair of conditions
\beq\label{413bis} ( \cL_X J(Z))^\cZ = (\cL_X J (\overline Z))^\cZ = 0\ .\eeq
However, since $X$ and $J$ are  real, 
these  equations are one  conjugate to the other, so that they are both satisfied  if and only if just the second one holds. Using the fact that $J|_{\cZ} = \Jst|_{\cZ}$ , this is in turn 
equivalent to 
\beq 
\begin{split}0 & = [X, J \overline Z]^\cZ -   (J[X, \overline Z ])^\cZ =   \\
& = [X, \Jst \overline Z]^\cZ - \Jst([X, \overline Z ]^\cZ)  =  - i  [X,  \overline Z]^\cZ - \Jst [X, \overline Z]^\cZ\ ,
\end{split}
\eeq
meaning that $ [X, \overline Z]^\cZ$ is a complex vector field taking values in  the anti-holomorphic distribution $T^{01}_{\Jst} (\OB \setminus \{0\})$ of the standard complex structure $\Jst$. The Lie bracket $ [X, \overline Z]$ can be easily computed using \eqref{1.4}. One gets that 
$ [X,  \overline Z]^\cZ$ is in $T^{01}_{\Jst}(\OB \setminus \{0\})$ if and only if  
 \beq\label{holomorphicity} \overline Z \cdot X^0 = \bar \z \frac{\p X^0}{\p \bar \z} = 0\ .\eeq
 This yields that, along each straight disk   $\D^{(c)}\= \{w  = c \}$, 
  the complex vector field  $X^{(c)} \= (X^0 Z)|_{\D^{(c)}\setminus \{0\}}$
 is  a vector field of type $(1,0)$,   holomorphic in the polar coordinate $\z$.  Since we are also assuming that $X$ is quasi-regular (hence with a continuous extension at $0$),  it follows that   $X^{(c)} = X^0|_{\D^{(c)}}  \z \frac{\partial}{\partial \z}|_{\D^{(c)}\setminus \{0\}} $ extends  continuously   at the origin and   that the function  $  a \= X^0\big|_{\D^{(c)}\setminus \{0\}}$  has  at most one pole of order $1$ at the origin. 
 On the other hand, since the vector field $X$ is tangent   to the boundary at the points of $\p \bB^n$,   the vector field $X^{(c)}$  is  tangent  to $\p \D^{(c)}$ at all boundary points. This implies that 
 \beq \label{boundarycond} X^{(c)}|_{\partial \D^{(c)}} + \overline{ X^{(c)}|_{\partial \D^{(c)}}}= s \bigg(i Z - i \overline Z\bigg)\bigg|_{\partial \D^{(c)}}\eeq
 for some smooth real function $s: \partial \D^{(c)} \longrightarrow \bR$.   This is equivalent to require that   $a\= X^0\big|_{\D^{(c)}\setminus \{0\}}$ is such that 
$a|_{\partial \D^{(c)}} = - \overline{a|_{\partial \D^{(c)}}}$.  Very standard arguments  show  that  the above two conditions 
imply that 
  $  a = \frac{\l}{\z} + i \k - \overline{\l} \z $
 for some complex number $\l$ and some real number $\k$ which depend only on  $c = (c^\a)$. From this we obtain that  $X^0$ has the form 
\eqref{623} for some appropriate function $Y^0(w, \bar w)$.  
\par
\medskip
 Consider now  a sequence of points  $x_k$   in $\OB \setminus \{0\}$ converging to $0$. With no loss of generality, 
we may assume that all points $x_k$ are in the domain of the  polar coordinates $\h_{(n)} = (w, \z)$ and that 
their expressions in  polar coordinates 
 $x_k =  (w^\a_k, \z_k)$ converge to an $n$-tuple $(w^\a_o, 0)$  with $0 \leq |w_o| <\infty$.   \par
 \medskip
  We now consider the sequences of complex  values 
  $C^\a(x_k) \= (X^\a \z)|_{x_k}$.
   Up to a subsequence, we may assume that the  limits $C^\a = \lim_{k \to \infty} C^\a_k$  exist  even if   they might not  be all finite. Then:\\[10pt]
(i) From  \eqref{vectorfields} and using  \eqref{coordtransf}  to relate   $(z{}^j)$ and  $(w{}^\a, \z)$, we have that  the sequence  of vectors $X^\a e_\b\big|_{x_k}$ converges   to the  vector in $T_0 \bB^n$ (here, some components might be equal to  $\infty$)
 \begin{multline}   \label{limitino} \frac{C^\b}{\sqrt{  \displaystyle \sum |w^\g_o|^2 +1}}\left( \d^\a_\b - w^\a_o\frac{\p\log \r^2}{\p w^\b} \bigg|_{(w_o, \bar w_o)}- \frac{\overline w^\b_o w^\a_o}{1 + \sum|w^\g_o|^2}\right)  \frac{\partial}{\partial z{}^\a}\bigg|_{0} -\\
 -\frac{C^\b}{\sqrt{  \displaystyle \sum |w^\g_o|^2 +1}}\left(\frac{\p \log \r^2}{\p w^\b}\bigg|_{(w_o, \bar w_o)} +  \frac{\overline w^\b_o}{1 + \sum|w^\g_o|^2}  \right) \frac{\partial}{\partial z{}^n}\bigg|_{0} 
  \ .
\end{multline}\\[10pt]
(ii) Using once again \eqref{incoord} and   the fact that  $X^0$ has the form \eqref{623},  the sequence of vectors $ X^0 Z\big|_{x_k} $ converges to the vector in $T_0 \bB^n$
 \beq \label{limitino-bis}   \ \frac{Y^0(w_o, \bar w_o)w^\a_o}{\sqrt{\displaystyle \sum |w^\g_o |^2 +1}}   \frac{\partial}{\partial z{}^\a}\bigg|_{0} +  \frac{Y^0(w_o, \bar w_o)}{\sqrt{\displaystyle \sum |w^\g_o |^2 +1}} \frac{\partial}{\partial z{}^n}\bigg|_{0} \ .\eeq
\\[10pt]
The condition that $X$ is  continuous at $0$   is equivalent to requiring that  for any  sequence $x_k$ as above,  the limit  $\lim_{x_k \to 0}  X|_{x_k} $ exists and does not depend on the choice of the   sequence. From  (i) and (ii)   we infer  that   the vector
$$ \frac{C^\b\d_\b^\a}{ \sqrt{  \displaystyle \sum |w^\g_o|^2 +1}} \frac{\partial}{\partial z{}^\a}\bigg|_{0} + \phantom{aa
aaaaaaaaaaaaaaaaaaa aaaa aaaa aaaa aaaa aaaa}$$
 \begin{multline} \label{limit-ter}
+ \frac{w^\a_o}{\sqrt{  \displaystyle \sum |w^\g_o|^2 +1}}\left(   Y^0(w_o, \bar w_o) - C^\b \frac{\p \log \r^2}{\p w^\b}\bigg|_{(w_o, \bar w_o)}\hskip - 7pt  -\frac{  C^\b\overline w^\b_o}{1 + \sum|w^\g_o|^2} \right)  \frac{\partial}{\partial z{}^\a}\bigg|_{0} +\\
 +\frac{1}{\sqrt{  \displaystyle \sum |w^\g_o|^2 +1}} \left( Y^0(w_o, \bar w_o) - C^\b  \frac{\p \log \r^2}{\p w^\b}\bigg|_{(w_o, \bar w_o)} -\frac{C^\b\overline w^\b_o}{1 + \sum|w^\g_o|^2}\right)   \frac{\partial}{\partial z{}^n}\bigg|_{0} 
\end{multline}
is  independent of  $w_o$. Necessary and sufficient conditions for this are: \\[10pt]
a) there is a constant $v^n$ such that 
 \beq \label{famina}  
 \begin{split} Y^0(w_o, \bar w_o) & - C^\b \frac{\p \log \r^2}{\p w^\b}\bigg|_{(w_o, \bar w_o)}  -\frac{C^\b\overline w^\b_o}{1 + \sum|w^\g_o|^2}  = 
 \\ = Y^0(w_o, \bar w_o) & - C^\b \frac{\p \left(\log  \r^2 + \log (1 +  \sum|w^\g_o|^2)\right)}{\p w^\b}\bigg|_{(w_o, \bar w_o)} = v^n \sqrt{  \displaystyle \sum |w^\g_o|^2 +1}
 \end{split}
 \eeq
and, in particular, all limits $C^\a$ are finite; \\
b)   the vector \eqref{limit-ter}, which can now be   written as
  $$\left( \frac{C^\a}{\sqrt{  \displaystyle \sum |w^\g_o|^2 +1}} + w^\b_o v^n\right) \frac{\partial}{\partial z^\b}\bigg|_{0} 
 +v^n   \frac{\partial}{\partial z^n}\bigg|_{0} \ ,$$
does not depend on $w_o$.\\[10pt]
 This last condition is equivalent to require that  there exist constants $v^\b$ such that 
 $C^\a = (v^\b - v^n w^\b_o) \sqrt{  \displaystyle \sum |w^\g_o|^2 +1} $. Replacing this into \eqref{famina} and  \eqref{limitino},  claims (1) and (2)  follow.
\end{pf}
The above lemma  motivates the following notion, which  enters as crucial ingredient in the proof of our main result. In what follows,  
$\r: \bC P^{n-1} \to \bR$ is  always the positive function,  determined by  the Kobayashi metric $\k_{x = 0}$ of $(\OB, J)$ and gives the Minkowski function of   the indicatrix $\overline{\cI_0} \subset T_0 \OB = \bC^n$.\par
\begin{definit} \label{specialA} {\rm
Consider a vector   $ v  = (v^i) \in \bC^n$,  a real valued function   $\s : \bC P^{n-1} \to \bR$ of class $\cC^{k,\a}$ and 
$n-1$ complex functions $\wt Y^\a: \TOB \to \bC$, which are of class $\cC^{k,\a}$ on $\OB \setminus \{0\}$ and with $\lim_{k \to \infty}  (\wt Y^\a \z)|_{x_k}  = 0$ for any sequence $x_k \to 0$ in  any set of adapted polar coordinates of $((\OB, J), \t_o)$.  We call {\it special vector field of $(\OB, J)$ determined by $v$,  $\s$ and $\wt Y^\a$} the  quasi-regular  vector field of $(\OB, J)$, pointed at $0$, whose    expansion in terms of an  adapted polar  frame field $(Z, e_\a)$ has the form
\beq\label{special}
\begin{split} X^{[\r](v,  \s, \wt Y^\a)} =  \bigg(\frac{Y^0}{\z} + & i \s -  \overline{Y^0} \z  \bigg)Z + \left(\frac{Y^\a}{\z} + \wt Y^\a\right) e_\a 
 +\\
 & + \bigg(\frac{\overline{Y^0}}{\bar \z}- i \s -Y^0 \bar \z  \bigg)\overline{Z} + \left(\frac{\overline{Y^\a}}{\overline \z} + \overline{\wt Y^\a}\right)\overline{e_{\a}}
 \end{split}\eeq
with $Y^0$ and $Y^\a$  as in \eqref{C7bis} and \eqref{C6bis}, respectively.}
\end{definit}
From the proof of Lemma \ref{lemmaA2}, the next   corollary  follows immediately.\par
\begin{cor} \label{corA4} A quasi-regular vector field  of class $\cC^{k, \a}$, pointed at $0$, on $(\OB, J)$,   with  $ X|_0 = v $ and satisfying  \eqref{conditiona} is a special vector field $X = X^{[\r](v,  \s, \wt Y^\a)}$ for some triple $(v,  \s, \wt Y^\a)$. Conversely, for  each choice of  a triple $(v, \s, \wt Y^\a)$, the corresponding special vector field $X^{[\r](v, \s, \wt Y^\a)}$ 
is a quasi-regular vector field  of class $\cC^{k, \a}$ pointed at $0$  on $(\OB, J)$,   with  $ X^{[\r](v,  \s, \wt Y^\a)}|_0 = v $ and satisfying  \eqref{conditiona}. 
\end{cor}
Moreover, we have the following technical remark.
\begin{lem}  \label{specialA-bis} 
In each set of polar coordinates $(\z, w^\a)$,   the functions $Y^0$, $Y^\a$,  defined in  \eqref{C7bis} and \eqref{C6bis},  satisfy the identities  
 \beq 
 \label{4.26}
 \begin{split}
e_{\b}\left(\frac{Y^0}{\z}  - \overline{Y^0} \z\right) &= \frac{H_\b}{\z} -
\overline{\left(H_{\bar \b} - Y^0 e_{\bar \b}(\log \r^2)\right)}\z \ , \\  e_{\bar \b}\left(\frac{Y^0}{\z}  - \overline{Y^0} \z\right) &= \frac{H_{\bar \b}}{\z}-
\overline{\left(H_{\b} - Y^0 e_{\bar \b}(\log \r^2)\right)}\z  \ , 
\end{split}
\eeq
where, if $h_{\a \b} \= \frac{\p \left( \log \r^2 + \log (1 +  \sum|w^\g_o|^2)\right)}{\p w^\a \p w^\b}$,
\beq \label{4.24} 
\begin{split} 
 \hskip-0,5cm H_{\b} &\=Y^\a \bigg(h_{\a \b} + \frac{\p \left(\log \r^2 + \log (1 +  \sum|w^\g|^2)\right)}{\p w^\a}\frac{\p \left( \r^2 + \log (1 +  \sum|w^\g|^2)\right)}{\p w^\b}\bigg)\\
 H_{\bar \b} &\=Y^\a g_{\a \bar \b}
 \end{split} \eeq
\end{lem}
\begin{pf} 
We first observe 
 that 
$$e_\b(Y^\a) =   - v^n \d^\a_\b  \sqrt{\textstyle 1 + \sum_\g|w^\g|^2} +  \frac{Y^\a \overline{w^\b}}{2 ({\textstyle 1 + \sum_\g|w^\g|^2)}}\ , \qquad e_{\bar \b} (Y^\a) = \frac{Y^\a w^\b}{2 ({\textstyle 1 + \sum_\g|w^\g|^2)}}\ .$$
Then
\begin{multline*} e_{\bar \b} \bigg(\frac{Y^0}{\z}\bigg) = \frac{1}{\z}e_{\bar \b} \left( v^n \sqrt{  \textstyle \sum_\g |w^\g|^2 +1} + Y^\a \frac{\p \left( \log r^2 + \log (1 +  \sum|w^\g|^2)\right)}{\p w^\a} \right) - \\
-  \frac{1}{2\z} \frac{\p \log\left(1 +  \sum_\g |w^\g|^2\right)}{\p \overline{w^\b}} \left( v^n \sqrt{  \textstyle \sum_\g |w^\g|^2 +1}\right) - 
\end{multline*}
$$
- \frac{1}{2\z} \frac{\p \log\left(1 +  \sum_\g |w^\g|^2\right)}{\p \overline{w^\b}}  \left( Y^\a \frac{\p \left( \log r^2 + \log (1 +  \sum|w^\g|^2)\right)}{\p w^\a} \right) = 
$$
\begin{multline*}
= \frac{1}{2 \z} v^n \frac{w^\b}{ \sqrt{  \textstyle \sum_\g |w^\g|^2 +1}}  + \frac{1}{\z} \frac{Y^\a w^\b}{2 ({\textstyle 1 + \sum_\g|w^\g|^2)}}\frac{\p \left( \log r^2 + \log (1 +  \sum|w^\g|^2)\right)}{\p w^\a} + \\
+ \frac{Y^\a g_{\a \bar \b}}{\z} - \frac{1}{2 \z} v^n \frac{w^\b}{ \sqrt{  \textstyle \sum_\g |w^\g|^2 +1}} -  \\
- \frac{1}{2 \z} \frac{Y^\a w^\b}{ ({\textstyle 1 + \sum_\g|w^\g|^2)}}\frac{\p \left( \log r^2 + \log (1 +  \sum|w^\g|^2)\right)}{\p w^\a}
= \frac{Y^\a g_{\a \bar \b}}{\z} = \frac{H_{\bar \b}}{\z}\ .
\end{multline*}
Similarly one gets  that  $e_{\b}\left(\frac{Y^0}{\z}\right) = \frac{H_\b}{\z}$.  From this and  $e_\b(\z\bar \z^2 \r^2) {= }e_{\bar \b}(\z\bar \z^2 \r^2) {=} 0$, one can also  derive the expressions for $e_{\b}(\overline {Y^0}\z)$,  $e_{\bar \b}(\overline {Y^0}\z)$ and get \eqref{4.26}. 
\end{pf}
\section{The proof of the Propagation of Regularity Theorem}
\label{section5}
\setcounter{equation}{0}
\subsection{The differential problem  characterizing fundamental pairs}
 In this section,  we show that the abstract  fundamental pairs $(J_t, X_t)$ of 
 a manifold in normal form $((\OB, J), \t_o)$   are precisely the solutions of  
 an appropriate   differential problem.  By Remark \ref{whatweneed} this reduces our main result  
 to  the proof of the existence of solutions to such  problem.
\begin{prop} \label{lemma44} Let  $(J_t, X_t)$ be an abstract fundamental pair of class $\cC^{k,\a}$ on a manifold in normal form $((\OB, J), \t_o)$, guided by  the curve  $v_t = X_t|_0$, $t \in [0, 1]$,  and let $(\z, w^\a)$
be a system of polar coordinates  with  associated   polar frame field $(Z, e_\a)$. 
Then 
 each $X_t$ is a special vector field $X_t = X^{[\r_t](v_t,  \s_t, \wt Y^\a_t)}$,  with   components 
\begin{multline} \label{crucial3} 
\!\!\!\!\! \wt Y^\a_t  {=}   i g^{\a \bar \b} e_{\bar \b} (\s_t)
  - \frac{i}{2}  (J_t - \Jst)^\a_{\bar \d} g^{ \bar  \d\g} \bigg(\frac{H_{t \g}}{\z} + \frac{\overline{H_{t \bar \g}}}{\bar \z} \bigg)  - \frac{i}{2}  (J_t - \Jst)^{\a}_{\d} g^{\d \bar  \g}\bigg(\frac{H_{t \bar \g}}{\z} + \frac{\overline{H_{t  \g}}}{\bar \z} \bigg) + \\
    + \frac{i}{2} (J_t - \Jst)^\a_{\bar \d} g^{\bar \d  \g} \left(
\overline{\left(H_{t \bar \g} - Y^0_t e_{\bar \b}(\log \r^2_t)\right)}\z  + \left(H_{t   \g} - Y^0_t e_{ \g}(\log \r^2_t)\right)\bar \z \right)  + \\
 +  \frac{i}{2}  (J_t - \Jst)^{\a}_{\d} g^{\d \bar  \g} \left(
\overline{\left(H_{t\g} - Y^0_t e_{\g}(\log \r_t^2)\right)}\z  + \left(H_{t\bar \g} - Y^0_t e_{\bar \g}(\log \r_t^2)\right)\bar \z \right)\ , 
 \end{multline}
where   $(g^{\a \overline \b}) $ is the inverse  matrix of   \eqref{att},  $H_{t \b}$,  $H_{t \bar \b}$ are  as in \eqref{4.24}, $(J_t - \Jst)^{\g}_{\bar \b}$,   $(J_t - \Jst)^{\bar \g}_{\bar \b}$ are the components of $(J_t - \Jst)(e_{\bar \b})$ in the frame $(e_\a, e_{\bar \a})$ and $Y^\a_t$, $Y^0_t$  are  as in \eqref{C6bis}, \eqref{C7bis}   with    $v = v_t$ and $\r_t$ determined by the Kobayashi metric at $0$ of $(\OB, J_t)$, as described above.\par
Conversely, if $X_t = X^{[\r_t](v_t, \s_t, \wt Y^\a_t)}$, $t \in [0,1]$,  is a one-parameter family of special vector fields,  with $ \wt Y^\a_t$  as in  \eqref{crucial3}, and if $J_t$ and $\r_t$  are  one-parameter families of complex structures and positive real functions on $\bC P^{n-1}$,   which satisfy    the differential problem given by \eqref{obstruction} and 
\beq \label{eiaeia} \frac{d J_t}{dt} = - \cL_{X_t} J_t\ ,\qquad J_0 = J\ ,\eeq
then $(J_t, X_t)$ is a fundamental pair guided by the curve $v_t$.
\end{prop}
\begin{pf}  Since $(J_t, X_t)$ is a fundamental pair,  the one-parameter family  $J_t$ consists of complex structures  that are pushed-down onto $\OB$  of $L$-complex structures of $\TOB$. This    implies   that 
$$ J_t|_{\cZ} =  \Jst|_{\cZ}\qquad\text{and}\qquad J_t(\cH) \subset \cH\ \quad \text{for each}\ t \in [0, T]\ .$$
Since $J_0 = J$,   this is tantamount  to say that 
\beq\label{6.7} \frac{d J_t}{dt}\bigg|_{\cZ} = 0\ ,\qquad \frac{d J_t}{dt}(\cH) \subset \cH\qquad \text{for all}\ t\ .\eeq
Let us focus on the first of these conditions. Combining it with the property $\cL_{X_t} J_t = - \frac{d J_t}{dt}$, this  is equivalent to 
 $\cL_{X_t} J_t(Z) = \cL_{X_t} J_t (\overline{Z}) = 0$ at all points of $\OB \setminus \{0\}$.
However, being $X_t$ and $J_t$  real, 
these  equations are one  conjugate to the other,  so that we may consider  just the second one. This is in turn 
equivalent to 
\beq \label{5.4} 0 = [X_t, J_t \overline{Z}] - J_t[X_t, \overline{Z} ] =[X_t, \Jst \overline{Z}] - J_t[X_t, \overline{Z} ] = - i  [X_t, \overline{Z}] - J_t  [X_t, \overline{Z}]\ ,\eeq
meaning that 
$ [X_t, \overline{Z}] $
 is a complex vector field taking values in  the $J_t$-anti-holomorphic distribution $T^{01}_{J_t} (\OB \setminus \{0\})$. \par
Consider now  the expansion of the vector fields $X_t$ in terms of a polar frame field $(Z, e_\a)$, associated with a set of polar coordinates
$$X_t = X^0_t Z + X^\a_t e_\a + \overline{X^0_t} \overline{Z} + \overline{X^\a_t} \overline{e_{\a}}\ .$$
Since  $T^{01}_{J_t} (\OB \setminus \{0\})|_x = \cZ^{01}_x \oplus \left(T^{01}_{J_t} (\OB \setminus \{0\})|_x \cap \cH^\bC_x\right)$, we have that  $[X_t, \overline{Z}]$ is in $T^{01}_{J_t} (\OB \setminus \{0\})$ if and only if 
\beq \label{6.15} 
\begin{split}
& \overline{Z}(X^0_t) = 0\ ,\\
& [\overline Z,  X^\a_t e_\a  +\overline{X^\a_t} \overline{e_{\a}} ] = \overline Z(X^\a_t) e_\a  + \overline Z(\overline{X^\a_t} )\overline{e_{\a}} \in T^{01}_{J_t} (\OB \setminus \{0\})\cap \cH^\bC\ .
\end{split}
\eeq
 The first part of the proof of Lemma \ref{lemmaA2} shows that the first condition in \eqref{6.15} is equivalent to the condition $(\cL_{X_t} J_t)^\cZ = 0$. Hence Corollary \ref{corA4} applies and, for each $t$, the vector field $X_t$ is a special vector field $X_t = X^{[\r_t](v_t,  \s_t, \wt Y^\a_t)}$. In particular,  $X^0_t$   is equal to 
\beq \label{623bis} X^0_t =  \frac{Y^0_t(w, \bar w)}{\z}  +  i \s_t(w,\bar w) - \overline{Y^0_t(w, \bar w)}\z  \eeq
with $Y^0_t$  defined in \eqref{C7bis}.\par
\smallskip
Let us now consider the second condition in  \eqref{6.15}. Since  the  $2$-form  $-\frac{i}{2 \t}d d^c_{J_t} \t_o|_{\cH^\bC \times \cH^\bC}$ is non-degenerate and  $J_t$-Hermitian, it  can be equivalently stated saying that for each $t \in [0,1]$, 
$$- \frac{i}{2 \t}d d^c_{J_t} \t_o( [\overline Z,  X^\a_t e_\a  +\overline{X^\a_t} \overline{e_{\a}} ],  F)_x{=} -\frac{i}{2 \t}d d^c_{\Jst} \t_o( [\overline Z,  X^\a_t e_\a  +\overline{X^\a_t} \overline{e_{\a}} ], F)_x = 0$$ 
 for any  complex vector field $ F$  in $\cH^{01}_t = \cH^\bC \cap T^{01}_{J_t} \OB$. Since  $\cL_{\overline Z} (\frac{i}{2 \t}d d^c_{\Jst} \t_o) = 0$,  this is equivalent to 
\beq\label{obstruction} -\frac{i}{2 \t} d d^c_{\Jst} \t_o(X^\a_t e_\a  +\overline{X^\a_t} \overline{e_{\a}} , [\overline Z, F]) = 0\qquad \text{for any}\ F \in \cH^{01}_t\ .\eeq
This relation is actually an identity that  is a consequence of  the  integrability of the  $J_t$,  a claim that can be directly checked  using   deformation tensors of complex structures (\cite{BD,PS}). For reader's convenience,  we give the details of such an argument  at the end of \S \ref{review} below.\par
\smallskip
 Let us now focus  on the   second part of  \eqref{6.7}, which  we did not consider yet. Once again, since $J_t$ and $X_t$ are related by  \eqref{42}, this is the same of 
 requiring that,  for any $E   \in \cH$,
\beq\label{6.15bis} (\cL_{X_t} J_t)(E) = [X_t, J_t E] - J_t[X_t, E] \in \cH\ .\eeq
This is also the same of saying that   for any $E   \in \cH$
\beq \label{417} d\t_o( [X_t, J_t E] ) - d \t_o(J_t[X_t, E])  = 0\ .\eeq
We now observe that 
\beq\label{www}
d\t_o( [X_t, J_t E] ) =  X_t \left(d\t_o(J_t E)\right)  - J_t E \left(d \t_o(X_t)\right) 
= - \t_o \left(J_t E (X^0_t) + J_tE(\overline{X^0_t}))\right) 
\eeq
and 
\begin{multline}\label{zzz}
 - d \t_o(J_t[X_t, E]) =  -  d d^c_{J_t} \t_o (X_t, E) - X_t \left(d\t_o(J_t E)\right)  +  E \left(d \t_o(J_t X_t)\right) = \\
=   - d d^c_{J_t}\t_o(X_t, E) + i \t_o \left(E(X^0_t)  -  E(\overline{X^0_t})\right) = \\
=   - d d^c_{J_t}\t_o(X^0_t Z + \overline{X^0_t}\overline{Z}, E) 
- d d^c_{J_t}\t_o(X^\a_t e_\a + \overline{X^\a_t}\overline{e_\a}, E) + i \t_o E(X^0_t  -  \overline{X^0_t})= \\
=- d d^c_{J_t}\t_o(X^\a_t e_\a + \overline{X^\a_t}\overline{e_\a}, E) + i \t_o E (X^0_t  -  \overline{X^0_t}) = \\
\overset{\text{Lemma \ref{crucial5}}}=  - d d^c_{\Jst}\t_o(X^\a_t e_\a + \overline{X^\a_t}\overline{e_\a}, E)+ i \t_o E(X^0_t -\overline{X^0_t})\end{multline}
From this and \eqref{www},  it follows that \eqref{417} is equivalent to 
\beq 
\label{tremoto}
\begin{split}
 \!\!\! d d^c_{\Jst}\t_o & (X^\a_t e_\a + \overline{X^\a_t}\overline{e_\a}, E) =  i\t_o \left((E+ i J_t E) (X^0_t) {-}(E - i J_t E)( \overline{X^0_t})\right) {=} \\
& =  i\t_o \left(d(X^0_t- \overline{X^0_t})  + i    d(X^0_t+ \overline{X^0_t}) \circ J_t\right) \big|_{\cH} (E) = \\
&= i\t_o \left(d(X^0_t- \overline{X^0_t}) + i d(X^0_t+ \overline{X^0_t}) \circ  \Jst \right) \big|_{\cH} (E)  - \\
& \hskip 4cm   - \t_o d(X^0_t+ \overline{X^0_t}) \circ  (J_ t- \Jst) \big|_{\cH} (E) \ .
  \end{split}
  \eeq
Moreover, denoting by $(d d^c_{\Jst} \t_o)^{-1}$ the $(2,0)$-type tensor field, with components given by the inverse matrix of the components of $d d^c_{\Jst} \t_o$, the  \eqref{tremoto} becomes
\beq\label{tremotobis}
\begin{split}
 \hskip - 0.27 cm X^\a_t e_\a + \overline{X^\a_t}\overline{e_\a} &=  i\t_o (d d^c_{\Jst} \t_o)^{-1} \left(d(X^0_t- \overline{X^0_t}) + i d(X^0_t+ \overline{X^0_t}) \circ  \Jst \right) \big|_{\cH} +\\
 & -    \t_o (d d^c_{\Jst} \t_o)^{-1} \left(d(X^0_t+ \overline{X^0_t}) \circ (J_t - \Jst)\big|_{\cH} \right)\ ,
\end{split}
  \eeq
  so that, by \eqref{att}  and \eqref{4.26},   
  \beq \label{tremototerra-1}
  \begin{split}
X^\a_t  &=  g^{\a \bar  \b}  \overline{e_{\b}}(X^0_t)   - \frac{i}{2} g^{\a \bar  \b} (J_t - \Jst)^\g_{\bar \b} e_{\g}(X^0_t+ \overline{X^0_t})  -
\\
& \hskip 5cm -  
\frac{i}{2} g^{\a \bar  \b}  (J_t - \Jst)^{\bar \g}_{\bar \b}  e_{\bar \g}(X^0_t+ \overline{X^0_t}) 
\\
& = \frac{Y^\a_t}{\z} + i g^{\a \bar \b} e_{\bar \b} (\s_t) 
  - \frac{i}{2} g^{\a \bar  \b} (J_t - \Jst)^\g_{\bar \b}\bigg(\frac{H_{t \g}}{\z} + \frac{\overline{H_{t \bar \g}}}{\bar \z} \bigg) -\\
& \hskip 2 cm  - \frac{i}{2} g^{\a \bar  \b} (J_t - \Jst)^{\bar \g}_{\bar \b}\bigg(\frac{H_{t \bar \g}}{\z} + \frac{\overline{H_{t  \g}}}{\bar \z} \bigg)  -\\
&  + \frac{i}{2} g^{\a \bar  \b} (J_t - \Jst)^\g_{\bar \b}\left(
\overline{\left(H_{t \bar \g} - Y^0_t e_{\bar \b}(\log \r^2)\right)}\z  + \left(H_{t   \g} - Y^0_t e_{ \g}(\log \r^2)\right)\bar \z \right)  + \\
& +  \frac{i}{2} g^{\a \bar  \b} (J_t - \Jst)^{\bar \g}_{\bar \b} \left(
\overline{\left(H_{t\g} - Y^0_t e_{\g}(\log \r^2)\right)}\z  + \left(H_{t\bar \g} - Y^0_t e_{\bar \g}(\log \r^2)\right)\bar \z \right)\ .
\end{split}
\eeq
It is now convenient to recall that, by Lemma \ref{crucial5}, for any $X, Y \in \cH$, 
$$d d^c \t_o(J_t X, Y) = d d^c_{J_t} \t_o(J_t X, Y) = - d d^c_{J_t} \t_o( X, J_t Y) = - d d^c \t_o( X, J_tY)\ .$$
This implies that $(J_t)^\g_\a g_{\g \bar \b} = - (J_t)^{\bar \g}_\a g_{\a  \bar \g}$ and $(J_t)^\g_{\bar \a} g_{\g \bar \b} = - (J_t)^{\g}_\a g_{\bar \a \g}$ and hence 
$$ g^{\a \bar \b} (J_t)^\g_{\bar \b}  = - (J_t)^{\a}_{\bar \b} g^{\bar \b   \g} \ ,\qquad  g^{\a \bar \b} (J_t)^{\bar \g}_{\bar \b} = -  (J_t)^{\a}_{\b} g^{\b \bar   \g}\ .$$
Using this property, we see that \eqref{tremototerra-1} can be also written  as 
\begin{multline*}
X^\a_t  
 = \frac{Y^\a_t}{\z} + i g^{\a \bar \b} e_{\bar \b} (\s_t) - \frac{i}{2}  (J_t - \Jst)^\a_{\bar \d} g^{ \bar  \d\g} \bigg(\frac{H_{t \g}}{\z} + \frac{\overline{H_{t \bar \g}}}{\bar \z} \bigg) -\\
  - \frac{i}{2}  (J_t - \Jst)^{\a}_{\d} g^{\d \bar  \g}\bigg(\frac{H_{t \bar \g}}{\z} + \frac{\overline{H_{t  \g}}}{\bar \z} \bigg) + \\
    + \frac{i}{2} (J_t - \Jst)^\a_{\bar \d} g^{\bar \d  \g} \left(
\overline{\left(H_{t \bar \g} - Y^0_t e_{\bar \b}(\log \r^2)\right)}\z  + \left(H_{t   \g} - Y^0_t e_{ \g}(\log \r^2)\right)\bar \z \right)  + \\
 +  \frac{i}{2}  (J_t - \Jst)^{\a}_{\d} g^{\d \bar  \g} \left(
\overline{\left(H_{t\g} - Y^0_t e_{\g}(\log \r^2)\right)}\z  + \left(H_{t\bar \g} - Y^0_t e_{\bar \g}(\log \r^2)\right)\bar \z \right)\ . 
 \end{multline*}
Since   $X_t$ must be a special vector field, determined by $v_t$ and some appropriate functions $\s_t$, $\wt Y^\a_t$, we see that the  $\wt Y^\a_t$ must be  equal to \eqref{crucial3}. Note also that such functions \eqref{crucial3}  automatically satisfy  the condition 
$\lim_{k \to \infty}(\wt Y^\a \z)|_{x_k}  = 0$ for each sequence   $x_k \to 0$. To see this it suffices to remember that  $\Jst$ is the standard complex structure of the indicatrix $\overline \cI_{x = 0}$ of $((\OB, J_t), \t_o)$ and   coincides with $J_t$ at $x = 0$. Hence  both the limits
$\lim_{k \to \infty}(J_t - \Jst)^{\a}_{\d}|_{x_k}$, $ \lim_{k \to \infty}(J_t - \Jst)^{\bar \a}_{\d}|_{x_k}$ are surely  $0$. 
\par
\smallskip
The last claim of the proposition is an immediate consequence  of  the fact that if $X_t$ is a special vector field satisfying \eqref{crucial3}  and $J_t$ satisfies the differential problem \eqref{eiaeia}, the first part of the proof implies that the $J_t$ are complex structures on $\OB$  induced by  L-complex structures on $\TOB$.
\end{pf}
 \subsection{Proof of  Theorem \ref{main}}\label{lastsect}
 The proof of Theorem \ref{main} is a consequence of the following   existence result on abstract fundamental pairs. 
 \begin{lem} \label{lemma52} Let $((\OB, J), \t_o)$ be a manifold of circular type in normal form of class $\cC^{k}$, $k\geq 3$, and   $v^o_s$, $s \in [0,1]$,    a fixed one-parameter family of vectors in $T_0\bB^n$ of class $\cC^{k}$  in $s$. 
 Then there is a maximal value $s_o \in (0, 1]$ such that for all $0 < s < s_o$ there exists  a  fundamental pair $(J_t, X_t)$ on $((\OB, J), \t_o)$ of class $\cC^{k-1,\a}$ for any $\a \in (0,1)$,  guided by the curve $v_t\= v^o_{s t}$, $t \in [0,1]$,  in which $X_t$
  is the special vector field $X_t  = X^{[\r_t](v_t, \s_t, \wt Y^\a_t)}$ with  $\s_t \equiv 0$, the $\wt Y^\a_t$ as  in   \eqref{crucial3} and $\r_t = \r$  for  each $t$.  The dependence of $X_t$ on $t$ is at least $\cC^{k-1, \a}$. \par
 In case the maximal value $s_o$ is strictly less than $1$,  there exists no fundamental pair $(J_t, X_t)$, $t \in [0,1]$, guided by the curve  $v_t\= v^o_{t s_o}$, $t \in [0,1]$. \par
\end{lem}
Indeed,  given a vector  $v \in \bC^n$ with $\|v\| = 1$,  one can  consider an associated curve of vectors $\mathbf v_s$, $s \in (0,1)$,  defined as  in \eqref{lallalla}. 
By  Theorem \ref{existenceanduniqueness}, Remark \ref{whatweneed}, Proposition \ref{lemma44} and Lemma \ref{lemma52},   for each $\l \in (0, \l_o\= s_o)$,  there is a one-parameter family of Monge-Amp\`ere exhaustions $\t^{(x_t)}: \OB \longrightarrow [0,1]$,  $t \in [0,1]$, centered at the points of the straight segment $x_t \= t (\l v)$, each  of class $\cC^{k-1,\a}$ for all $\a \in (0,1)$. The last claim of the lemma implies that,  if  $\l_o  < 1$, then  there is no one-parameter family of Monge-Amp\`ere exhaustions  centered at the points of the straight segment $x_t \= t (\l_o v)$ that  is well defined   also in  $x_{t = 1}$.\par
\smallskip
\begin{pf} Consider a set of polar coordinates $(\z, w^\a, \bar \z, \overline{w^\a})$,  an  adapted  polar frame field $(Z, e_\a = e^{(\r)}_\a)$ and the corresponding  dual coframe field $(Z^*, e^\a)$, both determined by the function $\r$ which gives the Minkowski function of the indicatrix $\overline{\cI}_{x = 0}$ of $(\OB, J)$.   All tensor fields $J_t$ of an abstract fundamental pair  have the form
\beq\label{tensoruccio} J_t = i Z \otimes Z^* - i \overline{Z \otimes Z^*} + J_{t \a}^{\ \b} e^\a \otimes e_\b + J_{t \bar \a}^{\ \b} \overline{e^\a} \otimes e_\b +  J_{t \a}^{\ \bar \b}
 e^\a \otimes \overline{e_\b} +  J_{t \bar \a}^{\ \bar \b} \overline{e^\a} \otimes \overline{e_\b}\ ,\eeq
 for  $\bC$-valued $\cC^{k-1,\a}$ functions $J^{\ A}_{t B}$  satisfying the reality conditions 
 $J_{t \bar \a}^{\ \bar \b} = \overline{J_{t \a}^{\b}}$,  $J_{t \a}^{\ \bar \b} = \overline{J_{t\bar \a}^{\ \b}}$. 
 Hence, the lemma corresponds to the existence of   complex functions $J_{t \a}^{\ \b}$, $J_{t \bar \a}^{\ \b}$, $t \in [0,1]$,  satisfying the required regularity conditions and  such that   tensors  \eqref{tensoruccio} verify: 
 \begin{itemize}
 \item[{\it a)}] $J_0 = J$; 
 \item[{\it b)}]  $J_t^2|_x = - Id_{T_x\TOB}$ for each $t \in [0,1]$ and $x\in \TOB$; 
 \item[{\it c)}]  the  Nijenhuis tensor  $N_{J_t}$ 
  is identically zero for each $t$; 
  \item[{\it d)}] for each $t \in [0,1]$ and for any vector field $E \in \cH$ 
 \begin{align}
 \label{5.13} & \frac{d J_t}{dt}(Z) +  (\cL_{X^{[\r] (v_t, 0, \wt Y_t^\a)}} J_t)(Z) = 0\ ,\\
\label{5.14}  &   \frac{d J_t}{dt}(E) +  (\cL_{X^{[\r] (v_t, 0, \wt Y_t^\a)}} J_t)(E) = 0\ .
\end{align}
Here, for  $Y^0_t$, $H_{t \bar \g}$,  $H_{t  \g}$  as in \eqref{C7bis}, \eqref{4.24}  for $v = v_t$,
\begin{multline} \label{crucial3bis} 
 \wt Y^\a_t  =   - \frac{i}{2}  (J_t - \Jst)^\a_{\bar \d} g^{ \bar  \d\g} \bigg(\frac{H_{t \g}}{\z} + \frac{\overline{H_{t \bar \g}}}{\bar \z} \bigg) -\frac{i}{2}  (J_t - \Jst)^{\a}_{\d} g^{\d \bar  \g}\bigg(\frac{H_{t \bar \g}}{\z} + \frac{\overline{H_{t  \g}}}{\bar \z} \bigg) + \\
    + \frac{i}{2} (J_t - \Jst)^\a_{\bar \d} g^{\bar \d  \g} \left(
\overline{\left(H_{t \bar \g} - Y^0_t e_{\bar \g}(\log \r^2)\right)}\z  + \left(H_{t   \g} - Y^0_t e_{ \g}(\log \r^2)\right)\bar \z \right)  + 
\end{multline}
$$
 +  \frac{i}{2}  (J_t - \Jst)^{\a}_{\d} g^{\d \bar  \g} \left(
\overline{\left(H_{t\g} - Y^0_t e_{\g}(\log \r^2)\right)}\z  + \left(H_{t\bar \g} - Y^0_t e_{\bar \g}(\log \r^2)\right)\bar \z \right)\ .$$
\end{itemize}
We now observe that   \eqref{5.13}  can be safely neglected, since 
$X^{[\r](v_t, 0, \wt Y^\a)}$ is a special vector field and hence, by the proof of Proposition \ref{lemma44},   such  equation is identically satisfied.  
It is also convenient  to decompose the vector field   $X_t \= X^{[\r](v_t, 0, \wt Y^\a)}$  into the  sum
\begin{align}
\nonumber & X_t  = \bX_t +  \bY_t\qquad \text{where}\ \ \bY_t  =   \wt Y^\a_t e_\a + \overline{\wt Y^\a_t} e_{\bar \a}\ \  \text{and}\\
\nonumber & \bX_t =  \left(\frac{Y^0_t}{\z} - \overline{Y^0_t} \z\right)  Z + \frac{Y^\a_t}{\z} e_\a + \left(\frac{\overline{Y^0_t}}{\bar \z}  - Y^0_t \bar \z\right) \overline Z + \frac{Y^{\bar \a}_t}{\bar \z}  e_{\bar \a}\ ,\\
\nonumber & \ \  \text{so that \eqref{5.14} 
can be written as}\\
\label{5.19} & \frac{d J_t}{dt}(E) +  \cL_{\bX_t + \bY_t} J_t(E) = 0\ ,
\end{align}
Note  that    $\bX_t$ is  independent of $J_t$,  while $\bY_t$   depends in a linear way on $J_t - \Jst$.   
\par
\medskip
Our  proof of  the existence of   $\cC^{k-1,\a}$  tensor fields  $J_t$ satisfying  {\it (a)} - {\it (d)} will  follow  from the  following three steps:\\[3pt]
{\it Step 1.} Translate the problem into an equivalent one on  the deformation tensors $\phi_t$ for the  L-complex structures $J_t$ (for the definition, see  \S \ref{review} below). \\[7pt]
{\it Step 2.}  Prove that all conditions  on the deformation tensors $\phi_t$  are  satisfied if just the   single  equation on $\phi_t$ corresponding to \eqref{5.19} is verified. \\[7pt]
{\it Step 3.}  Prove  the existence of  solutions to the equation corresponding to \eqref{5.19} for  $v_t = v^o_{\l t}$, $t \in [0,1]$, for all $\l$ in an appropriate open interval $(0, \l_o)$.\\[3pt]
According to this,  after  a short  introduction to the theory deformation tensors, the presentation  will be structured  in three parts, one per each of such steps.  After that,  we   discuss the  $\cC^\infty$ and $\cC^\o$ cases and make a short remark on examples.\par
\subsubsection{An introduction to  the deformation tensors of  L-complex structures}
\label{review}
Any tensor field  $J_t$ of the form \eqref{tensoruccio}  is  determined 
 by its restrictions to the distribution  $\cH$.  If we   also assume that it is an almost complex structure --  that is, condition {\it (b)} -- then it is actually    determined by   its associated  $J_t$-antiholomorphic distributions  $ \cH^{01}_t  \subset \cH^\bC$.  If the family  $J_t$ satisfies also   {\it (a)},  {\it (c)} and {\it (d)},  then, by Lemma \ref{lemma3.8},   it is a    one-parameter family of L-complex structure, according to Definition \ref{definition5.1}.  \par
  \smallskip
Consider now the   holomorphic and anti-holomorphic subdistributions $\cH^{10}$, $\cH^{01} = \overline{\cH^{10}}$ of $\cH^\bC \subset T\TOB$,  given by  the standard complex structure $\Jst$ of the indicatrix $\overline \cI_{x = 0}$ of $(\overline{\bB^n}, J_{t = 0} = J)$  and the analogous   distributions   $\cH^{10}_t$, $\cH^{01}_t = \overline{\cH^{10}_t} \subset\cH^\bC$,  given by the $\pm i$-eigenspaces in $\cH^\bC$ of the almost complex structures $J_t$, $t \in [0, 1]$.   By the results in  \cite{BD,PS},  each  distribution 
  $ \cH^{01}_t $  is  determined by   a unique     tensor field
 $\phi_t$ in  $\cH^{01*} \otimes \cH^{10} $,  which allows to express all subspace $\cH^{01}_{t}\big|_x \subset T_x\TOB$, $ x \in \TOB$,  in the form
 $$\cH^{01}_{t}\big|_x \= \{v \in \cH^\bC|_x\ :\ v = w + \phi_t(w)\ \text{for some}\ w^{01} \in \cH^{10}|_x\}\  .$$
 Such   $\phi_t$ is called {\it deformation tensor} of   $J_t$.  \par
 \smallskip
 Note that,  in each set  of adapted polar coordinates,  the components of  $J_t$ in the complex frames field $(Z, \wh e_\a \= e_\a + \overline{\phi_t(e_{\bar \a})}, \overline Z, \wh e_{\bar \a} = e_{\bar \a} + \phi_t(e_{\bar \a}))$
  are   constant and all equal  to $\pm i$ or $0$. This means that   the regularity of $J_t$ with respect to the  adapted polar coordinates  is always the same of   the regularity of    $\phi_t$ in those  coordinates. \par
 \medskip 
 The properties  that $J_t$ is an L-complex structure, i.e.   it is integrable or, equivalently, it  satisfies {\it (c)}, and that $\t_o$ is  Monge-Amp\`ere exhaustion for   $(\bB^n, J_t)$, corresponds to  the following   conditions on $\phi_t$, to  be   satisfied  for all  $0 \neq X, Y \in \cH^{01}$:
\begin{itemize}[itemsep=6pt plus 5pt minus 2pt, leftmargin=18pt]
\item[A)] for each $(t, y) \in [0,1] \times \TOB$,  one has that 
\beq \label{mannaggia} \ker \left(I_{T_x \TOB} - \overline \phi_t \circ \phi_t|_{y}\right)\qquad \text{is trivial}\ ;\eeq 
\item[B)]  
$ dd^c_{\Jst} \tau_o(\phi_t(X), \overline{\phi_t(X)}) {<} dd^c_{\Jst} \tau_o(\overline X, X)$ and $dd^c_{\Jst}\tau_o(\phi_t(X),$ $Y){+}dd^c_{\Jst} \tau_o(X, $ $\phi_t(Y))$ $= 0$;    
 \item[C)] 
 $ [\overline Z, X + \phi_t(X)] \in \cH^{01}_t$ or, equivalently,  
 \beq\label{int1-bis}  \cL_{\overline Z} \phi_t = 0\ ,\eeq
 where  $Z$ is the generator of  $\cZ^{10}$  defined in  \eqref{vectorfields};
 \item[D)]  $ [X + \phi_t(X), Y + \phi_t(Y)] \in \cH^{01}_t$. 
\end{itemize}
The geometrical interpretations of these conditions are the following. 
\begin{itemize}[itemsep=6pt plus 5pt minus 2pt, leftmargin=18pt]
\item[1)] Condition (A)  corresponds to the fact that the  tensors $\phi_t \in \cH^{01*} \otimes \cH^{10}$ determine a   direct sum decomposition  $\cH^\bC = \cH^{10}_t\oplus \cH^{01}_t $ with 
\beq \label{directsum}
\cH^{01}_t|_x \= \{v \in \cH^\bC|_x\, :\, v = w + \phi_t(w)\, ,w^{01} \in \cH^{10}|_x\}\ 
\text{and}\ \  \cH^{10}_t = \overline{\cH^{01}_t}
\eeq
 (see e.g. \cite{EO}).
If it holds, the one-parameter family of tensors $\phi_t$ defines an associated family of (possibly non-integrable) complex structures $J_t$ on the distribution $\cH$. 
\item[2)] Condition (B) expresses the condition that each level set of $\t_o$ is strongly pseudoconvex and hence that $\t_o$ is a  Monge-Amp\`ere exhaustion for $(\OB, J_t)$ if  $J_t$ is the complex structure determined  by the direct sum decomposition \eqref{directsum}. Note that (B)  implies (A), but  (A) does not imply (B). 
\item[3)] Conditions (C) and (D)  are equivalent  to the condition of integrability for the one-parameter family of complex structures  $J_t$ determined by the decomposition \eqref{directsum}.
\end{itemize}
  \par
\smallskip
The above  conditions are also sufficient in the following sense:  {\it if a  tensor field $\phi_t \in \cH^{01*} \otimes \cH^{01}$   satisfies {\rm (A) -- (D)} {\rm (actually, (B) -- (D) are enough, since  (B) implies (A))}, then the corresponding distribution $\cH^{01}_{t} \subset \cH^\bC$ uniquely determines an  L-complex structure $J_t$} (\cite{BD,PS}). 
 Finally we remark that conditions {\rm (A) -- (D)} are meaningful under very mild differentiability assumptions:  $\cC^{2}$-smoothness  on the data is sufficient.\par
\smallskip
We conclude this short review with  the  proof  that  the identity 
\eqref{obstruction} is  a  consequence of the integrability  of the complex structures $J_t$, as  claimed in the proof of Proposition \ref{lemma44}. Indeed,   using adapted polar frames, we may write  \eqref{obstruction} in the form
 \beq \label{obstruction-1} d d^c_{\Jst} \t_o\bigg( \frac{Y^\a_t}{\z} e_\a  +\frac{\overline{Y^\a_t}}{\bar \z} e_{\bar \a} + \bY_t , [\overline Z, e_{\bar \b}  + \phi_t(e_{\bar \b})]\bigg) = 0
\eeq
Since 
$\cL_{\bar Z}\phi_t = 0$ by    (C), we have     $[\bar Z, e_{\bar \a} + \phi_t(e_{\bar \a})] = 0$ and  \eqref{obstruction-1}  is  satisfied. \par
\subsubsection{Step 1 - Translation into a problem on deformations tensors}
\label{subsubsect1} 
By the previous  discussion, our  original problem 
translates into  a problem on  one-parameter families of $\cC^{k-1, \a}$ deformation tensors $\phi_t$ satisfying  (A) -- (D) and  the differential equations  \eqref{5.19} with  initial condition $\phi_0 = \phi_J$, where  $\phi_J$ is  the deformation tensor that allows to express    $J_{t = 0} = J$ as  deformation of  the complex structure $\Jst = \Jst{}^{(\rho)}$ of the indicatrix of $(\OB, J)$ at $x = 0$. Moreover, we claim that  condition (B)  is  satisfied whenever  all other conditions  hold and thus  it can be  neglected. This is because, if   the  $\phi_t$  satisfy (A),  (C),  (D) and  \eqref{5.19},   then the    family  of  integrable complex structures  $J_t$, determined  by means of the direct sum decomposition \eqref{directsum}, surely satisfies all hypotheses of   Lemma  \ref{lemma3.8}. This implies that   the exhaustion  $\t_o$ is  Monge-Amp\`ere  for  $(\overline{\bB^n}, J_t)$ for each $t$ and  that also  (B) holds.  \par
\medskip
In the remaining part of this section we derive a formulation  of \eqref{5.19}  as an explicit condition  on the $\phi_t$. For this,
recall that  for each   $E + \phi_t(E)  \in \cH^{01}_t$, we have that 
 $J_t(E + \phi_t(E)) = - i (E + \phi_t(E))$. Taking  Lie derivatives of  both sides with respect to the vector field $\frac{d}{dt} + X_t$ of $[0,1] \times \TOB$,   we get   
$$   \cL_{\frac{d}{dt} + X_t} (J_t) ( E + \phi_t(E))   = - i (I + i J_t) \cL_{\frac{d}{dt} + X_t}( E + \phi_t(E) ) \ .$$
This means that \eqref{5.19} holds if and only if  
$$ \cL_{\frac{d}{dt} + X_t}( E + \phi_t(E)|_x \in \ker (I + i J_t)|_x$$
 for   each  $(t,x)$,  i.e. if and only if 
\beq  \label{5.22} 
 \cL_{\frac{d}{dt} + X_t}  (E + \phi_t(E)) \in  \cH^{01}_t + \cZ^{01}\qquad
  \text{for all}\ E + \phi_t(E) \in \cH^{01}_t\ .
\eeq
On the other hand,  since  it  is the special vector field  described in Proposition  \ref{lemma44}, the vector field $X_t$ necessarily satisfies \eqref{6.15bis}. 
This implies that  $\cL_{ X_t}  (E + \phi_t(E))$ is always in  $ \cH^\bC$  and that  \eqref{5.22}
can be  replaced by the weaker condition
  \beq  \label{5.22bis} \cL_{\frac{d}{dt} + X_t}  (E + \phi_t(E)) \in \cH^{01}_t + \cZ^\bC, \qquad \text{for all}\ E \in \cH^{01}\ .\eeq
This condition can be stated as a system of p.d.e. as follows. Fix a (local) complex frame field $(e_{\bar \a})$ for $\cH^{01}$ and  for each $t$ consider the uniquely associated complex frame fields  $(\bE_{t|\bar \a} = e_{\bar \a} + \phi_t(e_{\bar \a}))$ and  $(\bE_{t \a} = \overline{\bE_{t|\bar \a}})$ for $\cH^{01}_t$ and  $\cH^{10}_t$, respectively.  Then, let $g^{(\phi_t)} \in \Hom( \cH \times \cH, \bR)$ be the only  $J_t$-invariant tensor field such that 
\beq \label{innerprod} g^{(\phi_t)}(\bE_{t|  \a}, \bE_{t|\bar \b}) = (I - \bar \phi_t \circ \phi_t)_\a^\g g_{\g \bar \b} \ ,\eeq
where $g_{\a \bar \b}$ are defined in \eqref{att} and where   $ (\cdot )_\a^\g, (\cdot )_{\bar \a}^{\bar \g}$  denote the components  of the considered $(1,1)$-tensors of $\cH^\bC$ in the frames $(e_\a, e_{\bar \a})$. Note that if $\phi_t$ satisfies (B) (hence, $\|\phi_t\|_{\cC^0} < 1$ for each $t \in [0, \ve]$),  then $g^{(\phi_t)}$ is  non-degenerate on  each space $\cH_x$. Observe also that, by construction and $J_t$-invariance, 
\beq \begin{split}
& g^{(\phi_t)}(E_{t|\bar \a}, \bE_{t|\bar \b}) = 0\ ,\\
& g^{(\phi_t)}(e_\a, \bE_{t|\bar \b})  = g^{(\phi_t)}(\bE_{t|\a}, \bE_{t|\bar \b}) - \phi_{t\a}^{\bar \g} g^{(\phi_t)}(e_{\bar \g}, \bE_{t|\bar \b}) = \\
& \hskip 0.5 cm = g^{(\phi_t)}(\bE_{t|\a}, \bE_{t|\bar \b}) - \bar \phi_{t\a}^{\bar \g} g^{(\phi_t)}(\bE_{t|\bar \g}, \bE_{t|\bar \b}) + (\bar \phi_t \circ \phi_t)_\a^\g g^{(\phi_t)}(e_\g, \bE_{t|\bar \b})  = \\ 
& \hskip 0.5 cm = g^{(\phi_t)}(\bE_{t|\a}, \bE_{t|\bar \b}) + (\bar \phi_t \circ \phi_t)_\a^\g g^{(\phi_t)}(e_\g, \bE_{t|\bar \b})\ . 
\end{split}
\eeq
In particular, from the last equality we get that 
\beq \label{convenient} g^{(\phi_t)}(e_\a, \bE_{t|\bar \b})  = ((I - \bar \phi_t \circ \phi_t)^{-1})_\a^\g g^{(\phi_t)}(\bE_{t|\g}, \bE_{t|\bar \b})  = g_{\a \bar \b}\ .\eeq
{\it This identity will appear to be  quite useful in the last part of the proof and it is the main motivation for considering the above definition for  $g^{(\phi_t)}$}. \par
\smallskip
Notice that using the tensor fields $g^{(\phi_t)}$,  condition  \eqref{5.22bis}  becomes equivalent to the system of p.d.e.'s  
  \beq  \label{5.22bis-1} g^{(\phi_t)}(\cL_{\frac{d}{dt} + X_t}  \bE_{t|\bar \a}, \bE_{t|\bar \b})  = 0\ .\eeq
\par
\medskip
We now recall that in \eqref{5.22bis}  the vector field $X_t = \bX_t + \bY_t$  {\it does depend} on the unknown $\phi_t$,  due to the fact that $\bY_t $  depends on $J_t - \Jst$.  We need to   make  such  dependence fully  explicit. 
For this, we first recall that, for each $t$, a vector field $Y$ in  $\cH^\bC$ uniquely  decomposes not only as a sum of holomorphic and anti-holomorphic  components with respect to  $\Jst$, but  also as a sum of holomorphic and  anti-holomorphic  components with respect to  $J_t$.   We  denote  such  two distinct   decompositions by 
 $ Y = Y^{10} + Y^{01}  = Y^{10}_t + Y^{01}_t$.
On the other hand, we know that  the components $Y^{10}_t$,  $Y^{01}_t$  have  the form 
$$ Y^{10}_t = \wh Y^{10(t)} + \overline{\phi_t}(\wh Y^{10(t)})\ ,\qquad Y^{01}_t = \wh Y^{01(t)} + \phi_t(\wh Y^{01(t)})$$
 for  appropriate  $\wh Y^{10(t)} \in \cH^{10}$ and $\wh Y^{01(t)} \in \cH^{01}$. A straightforward algebraic  computation shows that  such vectors are  expressed in terms of    $Y^{10}$,  $Y^{01}$ by 
 \beq \wh Y^{10(t)} \= (I -  \phi_t \circ \overline{\phi_t})^{-1}\left(Y^{10} - \phi_t(Y^{01})\right)\ ,$$
 $$  \wh Y^{01(t)}\= (I - \bar \phi_t \circ \phi_t)^{-1}\left(Y^{01} - \overline{\phi_t}(Y^{10})\right)\ , \eeq
 where each   $\phi_t|_x \in \Hom(\cH^{01}_x, \cH^{10}_x)$ is here considered  in  $\Hom(\cH^\bC_x, \cH^\bC_x)$,  acting  trivially on $\cH^{10}_x$. Note that, {\it due to condition (A), the linear operators $(I - \bar \phi_t \circ \phi_t)$ are invertible, so that the above expressions are  meaningful}.  It   follows that  the $J_t$-holomorphic and $J_t$-anti-holomorphic parts of   the  $e_\a \in \cH^{10}$ are 
 $$e_\a = \left\{(I -  \phi_t \circ \overline{\phi_t})^{-1}\left(e_\a\right) + \overline{\phi_t}((I -  \phi_t \circ \overline{\phi_t})^{-1}\left(e_\a\right))\right\} - $$
 $$ -  \left\{ (I - \bar \phi_t \circ \phi_t)^{-1}\left(\overline{\phi_t}(e_\a)\right) + \phi_t( (I - \bar \phi_t \circ \phi_t)^{-1}\left(  \overline{\phi_t}(e_\a)\right))\right\}\ . $$
Hence
 $$J_t(e_\a) =  i \left\{(I -  \phi_t \circ \overline{\phi_t})^{-1}\left(e_\a\right) + \overline{\phi_t}((I -  \phi_t \circ \overline{\phi_t})^{-1}\left(e_\a\right))\right\} + $$
 $$ +  i \left\{ (I - \bar \phi_t \circ \phi_t)^{-1}\left(\overline{\phi_t}(e_\a)\right) + \phi_t( (I - \bar \phi_t \circ \phi_t)^{-1}\left(  \overline{\phi_t}(e_\a)\right))\right\} $$
and  $(J_t - \Jst)(e_\a) = J_t(e_\a) - i e_{\a} $ is 
 $$(J_t - \Jst)(e_\a) =  2 i  \left\{ (I - \bar \phi_t \circ \phi_t)^{-1}\left(\overline{\phi_t}(e_\a)\right) + \phi_t( (I - \bar \phi_t \circ \phi_t)^{-1}\left(  \overline{\phi_t}(e_\a)\right))\right\} \ . $$
From the identity $\phi_t \circ (I - \bar \phi_t \circ \phi_t)^{-1} =  (I - \phi_t \circ \bar \phi_t)^{-1}  \circ\phi_t$ (it can be  checked  using the expansion $(I - \bar \phi_t \circ \phi_t)^{-1} = \sum_{r = 0}^\infty (\bar \phi_t \circ \phi_t)^r$) we may also say that 
 $$(J_t - \Jst)(e_\a) {=}  2 i  \big\{ \left((I - \bar \phi_t \circ \phi_t)^{-1}\hskip -0.1cm {\circ}\overline{\phi_t}\right)(e_\a) {+} \left((I -  \phi_t \circ \bar \phi_t)^{-1}\hskip -0.1cm{\circ} \left( \phi_t  \circ  \overline{\phi_t}\right) \right)(e_\a)\big\},  $$
 from which we get the components of $J-\Jst$  w.r.t. $(e_\a, e_{\bar \b})$. They form the matrix  
$$\left(
\smallmatrix 
 (J_t - \Jst)_\a^\b = 2 i ((I -  \phi_t \circ \bar \phi_t)^{-1})^{\b}_{\g} ( \phi_t \circ \bar  \phi_t)^{\g}_{\a} & 
  (J_t - \Jst)_\a^{\bar \b} = 2 i ((I - \bar \phi_t \circ \phi_t)^{-1})^{\bar \b}_{\bar \g} \bar \phi_t{}^{\bar \g}_\a  \\
 (J_t - \Jst)_{\bar \a}^{\bar \b} =  - 2 i ((I - \bar  \phi_t \circ  \phi_t)^{-1})^{\bar \b}_{\bar \g} ( \bar \phi_t \circ    \phi_t)^{\bar \g}_{\bar \a} &\ \ 
  (J_t - \Jst)_{\bar \a}^{\b} = - 2 i ((I -  \phi_t \circ \bar \phi_t)^{-1})^{\b}_{\g}  \phi_t{}^{\g}_{\bar \a} 
  \endsmallmatrix\right)\ .
  $$
Thus,  the explicit dependence of the $\wt Y^\a_t$ on the $\phi_t$  we are looking for  is 
\begin{multline} \label{5.27}
  \wt Y^\a_t  =  ((I -  \phi_t \circ \bar \phi_t)^{-1})^{\a}_{\g} ( \phi_t \circ \bar  \phi_t)^{\g}_{\d} g^{\d \bar  \h}\bigg(\frac{H_{t \bar \h}}{\z} + \frac{\overline{H_{t  \h}}}{\bar \z} + \\
 - \overline{\left(H_{t\h} - Y^0_t e_{\h}(\log \r^2)\right)}\z  - \left(H_{t\bar \h} - Y^0_t e_{\bar \h}(\log \r^2)\right)\bar \z\bigg)  +\\
  +((I -  \phi_t \circ \bar \phi_t)^{-1})^{\a}_{\g}  \phi_t{}^{\g}_{\bar \d} g^{ \bar  \d\h} \bigg(\frac{H_{t \h}}{\z} + 
  \frac{\overline{H_{t \bar \h}}}{\bar \z} - \\
 - \overline{\left(H_{t \bar \h} - Y^0_t e_{\bar \h}(\log \r^2)\right)}\z  - \left(H_{t   \h} - Y^0_t e_{ \h}(\log \r^2)\right)\bar \z  \bigg)  \ .
  \end{multline}
Recalling that 
$(I -  \phi_t \circ \bar \phi_t)^{-1} \circ  \phi_t  =
\phi_t \circ (I -  \bar \phi_t \circ  \phi_t)^{-1}$ and that 
$$ (I -  \phi_t \circ \bar \phi_t)^{-1} \circ  (\phi_t \circ \bar \phi_t) =  (I -  \phi_t \circ \bar \phi_t)^{-1}  - I\ ,$$ 
from   \eqref{5.27} we see   that  $\bY_t = \wt Y^\a_t e_\a + \overline{Y^\a_t} e_{\bar \a}$ has the form
$$ \bY_t = (\bY_{t(1)} + \phi(\bY_{t(1)})) + (\bY_{t(2)} + \phi(\bY_{t(2)})) + \bY'_{t}$$
where
\beq
\begin{split} &\bY_{t(1)} \= ((I -  \bar \phi_t \circ  \phi_t)^{-1} \bar \phi_t)^{\bar \b}_{\d} g^{\d\bar \h} \bigg(\frac{\overline{H_{t \h}}}{\bar \z} + 
  \frac{H_{t \bar \h}}{\z} - \hskip 5cm \\
&\hskip 2 cm - \left(H_{t \bar \h} - Y^0_t e_{\bar \h}(\log \r^2)\right) \bar \z  - \overline{\left(H_{t   \h} - Y^0_t e_{ \h}(\log \r^2)\right)} \z  \bigg) e_{\bar \b}\ ,
\\
&\bY_{t(2)} \=
 ((I -  \bar \phi_t \circ  \phi_t)^{-1})^{\bar \b}_{\bar \d} g^{\bar \d   \h}\bigg(\frac{H_{t  \h}}{\z} + \frac{\overline{H_{t \bar  \h}}}{\bar \z} + \hskip 5 cm\\
&\hskip 2 cm - \overline{\left(H_{t\bar \h} - Y^0_t e_{\bar \h}(\log \r^2)\right)}\z  - \left(H_{t \h} - Y^0_t e_{ \h}(\log \r^2)\right)\bar \z\bigg) e_{\bar \b}\ ,\\
&\bY'_{t} {\=}  
 - g^{\bar \b   \h}\bigg(\frac{H_{t  \h}}{\z}  + \frac{\overline{H_{t  \bar \h}}}{\bar \z}
- \overline{\left(H_{t\bar \h} - Y^0_t e_{\bar \h}(\log \r^2)\right)}\z -  \left(H_{t \h} - Y^0_t e_{\h}(\log \r^2)\right)\bar \z\bigg) e_{\bar \b}\ .
\end{split}
\eeq
Since $ \sum_{i = 1}^2 (\bY_{t(i)} + \phi(\bY_{t(i)}))$ is  in $\cH^{01}_t$ for each $t$,  if  (D) holds,  then  
also $\cL_{\sum_{i = 1}^2 (\bY_{t(i)} + \phi(\bY_{t(i)}))} (E + \phi_t(E))$ is in $\cH^{01}_t$.
This means that, due to  the integrability condition (D),   in the set of the three conditions,  the equation \eqref{5.22bis} involving  the (real) vector field $X_t = \bX_t + \bY_t$ can be  replaced by the 
same equation with the (complex) vector field $X'_t = \bX_t + \bY'_{t}$ in place of $X_t$.  
 The crucial advantage of such replacement comes from  the fact that   {\it the new vector field $X'_t$ is totally  independent of $\phi_t$}.\par
\medskip
Summing up, our problem is now  reduced to  proving {\it the existence of a one parameter family of tensor fields $\phi_t \in \cH^{01*} \otimes \cH^{10}$ of class $\cC^{k-1, \a}$ on $[0,1] \times \TOB$ satisfying the initial condition $\phi_{t = 0} = \phi_J$, the nondegeneracy condition \eqref{mannaggia} and the system of equations 
\begin{align}
\label{5.30-1} &[\overline Z, E + \phi_t(E)]  \in \cH^{01}_t\ ,\\[6pt] 
\label{5.30-2} &[E + \phi_t(E), F + \phi_t(F)] \in \cH^{01}_t\ ,\\[6pt]
\label{5.30}  & \cL_{\frac{d}{dt} + X'_t}  (E + \phi_t(E)) \in \cH^{01}_t + \cZ^\bC
\end{align}
for all  $E + \phi_t(E)  \in \cH^{01}_t$.} Note that \eqref{5.30}  can be equivalently stated as  \eqref{5.22bis-1} putting $X'_t$ in  the place of $X_t$.
\par
\subsubsection{Step 2 - Reduction to a single equation}
In this subsection we prove that if $\phi_t$ is a $\cC^{k-1, \a}$ solution of  \eqref{5.30} with initial condition $\phi_{t = 0} = \phi_J$, then \eqref{5.30-1} and \eqref{5.30-2} are automatically satisfied. 
In fact, we prove this claim  only in the case of   the  $\cC^{k-1, \a}$ solutions $\phi_t$, which can be obtained as  limits of sequences of real analytic solutions $\phi_t^{(n)}$. This weaker result  
is enough for our purposes, because the solutions, of which  we 
prove the existence in the final step,  are precisely  of such a  kind.\par
\smallskip
As usual, given a solution $\phi_t$ to \eqref{5.30}, let  $J_t$ be the corresponding one parameter family of (possibly non-integrable) complex structures  of the form \eqref{tensoruccio}, for which $\phi_t$ is the deformation tensor relatively to $\Jst$. Moreover, for  any $x \in  \TOB$,  consider a (locally defined) complex frame field $(Z, E_\a, \bar Z, E_{\bar \a} \= \overline{E_\a})$ on a neighborhood of $x$, in which the $E_{\bar \a}$ are generators  for the distribution $\cH^{01}$ (for instance, we may consider an adapted polar frame field). For each $t$,  we denote by $(Z, \bE_{t|\a}, \bar Z, \bE_{t|\bar \a} = \overline{\bE_{t|\a}})$ the associated complex  frame  field, where the  vector fields $\bE_{t|\bar \a}= E_{\bar \a} + \f_t(E_{\bar \a})$   generate    the  $J_t$-antiholomorphic distribution $\cH^{01}_t$ and  the $\bE_{t|\a}\= \overline{\bE_{t|\bar \a}}$  the distribution  $\cH^{10}_t$.    We also  denote by  $(Z^*, \bE_{t|}^\a, \bar Z^*, \bE_{t|}^{\bar \a})$ the dual coframes field of $(Z, \bE_{t|\a}, \bar Z, \bE_{t|\bar \a} = \overline{\bE_{t|\a}})$ and we set
$$ \bF_{t|\bar \a \bar \b} \=  [E_{\bar \a} + \phi_t(E_{\bar a}), E_{\bar \b} + \phi_t(E_{\bar \b})] = [\bE_{t|\bar\a}, \bE_{t|\bar \b}]\ . $$
Observe that, since they are tangent to the level sets of $\t_o$, they have  the form 
$$\bF_{t|\bar \a\bar \b} = f_{t|\bar \a \bar \b}^\g \bE_{t|\g} +  f_{t|\bar \a \bar \b}^{\bar \g} \bE_{t|\bar \g} +  f_{t|\bar \a \bar \b}^0 (Z - \bar Z) $$
with  $f_{t|\bar \a \bar \b}^\g = \bE_{t|}^\g(\bF_{t|\bar \a\bar \b} )$, $f_{t|\bar \a \bar \b}^{\bar \g} =  \bE_{t|}^{\bar \g}(\bF_{t|\bar \a\bar \b} ) $, $ f_{t|\bar \a \bar \b}^0 = Z^*(\bF_{t|\bar \a\bar \b})$. Note also  that, being $J$  integrable,  the  functions   $f_{t|\bar \a\bar\b}^\g\big|_{t = 0}$, $f_{t|\bar \a\bar\b}^0\big|_{t = 0}$ are zero, while \eqref{5.30-2} holds  
if and only if the functions $f_{t|\bar \a \bar \b}^\g$, $f_{t|\bar \a\bar\b}^0$  vanish identically  
for each $t \in [0,1]$. \par
\medskip 
Assume for the moment that the integrable complex structure $J_{t = 0} = J$,  the solution $\phi_t$ to \eqref{5.30} and the family $v_t$ we  consider,  are all real analytic, so that 
also the components of the complex vector field $X'_t$ are real analytic. We want   to show that,  under these assumptions,  the    $f_{t|\bar \a\bar\b}^\g$, $f_{t|\bar \a\bar\b}^0$ are identically   $0$ and that \eqref{5.30-2} holds. For this we need the following technical result. \par
\smallskip
\begin{sublem}\label{sublemma} Assume that  all components of the complex vector field $X'_t$ on $[0,1] \times \TOB$ are real analytic and consider  a system of partial differential equations  on $[0,1] \times \TOB$ for an  unknown  $U: [0,1] \times \TOB \to \bC^N$ having the form in   sets of polar coordinates
\beq \label{realanalytic} (\frac{d}{dt} + X'_t) (U^a) = \cF^a(t, \z, \bar \z, w^i, \overline{w^i}, U^b)\ ,\qquad 1 \leq a \leq N\ ,\eeq
 with $\cF^a$ real analytic in    $(t, \z, \bar \z, w^i, \overline{w^i})$ and polynomial in $U^b$. 
Then, for each  $t \in [0,1]$ and each  real analytic     $V: \TOB \to \bC^N$,  there is a unique local real analytic solution $U$ on $[t_o - \ve, t_o+ \ve] \times \TOB$  that satisfies \eqref{realanalytic} with  initial condition $U|_{t = t_o}  = V$. \par
Furthermore,  if $\cS$ is   a subset of  $\cC^\o(\TOB, \bC^N)$,   for which there exists an  {\rm a priori} upper bound for the $\cC^0$-norms of all  solutions to \eqref{realanalytic} with  initial conditions  in $\cS$, then for each $V \in \cS$ there exists a unique real analytic solution $U:[0,1] \times \TOB \to \bC^N$ to the system \eqref{realanalytic}  with $U|_{t = 0}  = V$.
\end{sublem}
\begin{pf} For each point $(t_o, x_o) \in [0,1] \times \TOB$, consider a neighborhood   $\cU \subset [0,1] \times \TOB$, on  which we may consider a system of real  coordinates  $(x^i)$,  which allow to identify $\cU$ with an open subset  of   $\bC^{2n+1 } \cap \{\Im(z^i) = 0\}$. Assume also that $\cU$ is sufficiently small so that all restrictions to $\cU$ of  the components of $X'_t$ and of the coefficients of the polynomials  $\cF^a(U, \overline U)$ extend as holomorphic functions of some open neighborhood   $\cW \subset \bC^{2n+1}$  of $\cU$. In this way,  the {\it complex}  components $X^i$ of the complex vector field $\frac{d}{dt} + X'_t$ on $[0,1] \times \TOB$ 
  can be taken as the restrictions to $\cU$ of  some  holomorphic functions of the form $A^i + i B^j$ on $\cW$. These complex functions are such that the $A^i$ and $B^j$ take  only real values  at the points of $\cU \subset \cW  \cap \{\Im(z^i) = 0\}$. In other words, 
$$\frac{d}{dt} + X'_t {=} \left(A^i \frac{\p}{\p x^i} + i B^j\frac{\p}{\p x^j}\right)\bigg|_{\cU} {=} \left(A^i \frac{\p}{\p z^i} + i B^j\frac{\p}{\p z^j} + A^i \frac{\p}{\p \bar z^i} + i B^j\frac{\p}{\p \bar z^j} \right)\bigg|_{\cU}\!\!\!.$$ 
Any local real analytic solution $ U = (U^a)$ of \eqref{realanalytic}   admits  a holomorphic extension  on  an open  neighborhood of $\cU \subset \bR^{2n+1}$ in $\bC^{2n+1}$.  Being holomorphic, such   extension of $U$  is solution to the system of differential equations 
 \beq \label{realanalytic-1}  \wh X(U^a) - F^\a_{\bar \b}( z^i,  U^b) = 0\qquad \text{with}\qquad \ \wh X \= A^i \frac{\p}{\p z^i} + i B^j\frac{\p}{\p z^j}\ .\eeq
The graphs in $\bC^{2n+1} \times \bC^{N}$  of the holomorphic solutions $U: 
\cW \to \bC^N$  to \eqref{realanalytic} coincide 
with the complex submanifolds of $\bC^{2n+1} \times \bC^N$ that are tangent to the holomorphic vector field    
\beq \label{vector} \bW \= \wh X|_z + \cF^a(z,  U) \frac{\p}{\p U^a}\ .\eeq
These graphs  can be constructed  as in the classical method of real characteristics, namely by taking  unions of complex holomorphic curves tangent to the holomorphic vector field $\bW$ (they do exist by  complex Frobenius Theorem - we call them {\it complex characteristics}),  passing through the points of the graph of an initial condition $V: \{t_o\} \times 
\wt \cW \subset \cW\to \bC^N$. This  method of construction gives local solutions of \eqref{realanalytic-1}. Restricting them  to $\cU = \cW \cap \{\Im(z^i) = 0\}$, we get the desired local real analytic solutions of \eqref{realanalytic}.  \par
By previous observations, any local real analytic solution of \eqref{realanalytic} can be obtained in this way and, by  construction, any such solution is uniquely 
determined by its values at some level set at $t = t_o$ for some fixed $t_o \in [0,1]$. This  implies the first claim of the sublemma.\par
\medskip
For the second claim, observe that by  compactness of $\TOB$, the above construction  allows to determine a unique real analytic solution $U$ to \eqref{realanalytic} for each initial value $U_{t = 0} = V \in \cS$  
and we may assume that $U$ is defined  on a set of the form $[0, t_o] \times \TOB$ for some $t_o \in (0,T]$. If there exists an {\it a priori} bound for the $\cC^0$ norm of such solution, the vector field \eqref{vector} has   components that are  {\it a priori} bounded  at the boundary  points of the graph of the solution $U: [0, t_o] \times \TOB \to \bC^N$. This implies that such solution can be  extended to a solution
 $U: [0, t_o{+}\ve] \times \TOB \to \bC^N$ for some $\ve > 0$. A standard 
 open and closed argument yields the existence of a solution to \eqref{realanalytic} on  
$[0,1] \times \TOB$.
\end{pf}
Let us now go back to our discussion on the functions $f_{t|\bar \a\bar \b}^\g$,  $f_{t|\bar \a\bar \b}^0$, under the assumption of real analyticity of the data. 
We recall that  we are assuming that $\phi_t$ is solution to \eqref{5.30}, i.e. that  
\beq \label{8marzo} \cL_{\frac{d}{dt} + X'_t} \bE_{t|\bar \h} = \s_{t|\bar \h}^{\bar \mu} \bE_{t|\bar \mu} + \s_{t|\bar \h}^0 Z +  \s_{t|\bar \h}^{\bar 0}\overline Z\eeq
for some complex functions  $ \s_{t|\bar \h}^{\bar \mu}$, $\s_{t|\bar \h}^0$, $\s_{t|\bar \h}^{\bar 0}$. Note also that, by \eqref{6.15}, 
 \beq \label{festa} \cL_{\frac{d}{dt} + X'_t} Z = \cL_{\frac{d}{dt} + X'_t} \overline Z = 0\ .\eeq
 Hence,  by duality,   
\beq \label{8marzo-bis} \cL_{\frac{d}{dt} + X'_t} \bE_{t}^{\h} =  \r_{t|\mu}^{\h} \bE_{t}^{\mu} +  \r_{t|0}^{\h} Z^* +  \r_{t|\bar 0}^{\h} \overline Z^*\ ,\qquad \cL_{\frac{d}{dt} + X'_t} Z^* =  \r_{t|\mu}^{0} \bE_{t}^{\mu} +  \r_{t|\bar \mu}^0\bE_{t}^{\bar \mu}\eeq
for appropriate complex functions $ \r_{t|\mu}^{\h} $, $ \r_{t|0}^{\h}$, $\r_{t|\bar 0}^{\h}$, $ \r_{t|\mu}^{0}$, $\r_{t|\bar \mu}^0$. \par
\smallskip
 We now claim that the functions $\s_{t|\bar \h}^0$ and $\s_{t|\bar \h}^{\bar 0}$ are actually linear combinations of the functions $f^{0}_{t|\bar \a \bar \b}$ defined above. Indeed, we recall that the (complex) vector field $X'_t$ differs from the (real) vector field $X_t$ by the vector field  $ \sum_{i = 1}^2 (\bY_{t(i)} + \phi(\bY_{t(i)}) \in \cH^{01}_t$. We may therefore write 
that $X'_t = X_t +  \bY^{\bar \a}_{t} \bE_{t|\bar \a}$ for some functions $ \bY^{\bar \a}_{t} $ and 
\beq \label{donna} \cL_{\frac{d}{dt} + X'_t} \bE_{t|\bar \h} = \cL_{\frac{d}{dt} + X_t} \bE_{t|\bar \h} + \bY^{\bar \mu}_{t}  \bF_{t|\bar \nu \bar \h} -  \bE_{t|\bar \h} (\bY^{\bar \mu}_{t}) \bE_{t|\bar \mu}\ .\eeq
Since $X_t$ is a special vector field as described in Proposition \ref{lemma44}, by \eqref{6.15bis}, we have that  $\cL_{\frac{d}{dt} + X_t} \bE_{t|\bar \h} \in \cH^\bC$.
Comparing \eqref{donna} with \eqref{8marzo} and using the general structure of the vectors $\bF_{t|\bar\a\bar \b}$, we get that 
\beq  \label{donna-bis} \s_{t|\bar \h}^0 =  \bY^{\bar \mu}_{t}  f^0_{t|\bar \mu \bar \h}\ ,\qquad  \s_{t|\bar \h}^{\bar 0}  = - \bY^{\bar \mu}_{t}  f^0_{t|\bar \mu \bar \h}\ ,\eeq
i.e.,   $\s_{t|\bar \h}^0$,  $\s_{t|\bar \h}^{\bar 0}$ are  linear  combinations of   $f^{0}_{t|\bar \a \bar \b}$,  as claimed.
\par
\smallskip
Now, combining \eqref{8marzo}, \eqref{festa}, \eqref{8marzo-bis} and \eqref{donna-bis}, we get  
\beq \label{intr}
\begin{split} 
\bigg(\frac{d}{dt} + &X'_t\bigg)(f_{t|\bar \a\bar \b}^\g) =     \r_{t|\mu}^{\g} f_{t|\bar \a\bar \b}^\mu +  (\r_{t|0}^{\g} -  \r_{t|\bar 0}^{\g}) f_{t|\bar \a\bar \b}^0 +
\bE_{t|}^\g([\cL_{\frac{d}{dt} +  X'_t}\bE_{t|\bar \a}, \bE_{t|\bar \b}]) +\\
 &+  \bE_{t|}^\g([\bE_{t|\bar \a}, \cL_{\frac{d}{dt} + X'_t}\bE_{t|\bar \b}]) = \\
 & =    \r_{t|\mu}^{\g} f_{t|\bar \a\bar \b}^\mu +  (\r_{t|0}^{\g} -  \r_{t|\bar 0}^{\g}) f_{t|\bar \a\bar \b}^0 +\\
& +  \s_{t|\bar \a}^{\bar \mu}f_{t|\bar \mu\bar \b}^\g + \s_{t|\bar \b}^{\bar \mu} f_{t|\bar \a \bar \mu}^\g +   \bY^{\bar \mu}_{t}  f^0_{t|\bar \mu \bar \a} (\bE_{t|}^\g([Z, \bE_{t|\bar \b}])-  \bE_{t|}^\g([\overline{Z}, \bE_{t|\bar \b}]))+ \\
& + \bY^{\bar \mu}_{t}  f^0_{t|\bar \mu \bar \b} (\bE_{t|}^\g( [\bE_{t|\bar \a}, Z])-  \bE_{t|}^\g([\bE_{t|\bar \a}, \overline Z])) \ ,
 \end{split}
 \eeq
 \beq \label{intr-1}
\begin{split} 
\bigg(\frac{d}{dt} + &X'_t\bigg)(f_{t|\bar \a\bar \b}^0) =     \r_{t|\bar \mu}^{0} f_{t|\bar \a\bar \b}^\mu + \r_{t| \bar \mu}^{0} f_{t|\bar \a\bar \b}^{\bar \mu} + 
Z^*([\cL_{\frac{d}{dt} +  X'_t}\bE_{t|\bar \a}, \bE_{t|\bar \b}]) +\hskip 1 cm \\
 &+  Z^*([\bE_{t|\bar \a}, \cL_{\frac{d}{dt} + X'_t}\bE_{t|\bar \b}]) = \\
 & =    \r_{t|\mu}^{0} f_{t|\bar \a\bar \b}^\mu  -   f^0_{t|\bar \nu \bar \mu} \bY^{\bar \nu}_{t} f_{t|\bar \a\bar \b}^{\bar \mu} +\\
&  + \s_{t|\bar \a}^{\bar \mu} f^0_{t|\bar \mu \bar \b} +  \bY^{\bar \mu}_{t}  f^0_{t|\bar \mu \bar \a} (Z^*([Z, \bE_{t|\bar \b}])-  Z^*([\overline{Z}, \bE_{t|\bar \b}]))
 \\
 & + \s_{t|\bar \b}^{\bar \mu} f^0_{t|\bar \a \bar \mu} +  \bY^{\bar \mu}_{t}  f^0_{t|\bar \mu \bar \b} (Z^*([\bE_{t|\bar \a}, Z])-  Z^*([\bE_{t|\bar \a}, \overline{Z}]))\ .
 \end{split}
 \eeq
 This shows that the functions $f^0_{t|\bar \a \bar \b}$ and $f^\g_{t|\bar \a \bar \b}$ satisfy a system of the form \eqref{realanalytic}.
 Note also that   \eqref{intr} and \eqref{intr-1} admit  the trivial  solution for the initial value problem   $ f_{t|\bar \a\bar\b}^\g\big|_{t = 0} =  f_{t|\bar \a\bar\b}^0\big|_{t = 0} = 0$. 
 By the uniqueness  of the local solutions, proved  in  Sublemma \ref{sublemma},  this  implies that,  under the above real analyticity  assumptions, the condition  \eqref{5.30-2} is   identically satisfied. 
 \par
\smallskip
Assume now that $J_{t = 0} = J$ and that   $\phi_t$ and $v_t$ are not real analytic, but  that nonetheless there is a sequence of real analytic curves $v^{(n)}_t$ converging in $\cC^{k}$-norm to $v_t$, a sequence of  real analytic  complex structures $J^{(n)}$, converging in $\cC^{k}$-norm to $J$, and  a corresponding  sequence of  real analytic solutions $\phi^{(n)}$ to  \eqref{5.30} with initial data $\phi^{(n)}_{t = 0} = \phi_{J_n}$,   converging in $\cC^{k}$-norm to the solution $\phi_t$.  Note  that, for each of  these sequences, the convergence   is also in  $\cC^{k-1,\a}$ for any $\a \in (0,1)$ and that the $\cC^{k-1,\a}$-norms of the elements of  the sequence converge to the finite value of the $\cC^{k-1,\a}$-norm  of the limit.  In this case, the  associated sequences  $f_{t|\bar \a\bar \b}^{(n) \g}$, $f_{t|\bar \a\bar \b}^{(n) 0}$ converge in $\cC^{k-2}$-norm to the functions  $f_{t|\bar \a\bar \b}^{\g}$, $f_{t|\bar \a\bar \b}^{0}$  determined by $\phi_t$, and their $\cC^{k-3,\a}$-norms converge to finite values. Since: \par
a) the equations satisfied by the $f_{t|\bar \a\bar \b}^{(n) \g}$, $f_{t|\bar \a\bar \b}^{(n) 0}$ are linear, \par 
b)  the coefficients of the systems  have $\cC^{k-3, \a}$-norms which tend  to the $\cC^{k-3, \a}$-norms of 
the coefficients of the equations determined by  the $\phi_t$, \par
c)  the initial data tend to the zero functions in $\cC^{k-3,\a}$-norm, \\
we conclude that also 
the limit functions $f_{t|\bar \a\bar \b}^{\g}$, $f_{t|\bar \a\bar \b}^{0}$ are identically vanishing and  that  \eqref{5.30-2} holds also in this case. \par
\medskip
It remains to prove that any   solution $\phi_t$ to \eqref{5.30},  with $\phi_{t = 0} = \phi_J$,  necessarily satisfies also the  integrability condition \eqref{5.30-1}. For doing this, we first observe that, since we just proved that  \eqref{5.30-2} is surely satisfied, by the discussion of \S \ref{subsubsect1},  we also have that the (real) special vector field $X_t \= X'_t + \bY_{t(1)} + \phi\left(\bY_{t(1)}\right) + \bY_{t(2)} + \phi\left(\bY_{t(2)}\right)$ is so that 
$\cL_{\frac{d}{dt} + X_t} J_t = 0$. 
Now,  let  $g_t$  be a one-parameter family of $J_t$-invariant  Riemannian metrics on $\TOB$ such that 
$\cL_{\frac{d}{dt} + X_t} g_t = 0$. Since  we know that $\cL_{\frac{d}{dt} + X_t} J_t = 0$, one can construct  such a one-parameter family $g_t$ by considering the flow $\Phi_t$ on $\bR \times \TOB$, determined by the vector field $\frac{d}{dt} + X_t$, and impose that,  for each pair  vector fields $W, W' \in T(\{t\} \times \TOB)$, 
$$g_t(W, W') \= - \frac{i}{2\t} d d^c \t_J(\Phi_{-t*}(W), J \Phi_{-t*}(W'))\ .$$
The family of these Riemannian metrics  is clearly invariant under the flow $\Phi_t$, but  it is also $J_t|_x$-invariant at each $(t, x)$,  since also the family $J_t$ is invariant under that   flow.  Finally, we  set
$$ \bF_{t|\bar 0 \bar \b} \=  [\bar Z, E_{\bar \b} + \phi_t(E_{\bar \b})] = [\bar Z, \bE_{t|\bar \b}]\ , \qquad f_{t|\bar 0\bar \b\bar \g} \= g_t(\bF_{t|\bar 0\bar \b}, \bE_{t|\bar \g})\ .$$
Note that all the vector fields $\bF_{t|\bar 0\bar \b}$ are in $\cH^\bC$ and that \eqref{5.30-1} holds if and only if   $f_{t|\bar 0\bar \b\bar \g} =   0$ for each $t$.
Notice  also  that    $f_{t|\bar 0\bar \b\bar \g}|_{t = 0} = 0$, being  $J_{t = 0} = J$  integrable. \par
\smallskip
Imitating the above arguments, our next  goal is  to prove that  the  functions  $f_{t|0\b\g}$ are constrained by  a differential problem, which has the trivial map  as    unique  solution. This would  imply our desired result. 
For proving this, we recall   that being  $\phi_t$ solution to \eqref{5.30},   it also satisfies \eqref{5.22bis-1} and, due to \eqref{5.4}, 
$$\cL_{X_t} \overline Z \in \cH^{01}_t\qquad \text{for each} \ t \in [0,1]\ .$$ 
Hence, from  $\cL_{\frac{d}{dt}} Z = 0$,  \eqref{5.30} and \eqref{5.30-2}, we obtain
\beq\label{5.31}  
\cL_{\frac{d}{dt} + X_t}(\bF_{t|\bar 0\bar \b}) =    [\cL_{X_t} \overline Z, \bE_{t|\bar \b}] + [\overline Z, \cL_{\frac{d}{dt} +X_t}(\bE_{t|\bar \b}] 
 \in \cH^{01} + [\cZ^{01}, 
\cH^{01}_t]\ .
\eeq
 This means that  these Lie derivatives can be written as linear combinations of the $\bE_{t|\bar \a}$ and of the  $\bF_{t|0\b}$.  
 Since  $\cL_{\frac{d}{dt} + X_t} g_t = 0$,     $\cL_{\frac{d}{dt} + X_t}(\bE_{t|\bar \g}) \subset \cH^{01}_t $ and  $g_t(\cH^{01}_t, \cH^{01}_t) = 0$,  we get 
\beq\label{5.33}
\begin{split}
\left(\frac{d}{dt} + X_t\right)(f_{t|\bar 0\bar \b\bar \g} ) &= g_t(\cL_{\frac{d}{dt} + X_t} \bF_{t|\bar 0\bar \b} , \bE_{t|\bar \g}) +  g_t( \bF_{t|\bar 0\bar \b} ,\cL_{\frac{d}{dt} + X_t} \bE_{t|\bar \g})= \\
& = C^{\bar \ve}_{t|\bar 0\bar \b\bar \g} g_t(\bF_{t|\bar 0\bar \ve} , \bE_{t|\bar \g}) +   D_{t|\bar \g}^{\bar \d} g_t(\bF_{t|\bar 0\bar \b}, \bE_{t|\bar \d})  = \\
& = C^{\bar \ve}_{t|\bar 0\bar \b\bar \g} f_{t|\bar 0\bar \ve \bar \g} +   D_{t|\bar \g}^{\bar \d} f_{t|\bar 0\bar \b\bar \d}
\end{split}
\eeq
for appropriate complex functions $C^{\bar \ve}_{t|\bar 0\bar \b\bar\g}$ and $D_{t|\bar \g}^{\bar \d}$.
Thus, along each integral curve of the $\cC^{k-1,\a}$ vector  field $\frac{d}{dt} + X_t$ of $[0,1] \times \TOB$, the  functions $f_{t|\bar 0\bar \b\bar \g}$ satisfy a system of  linear ordinary differential equations, which  admit a unique solution for each given initial value.  Since  $f_{t|\bar 0\bar \b\bar \g}|_{t = 0} = 0$  and the submanifold $\{0\}\times \TOB$ is  transversal to  
$\frac{d}{dt} + X_t$ at all points,  
the  functions $f_{t|\bar 0\bar \b\bar \g}$    must vanish  along each of the above  integral curves, hence at all points of $[0, 1]\times \TOB$ as claimed. \par
\subsubsection{Step 3 - Existence of a solution to  the  single equation  of  \eqref{5.30}}
 Let us now fix some new notation.  Given a one-parameter  family of   tensors  $\phi_t \in \cH^{01*} \otimes \cH^{10}$ and a system of polar coordinates $(\z, w^\a)$ on $\TOB$, we denote by $\phi_{t|\bar \a}^\b$ the components of $\phi_t$ in the associated adapted  polar frame field $(Z, e_\a, \bar Z, e_{\bar \a})$ and we write
 $\phi_t = \phi_{t|\bar \a}^\b e_\b\otimes e^{\bar \a}$.  We also set $\bE_{t|\bar \a}\= e_{\bar \a} + \phi_t(e_{\bar \a})$ and,   for simplicity,  we use the convenient notation  $<\cdot, \cdot >$ for the  $(0,2)$-tensor field $g^{(\phi_t)}$ introduced in \eqref{innerprod}.
 In this notation,    \eqref{5.30} (i.e., \eqref{5.22bis-1})  becomes
  $$  < \cL_{\frac{d}{dt} + X'_t}  (e_{\bar \a} + \phi_{t|\bar \a}^\b e_\b),  e_{\bar \g} + \phi_{t|\bar \g}^\d e_\d> = < \cL_{\frac{d}{dt} + X'_t}  (e_{\bar \a} + \phi_{t|\bar \a}^\b e_\b),  \bE_{t|\bar \g} >  = 0\ .$$
 Recalling that $< e_\b, \bE_{t|\bar \g}> = g_{\a \bar \g}$ (see  \eqref{convenient}) this  can be also written as
  \beq \label{5.37}
  \left(\frac{d}{dt}  + X'_t\right) (\phi_{t|\bar \a}^\b) g_{\b \bar \g}     =    - <[X'_t, e_{\bar \a}], \bE_{t|\bar \g}>   -  \phi_{t|\bar \a}^\b <[X'_t, e_{\b}], \bE_{t|\bar \g}> \ .
 \eeq
Multiplying both sides by the inverse matrix $(g^{\d \bar \g})$ of $(g_{\a \bar \b})$, we get that 
\beq \label{5.40}  \left(\frac{d}{dt}  + X'_t\right) \phi_{t|\bar \a}^\b = F^\a_{\bar \b}\left( t, x,  X'_t, \phi^\d_{t|\bar \g}, \overline{\phi^\d_{t|\bar \g}}\right)\ ,\eeq
for some appropriate real analytic functions $F^\a_{\bar \b}$ of   the points $(t,x) \in [0,1] \times \TOB$,  of the components of $X'_t$ and of the components of  $\phi_t$. \par
\smallskip
So, by the previous sections, our proof is  reduced to show the following: 
\begin{itemize}[itemsep=6pt plus 5pt minus 2pt, leftmargin=18pt]
\item[]{\it   there exists a maximum value   $s_o \in (0,1]$ such that for each $s \in (0, s_o)$ there is   a one-parameter family $\phi_t$ 
of  tensor fields in $\cH^{01*} \otimes \cH^{10}$, which is  solution to \eqref{5.40} for $t \in [0, s]$ with initial condition $\phi_{t = 0} = \phi_J$ and satisfies  the  nondegeneracy condition} (equivalent to \eqref{mannaggia}) 
\beq \label{mannaggiabis} \det\left( \d^{\a}_{\b}-  \phi^\a_{t|\bar \g} \overline{\phi^\g_{t|\bar \b}}\right) \neq 0\ ; \eeq 
{\it the value $s_o$ is  less then $1$ only in case such solutions are such that  $\lim_{t \to s_o} \det\left( \d^{\a}_{\b}-  \phi^\a_{t|\bar \g} \overline{\phi^\g_{t|\bar \b}}\right) =0$.}
\end{itemize}
\par
\medskip
In the case of  real analytic data, the   proof  is reached in a  direct way. Indeed, \par
\begin{prop} \label{prop5.4} Assume that  $J$ and $v_t$ are real analytic and, consequently,  that also    $\rho$ and   $ \phi_J$ are  real analytic. Then there exists a maximum value   $0 < s_o \leq 1$ such that  the system of partial differential equations    \eqref{5.40}  with initial condition $\phi_{t = 0} = \phi_J$ has a unique solution for all $t \in [0, s_o)$. The case $s_o < 1$ occurs only if 
  \beq \label{condicio}
  \lim_{t \to s_o} \det\left( \d^{\a}_{\b}-  \phi^\a_{t|\bar \g} \overline{\phi^\g_{t|\bar \b}}\right) =0
  \eeq
\end{prop}
\begin{pf}   By  Sublemma \ref{sublemma}, there 
exists a unique local solution to \eqref{5.40} with initial condition  $\phi_{t = 0} = \phi_J$ on a compact set of the form $[0, \ve] \times \TOB$.  Since $\phi_J$ is a deformation tensor of a complex structure, then it  
satisfies \eqref{mannaggiabis}. Thus,  there exists a sufficiently small $\ve$, such that the obtained solution satisfies also \eqref{mannaggiabis}  on the whole $[0, \ve] \times \TOB$.
 If we consider the subset  $B \subset (0,1]$
 \begin{multline} \label{5.51} B = \big\{\ s \in (0,1]\ :\ \text{there is a solution on} \ [0, s] \times \TOB\ \text{satisfying} \\
 \text{the initial condition}\ \phi_{t = 0} = \phi_J\ \text{and}\ \eqref{mannaggiabis}\ \big\}\ ,
 \end{multline}
then  $s_o = \sup B$ satisfies all requirements. By previous remark, $s_o > 0$. On the other hand, if $s_o < 1$ and the solutions on the intervals $[0, s]$
are such that $\lim_{t \to s_o} \det\left( \d^{\a}_{\b}-  \phi^\a_{t|\bar \g} \overline{\phi^\g_{t|\bar \b}}\right) \neq 0$, the same argument shows that 
such solution extends to a solution defined on a 
larger interval $[0, s_o + \ve']$, contradicting the definition of $s_o$. This implies the last claim of the proposition.
\end{pf}

\par
\smallskip
We  now  need to  prove  the existence of the maximal value   $0 < s_o \leq 1$ also in case we need to determine   $\cC^{k-1,\a}$ solutions on  intervals $[0, s] \subset [0, s_o)$, 
corresponding to   (not  real analytic) $\cC^k$  initial datum $\phi_J$. By considering the set $B$ defined in  \eqref{5.51}, it suffices to show that $B \neq \emptyset$ and then setting $s_o = \sup B$. In other words, we only need to show the existence of solutions with $\cC^{k}$  initial datum $\phi_J$ and satisfying \eqref{mannaggiabis} over some compact set of the form $[0, \ve] \times \TOB$ for some sufficiently small $\ve > 0$.  The method of construction of the required  solution will also show that,  in case $s_o < 1$, then necessarily $\lim_{t \to s_o} \det\left( \d^{\a}_{\b}-  \phi^\a_{t|\bar \g} \overline{\phi^\g_{t|\bar \b}}\right) =0$:  otherwise,  one might  determine  solutions on some larger interval $[0, s_o + \ve']$, in contrast with  the definition of $s_o$. \par
\smallskip
The strategy consists in  determining  the desired  solution $\phi_t$ as limit  of real analytic solutions $\phi_{(n)t}$  of \eqref{5.40}, whose  initial values  $\psi_{(n)} \= \phi_{(n) t= 0}$  are real analytic and approximate the  initial value $\phi_J$. 
For simplicity, we first work  under the additional assumption that  $\r$ and $v_t$ are  real analytic.
Consider  a finite covering of $\TOB$ by compact sets $K_i$, $i = 1, \ldots, N$, with smooth boundaries, on each of which the initial data $\phi_{J\bar \b}^\a\big|_{K_i}$ are  approximated in $\cC^{k}$-norm (thus, also in $\cC^{k-1, \a}$-norm for each $\a \in (0,1)$)  by a sequence  $\psi_{(n)\bar \b}^\a$ of real analytic functions. Such  approximating sequence $\psi_{(n)\bar \b}^\a$ surely exists and it can be even assumed to be formed by  rational polynomials  (see e.g. \cite{Fe}, Thm. 6.10). \par
\smallskip
Now,  given a  compact set $K =$ $K_i$ and an integer $n$,  let $\phi_{(n)t|\bar \b}^\a$   be  the unique solution   to \eqref{5.40} and \eqref{mannaggiabis} with initial value $\phi_{(n)t = 0|\bar \b}^\a  = \psi_{(n)\bar \b}^\a$ on $\{0\} \times K$,  defined on some compact set of the form $[0,t_n] \times K$,  with $t_n$ less than the maximal value $\l^{(n)}_o$ determined by Proposition \ref{prop5.4}.
We may  assume that $t_n$ is the infimum of the values $t'$, for which $\displaystyle \sup_{[0,t'] \times K} \|\phi_{(n) t|\bar \b}^\a\| \leq 1$.  By the method of construction of solutions using complex characteristics described in the proof of Sublemma \ref{sublemma}, 
such infimum $t_n$ is  bounded from below by some $t_o > 0$, independent of $n$ and determined by \\[-15pt]
\begin{itemize}[leftmargin=16pt]
\item[--] the value $ \sup_{ K} \|\psi_{(n)|\bar \b}^\a\|$, which we may assume  to be  less than or equal to 
$\sup_K  \|\psi_{\bar \b}^\a\| + \ve_o < 1$ for some $\ve_o > 0$  for all sufficiently large $n$; 
\item[--] the sup of the components of the holomorphic vector field $\wh X$, determined through \eqref{realanalytic-1} by the data in \eqref{5.40}, taken over the compact set of $\bC^{n+1} \times \bC^N$, which is the cartesian product $[0,1] \times K$ and the hypercube in $\bC^N$ spanned by  all possible values for  $\phi_{\bar \b}^\a$  with  norm bounded above by $2$. 
\end{itemize}
We may therefore assume that all solutions $\phi_{(n) \bar \b}^\a$ are defined on a compact set  $[0,t_o] \times K$, with  $t_o > 0$, 
and that they are  equibounded by $\displaystyle \sup_{[0,t_o] \times K} \|\phi_{(n) t|\bar \b}^\a\| \leq 1$. \par
\smallskip
We  now want to show that the solutions $\phi_{(n)t| \bar \b}^\a$ have also equibounded  first derivatives, so that they are equicontinuous on $[0,t_o] \times K$. To check this, for all complex vector field $e_{\a}$ and $e_{\bar \a}$ in $\cH^\bC|_{K}$, consider the unique real analytic extensions $\wh e_{t|\a}$, $\wh e_{t|\bar \a}$ in $\cH^\bC|_{[0,1] \times K}$   satisfying the differential problem 
\beq \label{starina} \cL_{\frac{d}{dt} + X'} \wh e_{t|\a} = \cL_{\frac{d}{dt} + X'} \wh e_{t|\bar \a} = 0\ ,\qquad   \wh e_{t = 0|\a}  = e_\a\ ,\   \wh e_{t = 0|\bar \a} = e_{\bar \a}\ .\eeq
Since these differential equations have  the form \eqref{realanalytic}, with right hand side linear in the unknowns $\wh e_{t|\a}$,  $\wh e_{t|\bar \a}$,   such extensions exist  on some set of the form $[0,\ve] \times K$, $\ve > 0$ (Sublemma \ref{sublemma}). They have  also $\cC^0$-norms that are uniformly bounded from above by the $\cC^0$-norm of the  initial values multiplied by some fixed constant. This  is indeed a consequence of the fact that 
the differential equations along the complex characteristics are linear. It actually  implies that the  $\wh e_{t|\a}$, $\wh e_{t|\bar \a}$ can be  extended over the whole $[0,1] \times K$. \par
\smallskip
We may now observe that, differentiating both sides of \eqref{5.37} along the vector fields $ \wh e_{t|\a}$ and  $\wh e_{t|\bar \a}$ and using the property that they  commute with   $\frac{d}{dt} + X'$ by \eqref{starina}, the derivatives $\wh e_{t|\a}(\phi_{(n)t| \bar \b}^\g)$ and  $\wh e_{t|\bar \a}(\phi_{(n)t| \bar \b}^\g)$ are solutions of a system of the form \eqref{realanalytic}, in which the right  hand side  is linear in the unknowns. Using once again Sublemma \ref{sublemma} and the linearity of the equations along the complex characteristic, we conclude that such derivatives are well defined at all points of $[0,t_o] \times K$, with $\cC^0$-norm  bounded above 
by the $\cC^0$-norms of their initial values times some fixed constant. Since  the initial values $\psi_{(n)\bar \b}^\g$ converge in $\cC^{k-1,\a}$-norm to $\phi_{J\bar \b}^\g$, with $k \geq 2$, we conclude that also the family of derivatives $\wh e_{t|\a}(\phi_{(n) t|\bar \b}^\g)$ and  $\wh e_{t|\bar \a}(\phi_{(n) t|\bar \b}^\g)$ are  equibounded. 
From this,  the equicontinuity   of the solutions $\phi_{(n) t|\bar \b}^\g$ follows, as claimed.\par
\smallskip
As direct consequence of Ascoli-Arzel\`a Theorem,  
the sequence of solutions $\phi_{(n) t|\bar \b}^\g$  converges in $\cC^0$-norm
to the components of  the family of deformation tensors $\phi_t$,  with initial value $\phi_{t = 0} = \phi_J$. An argument that  is basically the same as  the one used for proving equiboundedness of $\cC^1$-norms  shows that also all 
  derivatives up to order $k-1$ of the solutions $\phi_{(n) t|\bar \b}^\g$ are equibounded, proving in this way that the limit $\phi_t$ is  at least of  class $\cC^{k-1}$. In particular,  it follows that the limit is a solution to the first order differential problem \eqref{5.40}.\par
  \smallskip 
   However we claim that also the H\"older ratios of power $\a$ of the $(k-1)$-th order derivatives of the solutions $\phi_{(n) t|\bar \b}^\g$ are equibounded. This can be checked as follows. First, observe that the H\"older ratios of the $(k-1)$-th derivatives, evaluated for pairs of points that are in the same complex characteristic, are equibounded because  the restrictions to the complex characteristic of such  $(k-1)$-th derivatives are solutions of  the above described linear  system of ordinary differential equations with equibounded coefficients. This means that, in order to check the equiboundedness of  the H\"older ratios, we may restrict to considering those 
   evaluated at pair of points $(t, \wt x)$, $(t', \wt y) \in [0,1]\times K$,  with the same coordinate $t$. Now,  for each given pair of points $x, y \in K$, consider the  H\"older ratio of the $(k-1)$-th derivatives, evaluated per each $t$ at the  pair of points $(t, x_t, \phi_{t|\bar \b}^\g(x_t)), (t, y_t, \phi_{t|\bar \b}^\g(y_t))$, belonging  to the  two    complex characteristics originating  from  
 $(t = 0, x, \phi_{0|\bar \b}^\g(x))$ and $(t = 0, y, \phi_{0|\bar \b}^\g(y))$, respectively.  Up to multiplication by a positive constant, such H\"older ratios  are equibounded by the H\"older ratios at their uniquely associated pair of points  at $t = 0$, from which  the complex characteristics originate. Indeed, this is obtained by considering such H\"older ratios as functions of the  (complex)  variable $t$ of the pair of characteristics.  Using the differential equation satisfied by the $(k-1)$-th order derivatives  of $\phi_t$ in order to express  the first derivatives in $t$ of  the  H\"older ratios, one can observe that  they satisfy a simple linear system of  ordinary differential equations, with  coefficients,  given by rational functions  in which a) the  denominators are  bounded away from $0$ by constants that are independent of  $\phi_t$, b) 
  the  numerators depend on the functions $\phi_{(n)t|\b}^\g$ and their $(k-1)$-th derivatives,  and are equibounded. This implies the claimed existence of common upper bounds determined by  the H\"older ratio of the corresponding pair of points at $t = 0$.   \par
 \smallskip
 Since the initial values converge to $\phi_{J \bar \b}^\g$ in $\cC^{k-1,\a}$-norm, we derive in this way the  desired equiboundedness of the $\cC^{k-1,\a}$-norms of the  $\phi_{(n) t|\bar \b}^\g$
 and  the property that the limit deformation tensor  $\phi_t$ is actually of class $\cC^{k-1,\a}$. Note that, by the results of previous section, such  $\phi_t$ necessarily satisfies also the integrability conditions and  the condition (B) of the  deformation tensors of the L-complex structures.\par
 \smallskip  
By uniqueness of the  limits of the sequences of solutions $\phi_{(n) t|\bar \b}^\g$ on the intersections $K_i \cap K_j$ of two  compact sets of the considered covering $\{K_i\}$, we also obtain that all solutions constructed in this way on the sets  $[0,t_{o i}] \times K_i$ combine together 
and give a global solution over a compact set of the form $[0,\ve] \times \TOB$, for some appropriate $\ve > 0$,   under the considered   regularity assumptions. \par
\medskip
It remains to  prove the existence of $\cC^{k-1,\a}$ solutions to \eqref{5.40} and \eqref{mannaggiabis} under the assumption that also the functions $v_t$ and  $\r$ are not real analytic, but just of class  $\cC^{k}$. The method is precisely the same as before with the only difference that, now,  we have to approximate  with real analytic functions on compact sets $K_i$ not only  the initial datum $\phi_{J\bar \b}^\g|_{K_i}$  but also  the $\cC^{k}$ functions $v_t$ and  $\rho$.  A diagonal argument  implies  that also in this case the $\cC^{k-1,\a}$-norms of the  solutions $\phi_{(n)t|\bar \b}^\g$,  determined by the initial values $\psi_{(n)\bar \b}^\g$ and  real analytic vector fields $X'_{(n)t}$, determined by  the $v^{(n)}_t$ and  $\rho_n$,   are  equibounded. This implies the needed existence of a $\cC^{k-1,\a}$-solution also in this case. \end{pf}  
\subsubsection{The $\cC^{\infty}$ and the $\cC^{\omega}$ cases and Theorem \ref{main2}}
We end by observing that our arguments show that Theorem \ref{main} may be rephrased both in the  $\cC^{\infty}$  and the $\cC^{\omega}$ versions so that  also Theorem \ref{main2} follows. Indeed, as we already pointed out, the determination of the maximal value $\l_{v}$ in Theorem \ref{main} depends on the fulfillment of the conditions (A) -- (D) of Subsection 5.2.1,  which, as long as the data are at least of class  $\cC^{2}$, is independent on the degree of smoothness.  Because of the uniqueness of the Green pluripotential with a given pole, the $\cC^{\infty}$ version of Theorem \ref{main} is obtained applying the statement Theorem \ref{main}  for all $k\geq5$. 
Furthermore the steps of proof of the crucial Lemma \ref{lemma52} are provided first in the $\cC^{\omega}$ case (see Sublemma \ref{sublemma} and Proposition \ref{prop5.4}) and therefore also the   $\cC^{\omega}$ version of Theorem \ref{main} holds. 
\subsubsection{A Final Remark}
As we already mentioned,  the consideration  of  manifolds of circular type  was motivated by the relevance of two main classes of examples: the  smoothly bounded, strictly  convex domains  and the smoothly bounded strongly  circular domains of $\bC^n$.   However  the  manifolds of circular type 
are described in  abstract terms, an aspect    that   also allows  explicit construction procedures. 
In fact,  as we know,  each manifold of circular type is biholomorphic to a  manifold in normal form. This is nothing but the unit ball  $\bB^n$  equipped  with a non-standard complex structure  $J$, obtained by deforming   the standard one along special directions (Definition \ref{definition5.1}).   By \cite{BD,PS,PS2}, all such non-standard  complex structures $J$   are uniquely determined by their deformation tensors which, in turn, are parameterised in an explicit (although  non-trivial)  way  
 by a single freely specifiable complex function $h$.  By  choosing    appropriately such a function, one can construct   manifolds of circular type $(\bB^n, J)$ with  desired  properties. For instance, this idea has been successfully used in \cite{PS2} to prove the existence of  abstract   manifolds of circular type  with prescribed regularity properties at the center.   
The same technique  can be adopted to generate many  domains of circular type $(\overline \bB^n, J)$   for which there are  points that satisfy the (closed) condition (\ref{condicio}), provided that such condition is   translated, as we did in \cite{PS2},   into manageable differential equations for the deformation tensor  of  the complex structure. By our proof of Theorem \ref{main0},  the existence of such points  would imply that   the ``cloud'' of centers $\cU^{\text{(max)}}$ is a proper  subset of  the constructed domain. Further,  as in \cite{PS2}, such abstract examples would be embeddable as domains of $\bC^n$ as long as their complex structures are suitably  small deformations of the standard euclidean complex structure. 
Conversely, we also expect that a careful  study of the mentioned condition (\ref{condicio}) and of its  translation as a condition on the deformation tensor  would  suggest new biholomorphically invariant sufficient conditions  for the   regularity propagation  to occur globally. We  plan to  address these lines of research in the future.
\par
\begin{acknowledgements}
 This research was partially supported by the Project MIUR ``Real and Complex Manifolds: Geometry, Topology and  Harmonic Analysis'' and by GNSAGA of INdAM
\end{acknowledgements}



\end{document}